\documentclass[3p]{elsarticle}
\usepackage{amsmath}
	\numberwithin{equation}{section}

\usepackage[british]{babel}
\usepackage{tabularx}
\usepackage{multirow}
\usepackage{booktabs}
\usepackage{longtable}
\usepackage{rotating}
\usepackage{lineno}
\usepackage{pdflscape}
\usepackage[table]{xcolor}
\usepackage{subcaption}
\usepackage[export]{adjustbox}
\usepackage{amsmath} 
\usepackage{mathtools}
\usepackage{rotating}

\usepackage{algorithm}
\usepackage{algpseudocode}

\usepackage{letltxmacro}

\usepackage{amssymb}
\usepackage{amsmath}

\usepackage{amssymb}
\usepackage{amsthm}
\usepackage{mathtools}
\usepackage{braket}
\usepackage{bm}
\usepackage{bbm}
\usepackage{graphicx}
\usepackage{xcolor}
\usepackage{bbold}
\usepackage{blindtext}
\usepackage{graphicx}
\usepackage{graphics}
\usepackage{verbatim}   % for math
\usepackage{amsfonts}
\usepackage{adjustbox}
\usepackage{mathrsfs}  % % % % % % % % % % % % % % % % % % % % %

\usepackage{letltxmacro}

\usepackage{placeins}

\usepackage{siunitx}

\newcommand{\tens}[1]{\bm{\underline{#1}}}
\newcommand{\tend}[1]{\hbox{\oalign{$\bm{#1}$\crcr\hidewidth$\scriptscriptstyle\bm{\sim}$\hidewidth}}}
\newcommand{\tent}[1]{\hbox{\oalign{$\bm{#1}$\crcr\hidewidth$\scriptscriptstyle\bm{\simeq}$\hidewidth}}}
\newcommand{\tenq}[1]{\hbox{\oalign{$\bm{#1}$\crcr\hidewidth$\scriptscriptstyle\bm{\approx}$\hidewidth}}}

\newcommand{\R}{\bm{\underline{r}}}
\newcommand{\Q}{\bm{\underline{q}}}

\newcommand{\dep}{\bm{\underline{u}}}

\newcommand{\eps}{\hbox{\oalign{$\bm{\epsilon}$\crcr\hidewidth$\scriptscriptstyle\bm{\sim}$\hidewidth}}}
\newcommand{\deps}{\hbox{\oalign{$\bm{\delta\epsilon}$\crcr\hidewidth$\scriptscriptstyle\bm{\sim}$\hidewidth}}}
\newcommand{\pol}{\hbox{\oalign{$\bm{\tau}$\crcr\hidewidth$\scriptscriptstyle\bm{\sim}$\hidewidth}}}
\newcommand{\stress}{\hbox{\oalign{$\bm{\sigma}$\crcr\hidewidth$\scriptscriptstyle\bm{\sim}$\hidewidth}}}

\newcommand{\stiff}{\hbox{\oalign{$\bm{\lambda}$\crcr\hidewidth$\scriptscriptstyle\bm{\approx}$\hidewidth}}}

\usepackage[T1]{fontenc}
\usepackage[utf8]{inputenc}

\makeatletter
\newcommand\thefontsize{The current font size is: \f@size pt}
\makeatother

\usepackage{lineno,hyperref}
\modulolinenumbers[5]

\journal{ }

\begin{document}
%\layout
%\setuplayouts

\begin{frontmatter}

%\title{A study on variant selection during martensitic phase transition by atomistic simulations} 

%\title{\textcolor{black}{A new discrete FFT-based solver for micromechanics}}
\title{\textcolor{black}{A fast and robust discrete FFT-based solver for computational homogenization}}



\author[onera]{A. Finel \fnref{myfootnote}}

%\cortext[mycorrespondingauthor]{Corresponding author: alphonse.finel@onera.fr}

\fntext[myfootnote]{alphonse.finel@onera.fr}

\address[onera]{Université Paris-Saclay, ONERA, CNRS, Laboratoire d'Etude des Microsctructures (LEM), 92322, Ch\^{a}tillon, France}

%\

\begin{abstract}

We propose a new discrete FFT-based method for computational homogenization of micromechanics on a regular grid that is simple, fast and robust. The discretization scheme is based on a tetrahedral stencil that displays three crucial properties. First, and most importantly, the Fourier representation of the associated Green operator is defined for any finite q-vector generated by the periodic boundary conditions and that does not belong to the Reciprocal Lattice of the discrete grids. As shown in the paper, this property guaranties that, for any elastic contrats, even infinite, mechanical equilibrium is always mathematically stable, i.e. free of any unphysical patterns, such as oscillations, ringing or checkerboarding, a property which is not shared by the original Moulinec-Suquet method \cite{moulinec1994fast,moulinec1998numerical} nor by the rotated scheme proposed by Willot \cite{willot2015fourier}. Second, the components of tensorial quantities are all defined on the same location, which permits the use of any elastic anisotropy and any spatial variation of the material fields. Third, convergence to equilibrium using the simplest iterative scheme, the "basic scheme", is fast and the number of iterates stabilizes at high contrasts, so that infinite contrast is obtained  without additional computational cost.

\end{abstract}


\end{frontmatter}

\nolinenumbers

%%%%%%%%%%%%%%%%%%%%%%%%%%%%%%%%%%%%%%%%%%%%%%%%%%%%%%%%%%%%%%%%%%%%%%%%%%%%
\section{Introduction}\label{sec:intro}
%%%%%%%%%%%%%%%%%%%%%%%%%%%%%%%%%%%%%%%%%%%%%%%%%%%%%%%%%%%%%%%%%%%%%%%%%%%%
%

Computational micromechanics and homogeneization require the solution of the mechanical equilibrium of complex microstructures. Closed-form estimates of average properties, often based on the Eshelby inclusion problem, are of little use because real microstructures depart a lot from a single inclusion embedded in an infinite matrix. The alternative is to resort to full-field simulations of a complete and realistic microstructure. Among the available numerical techniques, Fourier-based methods have, during the last decades, gained an increasing popularity for several reasons. A full mesh discretization of the microstructure, which is a very complex operation in highly heterogeneous systems, is not required. The calculation is simply based on a regular grid which, in turn, may be directly identified with the voxels of an image obtained by 3D imaging techniques. Also, the underlying regular grid permits the use of very efficient librairies that implement Fast Fourier Transform (FFT) algorithms.  Finally, the method is easy to implement and produces impressive results for large-scale microstructures in a very short time.

Fourier-based methods have been introduced by Moulinec-Suquet \cite{moulinec1994fast,moulinec1998numerical} in the context of computational homogenisation (see \cite{schneider2021review} and \cite{lucarini2021fft} for recent reviews). They have also long been popular in the field of structural transformations since the pioneering work of Khachaturyan \cite{khachaturyan1983} on the elastic theory of solids and its extension to elastically inhomogeneous systems \cite{khachaturyan1995elastic} (see an earlier review in \cite{chen2002phase}). 

In fact, the original Moulinec-Suquet \cite{moulinec1994fast,moulinec1998numerical} and Khachaturyan \cite{khachaturyan1995elastic} methods for dealing with elastically inhomogeneous solids are based on the same foundations, consisting of a continuous formulation of mechanical equilibrium in a periodic unit cell, which is then solved in Fourier space using a simple iterative algorithm. Indeed, as noted for example in \cite{eyre1999fast}, the so-called  "basic scheme" introduced by Moulinec-Suquet may be interpreted as a Neumann expansion of the inhomogeneous elastic Green function associated to the microstructure. This is precisely the starting point of the perturbation method introduced by Khachaturyan, which can be recast as an iterative algorithm, as shown in \cite{hu2001phase}, which in turn happens to be equivalent to the basic scheme.

It has rapidly been observed that the original method, in its most basic form, suffers from two difficulties. First, the basic scheme does not converge for infinite contrast in the stiffness tensors, which has been associated to the non-convergence of the underlying Neumann series. Second, oscillations are often observed in the computed stress and strain fields in the presence of sharp material discontinuities. This is due to the numerical implementation of the method. Indeed, the mechanical equilibrium of the periodic unit cell and its Lippmann-Schwinger formulation are initially derived in a continuous setting. However, for the numerical implementation, a discretization of the real space must be considered. This discretization is associated, in the Fourier space, to a high frequency cut-off that generates oscillations and checkerboard patterns in real space.  

Accelerated algorithms have been proposed to overcome the convergence problems encountered with the basic scheme, such as the method proposed by Eyre and Milton \cite{eyre1999fast} in the context of electrical conductivity and further extended to elasticity \cite{michel2001computational,vinogradov2008accelerated}, the augmented Lagrangian method proposed by Michel et al. \cite{michel2000computational,michel2001computational} and the polarization method proposed by Monchiet and Bonnet \cite{monchiet2012polarization}. As shown in \cite{moulinec2003comparison} and \cite{moulinec2014comparison}, these methods converge considerably faster than the basic scheme. However, the simplicity of the original method is lost and, as shown in \cite{moulinec2014comparison},  still none of them converge when the contrast between the phases is infinite. 

Another FFT methodology, based on the Hashin-Shtrikman variational principle, has been proposed by Brisard and Dormieux \cite{brisard2010fft,brisard2012combining}. The main ingredient of this method is to consider piecewise constant approximations of the polarization field and ask the corresponding discrete variables to fulfil an energy principle instead of solving a discrete Lippmann-Schwinger equation, as it is done within the original Moulinec-Suquet method. As the associated algorithm iterates on successive approximations of the polarization field, the method is also refered to as a polarization-based method.

In it spirit, this polarization-based method is similar to the method introduced earlier by Rodney et al.  \cite{rodney2003phase} for modelling dislocations within a phase field approach. In this latter method, dislocation loops are decomposed into elementary loops, called \textit{loopons}, within which the associated phase field is piecewise constant. The elastic interactions between the loopons is then given by \textit{decorated} interactions matrices whose Fourier transforms are given by infinite sums over the discrete sites of the Reciprocal Lattice in Fourier space. The resulting Green operator is exactly the so-called \textit{periodized} Green operator of the polarization-based method proposed by Brisard and Dormieux \cite{brisard2010fft,brisard2012combining}. Thanks to its (piecewise) continuous character, the polarization-based method is automatically free of oscillatory artefacts and can handle infinite stiffness contrasts. However, as noted in \cite{brisard2012combining}, the infinite series involved in the calculation of the periodized Green operator converge very slowly, making their numerical implementation rather involved, especially for 3D cases. 

More recently, Eloh et al. \cite{eloh2019development} proposed a method which is also based on a piecewise constant approximation of local fields, but applied at the level of the Lippmann-Scnwinger equation. The method therefore also leads to a redefinition of the Green operator in the form of infinite series in reciprocal space and, consequently, suffers from the same weakness as the polarization-based method of Brisard and Dormieux.

All the FFT-based methods mentioned above have in common the formulation of the initial mechanical problem in a continuous setting, followed by a specific discretization strategy to generate a discrete Green operator. An alternative to this route is to readily start with a finite difference discretization of the relation between the strain and displacement fields, both defined on specific discrete grids. Besides the centered discretization scheme, which is known to generate strong oscillations, various discretization methods have been proposed. 

The standard staggered grid scheme, initially developed for fluid mechanics problems \cite{harlow1965numerical}, has been used in a Fourier-based setting by Geslin and Finel \cite{geslin2014investigation} and Schneider et al \cite{schneider2016computational}. Its main advantage is that it is free from any spurious oscillations. However, the price to pay is that different components of the displacement field and also of the strain field are stored on different locations, which generates difficulties in presence of sharp interfaces and in situations in which the constitutive law links together tensor components that do not sit on the same location. The discretization scheme proposed by Willot \cite{willot2015fourier}, based on a rotated staggered grid initially introduced for simulating the propagation of elastic waves \cite{saenger2000modeling}, avoids these difficulties. As it presents several positive properties, in particular the possibility to simulate porous materials and the reduction of the oscillations observed in the Moulinec-Suquet discretization, the rotated scheme has become the standard discretisation method for FFT-based numerical homogenization. However, as will be shown below, it is not exempt of spurious oscillations or chechkerboard patterning effects. Finally, the discretization scheme proposed by Ruffini and Finel \cite{ruffini2015phase} avoids the weaknesses of the two previous schemes: the components of a given tensor are all stored on the same location and it is free from any non-physical instability. However, this comes at the expense of introducing a subgrid discretization in which each initial voxel is subdivided into $2^D$ subvoxels, where $D$ is the spatial dimension, which increases in the same proportion the computational load.

In the present work, we propose a new discretization method, based on a tetrahedral stencil between the strain and displacement fields. This stencil avoids the difficulties encountered with the discrete schemes discussed above. It is mathematically free from any non-physical instabilities, does not resort on subvoxels defined on a finer scale and the components of a given tensor are all defined on the same location. Finally, we show through numerical applications that convergence to mechanical equilibrium using the initial Moulinec-Suquet basic scheme is fast, even in the presence of voids.

The paper is organized as follows. In section \ref{sec:method}, we develop the model and its numerical implementation. In section \ref{sec:numerical_exp}, we analyse the performance of the tetrahedral stencil method for different applications, including comparisons with the original Moulinec-Suquet method and the rotated finite difference scheme and we conclude in section \ref{sec:conclusion}.

The following notation is used throughout the paper. $a$ refers to a scalar quantity whereas bold symbols $\tens{a}$, $\tend{a}$, $\tent{a}$ and $\tenq{a}$ denote first-, second-, third- and fourth-order tensor quantities, respectively.  $\overline{A}$ is the complex conjugate of $A$, $f(\Q)$ represents the Fourier transform of $f(\R)$ and, as usual, single dot ($.$) and double-dot ($:$) represent contraction and double contraction, respectively.

%%%%%%%%%%%%%%%%%%%%%%%%%%%%%%%%%%%%%%%%%%%%%%%%%%%%%%%%%%%%%%%%%%%%%%%%%%%%
\section{Methods : a tetrahedral stencil}\label{sec:method}
%%%%%%%%%%%%%%%%%%%%%%%%%%%%%%%%%%%%%%%%%%%%%%%%%%%%%%%%%%%%%%%%%%%%%%%%%%%%

The problem at hand is to solve for static mechanical equilibrium in a 3D cell containing a representative volume of a microstructure composed of one or more phases embedded into a matrix and subjected to an applied overall strain or external stress. We use linear elasticity at small strains in its standard form as our point of departure for discussing the discretization schemes and the solution methods. However, the discretization method proposed below is general and may be used in more complex situations, such as nonlinear elasticity or plastic flow. The microstructure is described by a stiffness distribution $\stiff(\R)$ and by an eigenstrain distribution $\eps^0(\R)$ that measures a free-stress deformation field with respect to the matrix. The local stress $\stress(\R)$ is therefore related to the local strain $\eps(\R)$ by
\begin{equation}
 \label{eq:linear_elasticity}
 \stress(\R) = \stiff(\R):(\eps(\R) -\eps^0(\R)).
\end{equation}
The central point of our treatment is a proper choice of a discretization scheme that relates the strain and displacement fields. We proceed as follows.
%
%%%%%%%%%%%%%%%%%%%%%%%%%%%%%%%%%%%%%%%%%%%%%%%%%%%%%%%%%%%%%%
\subsection{Discretization scheme}\label{sec:discretization}
\begin{figure}[h]
	\centering
	\includegraphics[width=0.79\textwidth]{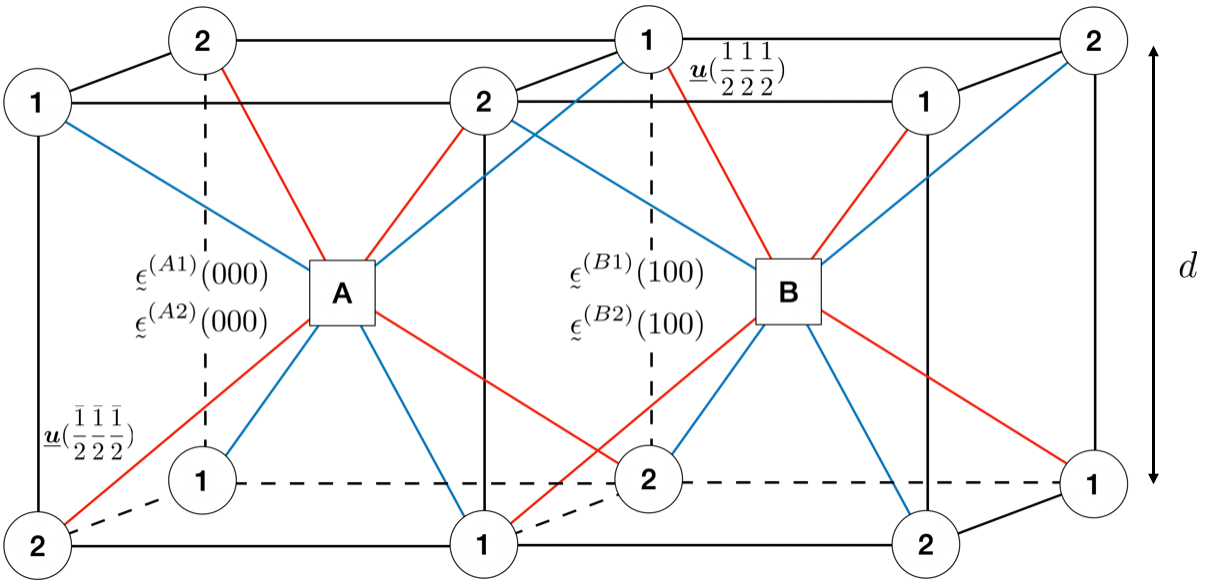}
	\caption{ \textcolor{blue}{Tetrahedral stencil for a simple cubic discretization. Only two voxels along axis (100) are shown. The simple cubic displacement grid, whose sites are located at the corners of the voxels, is split into two FCC subgrids whose sites are labelled by the indexes 1 and 2, respectively. Similarly, the simple cubic strain grid, whose sites are located at the centers of the voxels, is split into two FCC subgrids whose sites are labelled by the indexes A and B, respectively. The four translations that relate a strain site of type A to the four sites of a tetrahedron of type 1 (blue lines that emerge from a site A) are related through a central symmetry operation to the four translations that relate a strain site of type B to the four sites of a tetrahedron of the same type 1 (red lines that emerge from a site B). This central symmetry relation is reflected in the corresponding finite difference operators $\tens{D}^{(A1)}$ and $\tens{D}^{(B1)}$ that link the strains on subgrids A and B to the displacements on subgrid 1, as shown in Eq.~(\ref{eq:DB1}). Similarly, the translation sets that relate strain sites of type A and B to tetrahedra of type 2 (red and blue lines that emerge from sites of type A and B, respectively) are related by the central symmetry operation mentioned above. This is reflected into the relation between the finite difference operators $\tens{D}^{(A2)}$ and $\tens{D}^{(B2)}$ given in Eq.~(\ref{eq:DA2}). Finally, we observe that the four translations that relate a strain site A to displacement sites of type 1 are identical to the ones that relate a strain site B to displacement sites of type 2, which implies the equality of the operators $\tens{D}^{(A1)}$ and $\tens{D}^{(B2)}$, as stated in Eq.~(\ref{eq:DB2}). }}
\label{fig:tetrahedral_stencil}
\end{figure}

%%%%%%%%%%%%%%%%%%%%%%%%%%%%%%%%%%%%%%%%%%%%%%%%%%%%%%%%%%%%%%
%
We consider a 3D periodic rectangular cell discretized into $N_1\times N_2\times N_3=N$ cubic voxels. The integer $N_i$, $i=1,3$, are even. The edges of the cell are parallel to the cartesian coordinates and the voxel edge length is equal to $d$ (see Fig.~(\ref{fig:tetrahedral_stencil})). We refer to voxel $(l_1,l_2,l_3)$ as the voxel whose center sits at position $\R$ given (in units of $d$) by 
\begin{equation}
\label{eq:r-vectors}
\begin{split}
\R &= l_1\,\tens{a}_1 + l_2\,\tens{a}_2 +  l_3\,\tens{a}_3, \\
 l_i &\in[0,...,N_i-1]    \quad ( i=1,2,3),
\end{split}
\end{equation}
where $\tens{a_i}, i=1,2,3,$ are unit vectors along the cartesian axis. The centers of the voxels form a simple cubic (SC) lattice with a grid spacing equal to $d$. The eight corners of cell $(l_1,l_2,l_3)$ sit on sites whose coordinates are $(l_1\pm 1/2,l_2\pm 1/2,l_3\pm 1/2)$. Altogether, the corners of all cells form another SC lattice. 

The first step is to choose the locations where the different tensors are defined. In order to handle any elastic anisotropy and any spatial variation of the material fields, the discrete version of the constitutive relation that link stresses and strains must involve tensorial quantities all defined on the same location.  Here, all the 2nd-order and 4th-order tensors, i.e. the strain and stress fields as well as the material stiffness and eigenstrain, will be defined on the SC voxel grid, hereafter referred to as the SC strain grid, whereas the displacement field will live on the SC corner grid, hereafter noted as the SC displacement grid. In particular, all the components (diagonal and non-diagonal) of a tensor live on the same location, which avoids the difficulty that arises when the discretization of a constitutive laws relates a component of a tensor to components of another one defined on different locations, as it is the case in the usual SC staggered scheme mentionned above (see for example \cite{schneider2016computational}). 

The key point now is to select a judicious stencil that relates the strain field to the displacement field. In order to automatically avoid the instabilities often associated to finite difference schemes, we propose a tetrahedral stencil, in which the strain at the center of a voxel is linked to the displacement vectors that sit at the fourth corners of one of the two regular tetrahedra that can be constructed out of the eight corners of that voxel (see Fig.~(\ref{fig:tetrahedral_stencil})). In other words, we split the SC displacement grid into two face centered cubic (FCC) subgrids, that we label with the indexes 1 and 2 used in Fig.~(\ref{fig:tetrahedral_stencil}) to represent the corresponding sites and that we also used as superscripts when we refer to a quantity defined on the corresponding subgrid. Now, considering one of the FCC displacement subgrid, say subgrid 1, we observe that any site of the SC strain grid sits at the center of one of the tetrahedra formed by the nodes of this subgrid. This offer the possibility to define two tetrahedral discretization schemes, associated to displacement subgrids 1 and 2, respectively.

Obviously, the discretization scheme associated to one of the FCC displacement subgrids does not involve the displacements on sites that belongs to the other subgrid. 
As a result, if we select only one of the two FCC discretization schemes, some symmetry elements on the SC strain grid that would be preserved if we had used a full SC displacement grid are now lost. In order to recover the full cubic symmetry on the SC strain grid, we consider simultaneously the two discretization schemes.  The final strain field will be simply defined as the average of the strain fields generated by the two discretization schemes. However, since  this averaging is "on site", it will not generate a strain field less accurate than the two initial fields, as would be observed if some spatial averaging process had been used (see for example \cite{Willot2008FastFT}). Also, although one might naively expect that using two discretisation schemes would double the computational effort in terms of memory and CPU time, we will see below that this is not the case, as the two mechanical equilibria that emerge from the two discretisation schemes can be collapsed into a single one.

Next, considering one of these FCC dicretisation schemes, say discretization scheme 1, we observe that the translations that relates the four vertices of a tetrahedron, which sit on subgrid 1,  to its center, which sits on the SC strain grid, alternate between two different sets, represented in red and blue, respectively,  in Fig.~(\ref{fig:tetrahedral_stencil}). In other words, the SC strain grid must also be split into two FCC subgrids, that we label with the indexes $A$ and $B$ used in Fig.~(\ref{fig:tetrahedral_stencil}) to represent the corresponding sites and that we also use as superscripts to label a quantity that is defined on the corresponding strain subgrid.

We are now in position to define the finite difference schemes that relate the strains defined on the FCC strain subgrids to the displacements defined on the FCC displacement subgrids. We first split the strains into their spatial averages and fluctuating terms, i.e. for discretization schemes $k=1,2$ and for strain subgrids  $K=A,B$, we write:
\begin{linenomath*}
\begin{align}
\text{if } \R \in K ,\quad  \eps^{(Kk)}(\R)= \bar \eps + \deps^{(Kk)}(\R),
\label{eq:eps}
\end{align}
\end{linenomath*}
where the notation $ A^{(Kk)}(\R)$ refers to a quantity that lives on a site $\R$ that belongs to strain subgrid $\rm K$ and that is defined using the discretization scheme associated to displacement subgrid $k$.
For strain subgrids $K=A,B$ and discretization schemes $k=1,2$, the fluctuating strain fields $\deps^{(Kk)}(\R)$, whose spatial averages are equal to zero, are linked to the fluctuating displacement field $\dep^{(k)}(\R)$ by:
\begin{equation}
\text{if } \R \in K ,\quad  \delta \epsilon^{(Kk)}_{ij}(\R) = \frac12 \{D^{(Kk)}_i[u^{(k)}_j](\R)+D^{(Kk)}_j[u^{(k)}_i](\R)\},
\label{eq:deps}
\end{equation}
where the dimensionless quantities $u^{(k)}_i$, $i=1,2,3$, represent the coordinates of the displacement field in units of the cell length $d$. Following the lines of the tetrahedral stencil described above, inspection of Fig.~(\ref{fig:tetrahedral_stencil}) shows that the operators $D^{(Kk)}_i$, $i=1,2,3$, act on a discrete scalar field $f$ according to
\begin{subequations}
\label{eq:DA1}
\begin{align}
D_1^{(A1)}[f](lmn) =& \frac12\{ f(lmn +\frac{111}{2}) + f(lmn + \frac{1\bar1\bar1}{2}) - f(lmn + \frac{\bar1\bar11}{2}) - f(lmn + \frac{\bar11\bar1}{2})\}, 
\label{eq:DA1_1} \\
D_2^{(A1)}[f](lmn) =& \frac12\{ f(lmn + \frac{111}{2}) + f(lmn + \frac{\bar11\bar1}{2}) - f(lmn + \frac{\bar1\bar11}{2}) - f(lmn + \frac{1\bar1\bar1}{2})\},
\label{eq:DA1_2} \\
D_3^{(A1)}[f](lmn) =& \frac12\{ f(lmn +\frac{111}{2}) + f(lmn + \frac{\bar1\bar11}{2})  - f(lmn + \frac{\bar11\bar1}{2}) - f(lmn + \frac{1\bar1\bar1}{2})\},
\label{eq:DA3_3}
\end{align}
\end{subequations}
and, for $i=1,2,3$,
\begin{align}
D_i^{(B1)}[f](lmn)&= -C_{lmn}\{D_i^{(A1)}[f](lmn)\}   \label{eq:DB1},  \\
D_i^{(B2)}[f](lmn)&= D_i^{(A1)}[f](lmn)                      \label{eq:DB2}. \\
D_i^{(A2)}[f](lmn) &= -C_{lmn}\{D_i^{(B2)}[f](lmn)\}  \label{eq:DA2}.
\end{align}
The short hand notation $(lmn)$ refers to a site that belongs to the strain grid and whose coordinates are $(l,m,n)$, $(lmn+\frac{abc}{2})$ refers to a site that sits on the displacement grid and whose coordinates are $(l+\frac{a}{2},m+\frac{b}{2},n+\frac{c}{2})$ and $\bar 1$ is equal to $-1$. The operator $C_{lmn}$ in Eqs.~(\ref{eq:DB1}) and (\ref{eq:DA2}) represents the central symmetry around site $(lmn)$. Finally, we impose periodic boundary conditions (PBC) on the microstructure, i.e. on the stiffness $\stiff(\R)$ and eigenstrain $\eps^0(\R)$ which are defined on the SC strain grid, and on the fluctuating displacement fields $\dep^{(k)}(\R)$, $k=1,2$, which are defined on the SC displacement grid:
\begin{align}
\label{eq:PBC}
\dep^{(k)}(l+N_1,m,n) = \dep^{(k)}(l,m+N_2,n) = \dep^{(k)}(l,m,n+N_3) = \dep^{(k)}(l,m,n) \quad ( k=1,2),
\end{align}
and similar conditions on $\stiff(\R)$ and $\eps^0(\R)$. This automatically guaranties that the averages of the fluctuating strain fields $\deps^{(Kk)}(\R)$ defined in Eq.~(\ref{eq:deps}) are equal to zero. Also, the PBC are transmitted to $\deps^{(Kk)}(\R)$, and therefore to $\tend{\epsilon}^{(Kk)}(\R)$, through Eq.~(\ref{eq:deps}) and, consequently, to the stress field $ \stress^{(Kk)}(\R)$ through the constitutive law, here given by Eq.~(\ref{eq:linear_elasticity}).

Next, in order to write the equilibrium equations in Fourier space, we define the discrete Fourier transforms associated to our FCC subgrids.
%
%%%%%%%%%%%%%%%%%%%%%%%%%%%%%%%%%%%%%%%%%%%%%%%%%%%%%%%%%%%%%%
\subsection{Discrete Fourier transforms}\label{subsec:FT}
%%%%%%%%%%%%%%%%%%%%%%%%%%%%%%%%%%%%%%%%%%%%%%%%%%%%%%%%%%%%%%
%
Fourier transforms on FCC displacement subgrids $k=1,2$ and FCC strain subgrids  $K=A,B$, are defined by
\begin{subequations}
 \label{eq:FT}
 \begin{align}
\dep^{(k)}(\Q)&=\frac1N\sum_{\R \in k } \dep^{(k)}(\R) \exp(-i\Q.\R)   \label{eq:FT_subgrid_dep},  \\
\tend{A}^{(Kk)}(\Q)&=\frac1N\sum_{\R \in K} \tend{A}^{(Kk)}(\R) \exp(-i\Q.\R)   \label{eq:FT_subgrid_eps},
\end{align}
\end{subequations}
where $\tend{A}^{(Kk)}(\R)$ stands for any 2nd-order tensor field defined on subgrid $K$ and associated to discretization scheme $k$. These Fourier transforms are defined for all q-vectors that account for the PBC on the SC grids defined in Eq.~(\ref{eq:PBC}), i.e for the $N=N_1\times N_2 \times N_3$ q-vectors defined (in unit of $1/d$) by
\begin{equation}
\label{eq:q-vectors}
\begin{split}
\Q &= q_1\,\tens{a}_1^* + q_2\,\tens{a}_2^* +  q_3\,\tens{a}_3^*, \\
 q_i &= 2\pi  \frac{h_i}{N_i} \,\text{ and } \, h_i = 0,...,N_i-1,    \quad ( i=1,2,3),
\end{split}
\end{equation}
where $(\tens{a_1}^*,\tens{a_2}^*,\tens{a_3}^*)$ is the reciprocal basis associated to $(\tens{a_1},\tens{a_2},\tens{a_3})$\footnote{ Basis $(\tens{a_1}^*,\tens{a_2}^*,\tens{a_3}^*)$ and $(\tens{a_1},\tens{a_2},\tens{a_3})$ are linked by $\tens{a_i}^*.\tens{a_j}=\delta_{ij}$ where $\delta_{ij}$ is the Kronecker symbol. In the present situation, as $(\tens{a_1},\tens{a_2},\tens{a_3})$ is orthonormal, we have $\tens{a_i}^*=\tens{a_i}, i=1,2,3$.}. In the following, the short hand notation $(q_1q_2q_3)^*$ will refer to vector $\Q=q_1\tens{a_1}^*+q_2\tens{a_2}^*+q_3\tens{a_3}^*$. We note that a sum involved in Eq.~(\ref{eq:FT_subgrid_dep}) or  Eq.~(\ref{eq:FT_subgrid_eps}) does not run on a full SC grid, but only on the FCC subgrid on which the corresponding field is defined. Consequently, as an FCC subgrid contains only $\frac{N}{2}$ nodes, the $N$ Fourier components defined in Eq.~(\ref{eq:FT_subgrid_dep}) or in Eq.~(\ref{eq:FT_subgrid_eps}) are linked by relations that make them dependent\footnote{We mean here relations other than the conjugaison relations $\overline{\dep^{(k)}(\Q_1)}=\dep^{(k)}(\Q_2)$ and $\overline{\tend{A}^{(Kk)}(\Q_1)}=\tend{A}^{(Kk)}(\Q_2)$ where $\Q_1+\Q_2=\Q_0$ with $\Q_0=2(\pi00)^*$, $2(0\pi0)^*$ or $2(00\pi)^*$, which are due to the fact that the  fields $\dep^{(k)}(\R)$ and $\tend{A}^{(Kk)}(\R)$ are real.}. This point will be discussed later on. 

The inverse Fourier transforms associated to the direct Fourier transforms defined in Eqs.~(\ref{eq:FT_subgrid_dep}) and (\ref{eq:FT_subgrid_eps}) are given by
\begin{subequations}
\label{eq:IFT_1}
\begin{align}
&\R \in k,  \quad   \dep^{(k)}(\R) =\sum_{\Q} \dep^{(k)}(\Q) \exp(i\Q.\R), \label{eq:IFT_1_dep} \\
&\R \in K, \,\, \tend{A}^{(Kk)}(\R) =\sum_{\Q} \tend{A}^{(Kk)}(\Q) \exp(i\Q.\R),   \label{eq:IFT_1_A}
\end{align}
\end{subequations}
where the sums run on the q-vectors defined in Eq.~(\ref{eq:q-vectors}). Note that these inverse Fourier transforms are defined only for the $\frac{N}{2}$ sites that belong to the relevant FCC subgrid. In other words, an inverse Fourier transform computed on a site $\R$ that does not belongs to the relevant FCC subgrid is equal to zero:
\begin{subequations}
\label{eq:IFT_2}
\begin{align}
&\R \notin k:  \quad   \sum_{\Q} \dep^{(k)}(\Q) \exp(i\Q.\R)  = 0,\label{eq:IFT_2_dep} \\
&\R \notin K : \sum_{\Q} \tend{A}^{(Kk)}(\Q) \exp(i\Q.\R) = 0. \label{eq:IFT_2_A}
\end{align}
\end{subequations}

Because of the FCC nature of the subgrids, Fourier transforms defined in Eq.~(\ref{eq:FT}) and their inverse counterparts defined in Eq.~(\ref{eq:IFT_1}) have specific properties that we discuss now. Let $\Omega$ be the set of q-vectors defined by 
\begin{equation}
\label{eq:Omega}
\Omega = \{ \tens{Q} : \tens{Q}=\pi(s_1s_2s_3)^*; \,s_i=-1,1; \,i=1,2,3\}.
\end{equation}
Any vector $\tens{Q}$ that belongs to $\Omega$ belongs also to the Reciprocal Lattice of the FCC subgrids, i.e. fulfils the relation $\exp i\tens{Q}.(\R-\R')=1$ for any vectors $\R$ and $\R'$ that sit on the same FCC subgrid. In other words, for a given FCC subgrid, say for exemple the  FFC displacement subgrid $k$, the Fourier basis function $\exp(i\tens{Q}.\R)$ does not depend on $\R$ when $\tens{Q}$ belongs to $\Omega$. This means that there exist two functions $C^{(k)}(\tens{Q})$, $k=1,2$, such that 
\begin{align}
\label{eq:functions_C}
\forall \R \in k ,\quad \forall \tens{Q} \in \Omega,\quad    \exp(-i\tens{Q}.\R) = C^{(k)}(\tens{Q}),
\end{align}
and two functions $E^{(K)}(\tens{Q})$, $K=A,B$, such that
\begin{align}
\label{eq:functions_E}
\forall \R \in K ,\quad \forall \tens{Q} \in \Omega,\quad    \exp(-i\tens{Q}.\R) = E^{(K)}(\tens{Q}).
\end{align}
It is easily seen that, for any $\tens{Q}$ in $\Omega$, i.e. for $\tens{Q}=\pi(s_1s_2s_3)^*$ with $s_i=\pm 1$, functions $C^{(k)}(\tens{Q})$ and $E^{(K)}(\tens{Q})$  are given by \footnote{Functions $C^{(1)}(\tens{Q})$, $C^{(2)}(\tens{Q})$ and $E^{(B)}(\tens{Q})$ differ from 1 because the origin of real space is not located on a site of subgrids $k=1$,  $k=2$ or  $K=B$ (see Fig.~(\ref{fig:tetrahedral_stencil})), but on a site that belongs to subgrid $K=A$, which is the reason why $E^{(A)}(\tens{Q})$ is equal to 1. }
\begin{subequations}
\label{eq:functions_C_and_E}
\begin{align}
C^{(1)}(\pi(s_1s_2s_3)^*) &= \exp(-\frac{i}{2} \pi (s_1+s_2+s_3))              \label{eq:functions_C1} \\
C^{(2)}(\pi(s_1s_2s_3)^*) &= \exp(+\frac{i}{2} \pi (s_1+s_2+s_3))             \\
E^{(A)}(\pi(s_1s_2s_3)^*) &= +1                                                                \\
E^{(B)}(\pi(s_1s_2s_3)^*) &= -1                    
\end{align}
\end{subequations}

The first consequence of properties (\ref{eq:functions_C}) and (\ref{eq:functions_E})  is that, if two vectors $\Q_1$ and $\Q_2$ defined in (\ref{eq:q-vectors}) are related by $\Q_2=\Q_1+\tens{Q}$ where $\tens{Q}$ belongs to $\Omega$, Fourier components $\dep^{(k)}(\Q_2)$ and $\dep^{(k)}(\Q_1)$ are linked to each other:
\begin{align}
\label{eq:Fourier_component_relations_1}
\text{ if } \Q_2=\Q_1+\tens{Q} \text{ with }  \tens{Q} \in \Omega, \quad \dep^{(k)}(\Q_2) = \dep^{(k)}(\Q_1)\, C^{(k)}(\tens{Q}).
\end{align}
Siimilarly, for any 2nd-order tensor $\tend{A}^{(Kk)}$ defined on subgrid $K$ and associated to discretization scheme $k$, we have:
\begin{align}
\label{eq:Fourier_component_relations_2}
\text{ if } \Q_2=\Q_1+\tens{Q} \text{ with }  \tens{Q} \in \Omega, \quad \tend{A}^{(Kk)}(\Q_2) = \tend{A}^{(Kk)}(\Q_1)\, E^{(K)}(\tens{Q}).
\end{align}
Now,  it is straightforward to realize that, within the set of q-vectors defined in (\ref{eq:q-vectors}),  any vector $\Q_1$ is uniquely related to a vector $\Q_2$ such that $\Q_2-\Q_1$ belongs to $\Omega$. This, together with Eq.~(\ref{eq:Fourier_component_relations_1}) (respectively (\ref{eq:Fourier_component_relations_2})), generates exactly $\frac{N}{2}$ relations that linked together the $N$ Fourier components defined in Eq.~(\ref{eq:FT_subgrid_dep}) (respectively in Eq.~(\ref{eq:FT_subgrid_eps})), as anticipated above.

Next, another important consequence of property (\ref{eq:functions_C}) is that, for $\Q=(\pi\pi\pi)^*$, which is the unique q-vector within the set (\ref{eq:q-vectors}) that also belongs to the set $\Omega$, and therefore to the Reciprocal Lattice of the FCC subgrids, the contribution, in the inverse Fourier transform (\ref{eq:IFT_1_dep}), of the Fourier component $\dep^{(k)}(\Q)$ to the displacement field $\dep^{(k)}(\R)$ is independent of $\R$, i.e., we have, for $k=1,2$:
\begin{equation}
\label{eq:q_0}
\text{ if }  \R \text{ and } \R' \in k \text{ and } \Q=(\pi\pi\pi)^* , \quad \dep^{(k)}(\Q)\exp(i\Q.\R)=\dep^{(k)}(\Q)\exp(i\Q.\R').
\end{equation}
In other words, the mode $\Q=(\pi\pi\pi)^*$ corresponds simply to uniform translations in real space, as does the origin $\Q=(000)^*$ of Fourier space. We will refer to this result later when we discuss the Fourier transforms of the differential operators associated to our stencil and again when we discuss the stability of the mechanical equilibrium generated by the corresponding tetrahedral discretisation scheme.

We may now write the Fourier transform versions of Eq.~(\ref{eq:deps}). We obtain, for $k=1,2$ and $K=A,B$,
\begin{align}
\deps^{(Kk)}(\Q)& = \tens{D}^{(Kk)}(\Q) \otimes_s \dep^{(k)}(\Q)),
\label{eq:deps_Q}
\end{align}
where the symbols $\otimes_s$ represents the symmetrized tensor product. The functions $\tens{D}^{(Kk)}(\Q)$, for $k=1,2$ and $K=A,B$, are the Fourier transforms of the discrete operators defined in Eqs.~(\ref{eq:DA1}-\ref{eq:DA2}). They are given by, for $i=1,2,3$,
\begin{subequations}
\label{eq:DKkq}
\begin{align}
D_i^{(A1)}(\Q) &= D_i(\Q)                       \label{eq:DA1q}, \\
D_i^{(B1)}(\Q) &= -\overline{D_i(\Q)}      \label{eq:DB1q}, \\
D_i^{(B2)}(\Q) &= D_i(\Q)                       \label{eq:DB2q}, \\
D_i^{(A2)}(\Q) &= -\overline{D_i(\Q)}     \label{eq:DA2q}.
\end{align}
\end{subequations}
where the operators $D_i(\Q)$, $i=1,2,3$, are given by
\begin{subequations}
\label{eq:Dq}
\begin{align}
D_1(\Q) &= \frac12 \{ \exp(\frac{i}{2}(q_1+q_2+q_3) + \exp(\frac{i}{2}(q_1-q_2-q_3)- \exp(\frac{i}{2}(-q_1+q_2-q_3) - \exp(\frac{i}{2}(-q_1-q_2+q_3) \}  \label{eq:Dq_1},\\ 
D_2(\Q) &= \frac12 \{ \exp(\frac{i}{2}(q_1+q_2+q_3) - \exp(\frac{i}{2}(q_1-q_2-q_3) + \exp(\frac{i}{2}(-q_1+q_2-q_3) - \exp(\frac{i}{2}(-q_1-q_2+q_3) \} \label{eq:Dq_2}, \\
D_3(\Q) &=\frac12 \{ \exp(\frac{i}{2}(q_1+q_2+q_3) - \exp(\frac{i}{2}(q_1-q_2-q_3) - \exp(\frac{i}{2}(-q_1+q_2-q_3) + \exp(\frac{i}{2}(-q_1-q_2+q_3)\},  \label{eq:Dq_3}
\end{align}
\end{subequations}
for all vectors $\Q = q_1\,\tens{a}_1^* + q_2\,\tens{a}_2^* +  q_3\,\tens{a}_3^* $ defined in Eq.~(\ref{eq:q-vectors}). For latter reference, we mention that these operators fulfil the following property:
\begin{align}
\label{eq:Fourier_component_relations_3}
\text{ if } \Q_2=\Q_1+\tens{Q} \text{ with }  \tens{Q} \in \Omega, \quad \tens{D}(\Q_2) = \tens{D}(\Q_1)\, \overline{C^{(1)}(\tens{Q})}.
\end{align}
where the function $C^{(1)}(\tens{Q})$ has been given in Eq. (\ref{eq:functions_C1}) and $\Omega$ is the set defined in (\ref{eq:Omega}).

Now, we note that, as expected, the operators $D_i(\Q)$, $i=1,2,3$, are equal to zero for $\Q = (000)^*$. This results from the fact that mode $\Q=(000)^*$, which corresponds to a uniform translation in real space, does not generate any strain. Therefore, the Fourier transform of any differential operator, whether discrete or continuous, vanishes at the origin of Fourier space.

Meanwhile, we also observe that, in the present situation, the operators $D_i(\Q)$, $i=1,2,3$, are equal to zero for $\Q = (\pi\pi\pi)^*$, which implies that, within our tetrahedral stencil, mode $\Q = (\pi\pi\pi)^*$ also cannot generate any strain. However, this is not problematic because as shown above (see Eq.~(\ref{eq:q_0}) and corresponding discussion), within our discretisation scheme which is based on FCC subgrids,  $\Q = (\pi\pi\pi)^*$ belongs to the Reciprocal Lattice of the subgrids and, therefore, corresponds to uniform translations of the discrete displacements subgrids, which cannot generate any strain. In other words, the stencil used here has the important property that the Fourier transforms of the associated differential operators vanish only when they should, i.e. for q-modes that corresponds to uniform translations.

%
%%%%%%%%%%%%%%%%%%%%%%%%%%%%%%%%%%%%%%%%%%%%%%%%%%%%%%%%%%%%%%
\subsection{Mechanical equilibrium and Lippmann-Schwinger equations in Fourier space}\label{subsec:Equi_in_FS}
%%%%%%%%%%%%%%%%%%%%%%%%%%%%%%%%%%%%%%%%%%%%%%%%%%%%%%%%%%%%%%
%
We have now the material to write the Fourier space version of the equilibrium equations associated to our stencil. As discussed above, we consider simultaneously the discretization schemes associated to the two displacement subgrids $k=1,2$. Formally, the discretization scheme $k$ is associated to the elastic energy
\begin{equation}
E^{(k)} = \frac{d^3}{2} \sum_{K=A,B} \,\, \sum_{\R \in K} \,(\eps^{(Kk)}(\R)-\eps^0(\R)) : \stiff(\R) : (\eps^{(Kk)}(\R)-\eps^0(\R)),
\label{eq:Elast_energy_k}
\end{equation}
which corresponds to the initial constitutive law considered here, i.e.
\begin{equation}
\label{eq:linear_elasticity_Kk}
\stress^{(Kk)}(\R) = \frac{1}{d^3} \frac{\partial E^{(k)}}{\partial \eps^{(Kk)}(\R)} =  \stiff(\R):(\eps^{(Kk)}(\R) -\eps^0(\R)).
\end{equation}
When the system is subjected to an overall applied strain, the average strain $\bar \eps$ introduced in Eq.~(\ref{eq:eps}) is simply equal to the applied strain. Within this boundary condition, the mechanical equilibrium corresponds solely to the minimization of the discrete energy functional with respect to the fluctuating displacement field. Boundary conditions corresponding to an applied stress will be discuss later. Until otherwise stated, we consider here that the system is subjected to an applied strain. Mechanical equilibrium is then given by:
\begin{equation}
\forall \R \in k,  \quad  \frac{\partial E^{(k)} }{\partial \dep^{(k)}(\R)} = 0,  \quad  k=1,2.
\label{eq:Equilibrium_k_real}
\end{equation}
Using the inverse Fourier transform (\ref{eq:IFT_1_dep}) and the chain rule, we immediately get that, in Fourier space, these equations imply the following ones :
\begin{equation}
\forall \Q, \quad  \frac{\partial E^{(k)} }{\partial \dep^{(k)}(\Q)} = 0, \quad  k=1,2.
\label{eq:Equilibrium_k_Q}
\end{equation}
A remark is necessary. As an FCC subgrid contains only $\frac{N}{2}$ nodes, Eq.~(\ref{eq:Equilibrium_k_real}) represents $\frac{N}{2}$  independent equations for each discretisation scheme $k$, whereas Eq.~(\ref{eq:Equilibrium_k_Q}), which is defined for each of the $N$ q-vectors defined in Eq.~(\ref{eq:q-vectors}), represents $N$ equations. However, as discussed in the previous section (see Eqs.~(\ref{eq:Fourier_component_relations_1})-(\ref{eq:Fourier_component_relations_2})) and the corresponding discussion), there are in fact only $\frac{N}{2}$ independent Fourier components $\dep^{(k)}(\Q)$, which reduces the number of independent equations in Eq.~(\ref{eq:Equilibrium_k_Q}) to $\frac{N}{2}$. We return to this point at the end of this section, where we we will make explicit the dependancies between the Fourier components $\dep^{(k)}(\Q)$.

Now, a straightforward calculation (see \ref{sec:appendixA}) shows that the derivative of energy $E^{(k)}$ with respect to Fourier component $\dep^{(k)}(\Q)$ is given by
\begin{equation}
\quad \frac{\partial E^{(k)} }{\partial \dep^{(k)}(\Q)} = N d^3 \sum_{K=A,B}   \overline{\stress^{(Kk)}(\Q)}.{\tens{D}^{(Kk)}(\Q)},  \quad  k=1,2,
\label{eq:delta_Energy_k_wt_u_Q}
\end{equation}
where $\stress^{(Kk)}(\Q)$ is the Fourier transform of the stress field $\stress^{(Kk)}(\R)$. Therefore, the equilibrium given by Eq.~(\ref{eq:Equilibrium_k_Q}) may be written as:
\begin{equation}
\forall \Q \ne (000)^*, \quad  \sum_{K=A,B}  \tend{\sigma}^{(Kk)}(\Q).\overline{\tens{D}^{(Kk)}(\Q)}  = 0, \quad  k=1,2.
\label{eq:Equilibrium_k_Q_div}
\end{equation}
The point $\Q=(000)^*$ has been excluded because this equation is trivially verified at the origin of the Fourier space, which is a consequence of the translational invariance associated to $\Q=(000)^*$. Now, following the route proposed by Moulinec-Suquet \cite{moulinec1994fast,moulinec1998numerical} and by Khachaturyan \cite{khachaturyan1995elastic}, we rewrite Eq.~(\ref{eq:Equilibrium_k_Q_div}) in a form that allows it to be solved using an iterative algorithm. The strategy consists in introducing a homogeneous reference medium with constant stiffness tensor $\tenq{\lambda}^0$ and the associated polarization tensors, which are defined on the same nodal points than the strain and stress fields. Because the corresponding SC grid has been split into two FCC subgrids $K=A,B$, there are two polarization tensors for each discretization $k=1,2$. In Fourier space, they are defined by
\begin{align}
 \pol^{(Kk)}(\Q)  &= \stress^{(Kk)}(\Q) - \tenq{\lambda}^0:\eps^{(Kk)}(\Q), \quad k=1,2, \quad K=A,B.
\label{eq:pol}
\end{align}
Inserting these polarization tensors into Eq.~(\ref{eq:Equilibrium_k_Q_div}) and using Eqs.~(\ref{eq:deps_Q}) leads to the Fourier space version of the Lippmann-Schwinger equations for the displacement fields,
\begin{equation}
\forall \Q \ne (000)^*, \quad  \tend{\Omega}^{(k)}(\Q)^{-1}.\dep^{(k)}(\Q) = - \sum_{K=A,B} \pol^{(Kk)}(\Q).\overline{\tens{D}^{(Kk)}(\Q)}, \quad  k=1,2,
\label{eq:Lippmann_Schwinger_1}
\end{equation}
where $\tend{\Omega}^{(k)}(\Q)^{-1}$, the inverse of the discrete Green function associated with the discretization scheme $k$, is given by
\begin{equation}
 \tend{\Omega}^{(k)}(\Q)^{-1} = \sum_{K=A,B} \overline{\tens{D}^{(Kk)}(\Q)}.\tenq{\lambda}^0. \tens{D}^{(Kk)}(\Q).
\label{eq:Green_tensor_k}
\end{equation}
 The inverse of the two Green functions, $\tend{\Omega}^{(k)}(\Q)^{-1}$ for $k=1,2$, are in fact identical. Indeed, for $k=1$, Eq.~(\ref{eq:Green_tensor_k}) reads
\begin{equation}
\tend{\Omega}^{(1)}(\Q)^{-1} = \overline{\tens{D}^{(A1)}(\Q)}.\tenq{\lambda}^0. \tens{D}^{(A1)}(\Q) +  \overline{\tens{D}^{(B1)}(\Q)}.\tenq{\lambda}^0. \tens{D}^{(B1)}(\Q),
\end{equation}
which, using Eqs.~(\ref{eq:DA1q}) and (\ref{eq:DB1q}), leads to
\begin{equation}
\tend{\Omega}^{(1)}(\Q)^{-1} = \overline{\tens{D}(\Q)}.\tenq{\lambda}^0. \tens{D}(\Q) +  \tens{D}(\Q).\tenq{\lambda}^0. \overline{\tens{D}(\Q)},
\label{eq:Green_tensor_1}
\end{equation}
whereas, for $k=2$, Eq.~(\ref{eq:Green_tensor_k}) reads
\begin{equation}
\tend{\Omega}^{(2)}(\Q)^{-1} = \overline{\tens{D}^{(A2)}(\Q)}.\tenq{\lambda}^0. \tens{D}^{(A2)}(\Q) +  \overline{\tens{D}^{(B2)}(\Q)}.\tenq{\lambda}^0. \tens{D}^{(B2)}(\Q),
\end{equation}
which, using Eqs.~(\ref{eq:DB2q}) and (\ref{eq:DA2q}), leads to
\begin{equation}
\tend{\Omega}^{(2)}(\Q)^{-1} = \overline{\tens{D}(\Q)}.\tenq{\lambda}^0. \tens{D}(\Q) +  \tens{D}(\Q).\tenq{\lambda}^0. \overline{\tens{D}(\Q)}.
\label{eq:Green_tensor_2}
\end{equation}
When they are defined, the Green functions $\tend{\Omega}^{(1)}(\Q)$ and $\tend{\Omega}^{(2)}(\Q)$ are therefore identical and noted in the following $\tend{\Omega}(\Q)$. Obviously, $\tend{\Omega}(\Q)$ is real and symmetric:
\begin{align}
\label{eq:displacement_Green_function_symmetry}
\Omega_{ij}(\Q) = \overline{\Omega_{ij}(\Q)}= \Omega_{ji}(\Q).
\end{align}
For latter reference, we mention also that, as an immediate consequence of property (\ref{eq:Fourier_component_relations_3}), $\tend{\Omega}(\Q)$ fulfils the following property:
\begin{align}
\label{eq:Fourier_component_relations_4}
\text{ if } \Q_2=\Q_1+\tens{Q} \text{ with }  \tens{Q} \in \Omega, \quad \tend{\Omega}(\Q_2) = \tend{\Omega}(\Q_1).
\end{align}

Now, we recall that, for $\Q=(\pi\pi\pi)^*$, the operator $\tens{D}(\Q)$ is equal to zero, which makes the inverse of the Green functions also equal to zero, leaving the Green function  $\tend{\Omega}(\Q)$  undefined for this mode. However, as explained at the end of section (\ref{subsec:FT}), this mode is irrelevant as, within our discretisation scheme, it corresponds to uniform translations of the discrete displacement subgrids and, therefore, do not generate any strain. Hence, we can safely set  $\tend{\Omega}(\Q)=\tend{0}$ for $\Q=(\pi\pi\pi)^*$ and, obviously, also for $\Q=(000)^*$, and write the Lippmann-Schwinger equations (\ref{eq:Lippmann_Schwinger_1}) under the form:
\begin{equation}
\forall \Q, \quad  \dep^{(k)}(\Q) = - \tend{\Omega}(\Q) \, . \,\{\sum_{K=A,B} \pol^{(Kk)}(\Q).\overline{\tens{D}^{(Kk)}(\Q)}\,\}, \quad k=1,2.
\label{eq:Lippmann_Schwinger_2}
\end{equation}
We stress that the situation for the point $\Q=(\pi\pi\pi)^*$ discussed above is different from the one encountered, for example, when the discretization method is based on the centered scheme or on the rotated scheme proposed by Willot \cite{willot2015fourier}. In these discretization schemes, the differential operators are also equal to zero for some non-zero q-vector. However, contrary to the present situation, these q-vectors do not sit on the Reciprocal Lattice of the chosen real space displacement grid. Consequently, cancellation of the differential operators at these q-points is problematic, as it makes the Green functions and, therefore, the corresponding Fourier components of the displacement and strain fields undefined, even though the Fourier components of the strain remain finite in the neighbourhood of these q-points. A recipe proposed by Willot consists in setting to zero the Green functions when the differential operators vanish, which is equivalent to setting to zero the correspondant Fourier components of the strain. However, it does not solve the problem, because, as shown in \ref{sec:appendixB}, it leaves the Fourier transforms of the Green function, and therefore of the strain and stress fields, non-analytic at the corresponding q-points. This, as we will illustrate below, generates inevitably non-physical instabilities in real space, such as oscillations and cherkerboard effects on the strain and stress fields. In contrast, the tetrahedral scheme proposed here, which guarantees that the operators $D_i(\Q)$ are different from zero for any q-vector that does not sit on the Reciprocal Lattice of the real-space subgrids, is automatically exempt of these artefacts.

At this stage, we have two equilibrium equations (\ref{eq:Lippmann_Schwinger_2}), each associated with one of the discretization schemes $k=1,2$. In principle, we could solve these two Lippmann-Schwinger equations using some iterative scheme, such as the one proposed by Moulinec-Suquet. However, such a procedure would not be optimised from the computational point of view. Indeed, the Fourier transforms are here defined on FCC subgrids (see Eqs.~(\ref{eq:FT_subgrid_dep} and \ref{eq:FT_subgrid_eps})). Consequently, the 3D-lattice sums that appear in these equations cannot be split into three one-dimensional sums, each running along one index. In other words, if site $(lmn)$ appears in a sum, site $(l+1,m,n)$ does not, contrary to what is required by the numerical libraries implementing fast Fourier transform algorithms. A straightforward way to circumvent this difficulty is simply to set to zero the value of a field on sites that do not belong to the relevant subgrid (as exemple, set $\dep^{(k)}(\R)=\tens{0}$ if $\R \notin k$) and perform all Fourier transforms on full SC grids. However, this will obviously double the necessary computational and memory load.

In fact, we show in the next paragraph that the two Lippmann-Schwinger equations (\ref{eq:Lippmann_Schwinger_2}) can be combined into a single one. Also, this procedure will result in Fourier transforms involving complete SC grids, which naturally lend themselves to FFT algorithms, without the need to fill unnecessary information with zeros.

But before, we return to the remark we made after Eq.~(\ref{eq:Equilibrium_k_Q}), when we stressed that, for a given FCC displacement subgrid $k$, there are only $\frac{N}{2}$ independent Fourier coefficients $\dep^{(k)}(\Q)$ among the $N$ ones that are associated to the set of q-vectors defined in (\ref{eq:q-vectors}). This property should in fact be apparent from the Lippmann-Schwinger equations (\ref{eq:Lippmann_Schwinger_2}). Indeed, for the vector $\Q_1 +\tens{Q}$, where $\Q_1$ is any vector within the set defined in (\ref{eq:q-vectors}) and $\tens{Q}$ the unique vector within the set $\Omega$ defined in (\ref{eq:Omega}) such that $\Q_1 +\tens{Q}$ belongs also to the set  (\ref{eq:q-vectors}), the Lippmann-Schwinger equation
\begin{align}
\dep^{(k)}(\Q_1 + \tens{Q}) &= - \tend{\Omega}(\Q_1 + \tens{Q}) \, . \,\{ \pol^{(Ak)}(\Q_1 + \tens{Q}).\overline{\tens{D}^{(Ak)}(\Q_1 + \tens{Q})} +  \pol^{(Bk)}(\Q_1 + \tens{Q}).\overline{\tens{D}^{(Bk)}(\Q_1 + \tens{Q})}\,\}
\end{align}
becomes
\begin{align}
\label{eq:LS_equivalence}
\dep^{(k)}(\Q_1 + \tens{Q}) &= - \tend{\Omega}(\Q_1) \, . \,\{ \pol^{(Ak)}(\Q_1).\overline{\tens{D}^{(Ak)}(\Q_1)}E^{(A)}(\tens{Q})C^{(k)}(\tens{Q}) +  \pol^{(Bk)}(\Q_1).\overline{\tens{D}^{(Bk)}(\Q_1)}E^{(B)}(\tens{Q})\overline{C^{(k)}(\tens{Q})}\,\}, 
\end{align}
where we have used property (\ref{eq:Fourier_component_relations_4}), which is fulfilled by the Green function $\tend{\Omega}(\Q)$, property (\ref{eq:Fourier_component_relations_2}), which is here used for the polarization tensors $\pol^{(Kk)}(\Q)$ and, finally, the definitions of the operators $\tens{D}^{(Kk)}(\Q)$ given in Eqs.( \ref{eq:DKkq}) together with property (\ref{eq:Fourier_component_relations_3}). From the definitions of the functions $C^{(k)}$ and $E^{(K)}$ given in Eq. (\ref{eq:functions_C_and_E}), it is easy to see that 
\begin{equation}
 \forall \tens{Q} \in \Omega,\quad   E^{(A)}(\tens{Q})C^{(k)}(\tens{Q}) =  E^{(B)}(\tens{Q})\overline{C^{(k)}((\tens{Q}))} = C^{(k)}(\tens{Q}).
\end{equation}
Therefore, equation (\ref{eq:LS_equivalence}) becomes 
\begin{equation}
\begin{split}
\dep^{(k)}(\Q_1 + \tens{Q}) &= - \tend{\Omega}(\Q_1) \, . \,\{ \pol^{(Ak)}(\Q_1).\overline{\tens{D}^{(Ak)}(\Q_1)}+  \pol^{(Bk)}(\Q_1).\overline{\tens{D}^{(Bk)}(\Q_1)}\,\} \,C^{(k)}(\tens{Q}), \\
&= \dep^{(k)}(\Q_1)\,C^{(k)}(\tens{Q}).
\end{split}
\end{equation}
As expected, we recover the dependancies between the Fourier components of the displacement fields derived previously in (\ref{eq:Fourier_component_relations_1}) as consequences of the FCC nature of the subgrids and that let, for each subgrid $k$, only $\frac{N}{2}$ independent Fourier component $\dep^{(k)}(\Q)$.

In the next paragraph, we reformulate the two mechanical equilibria and corresponding Lippmann-Schwinger equations, which each involves only $\frac{N}{2}$ independent Fourier components, into a single one that simultaneously involves $N$ independent Fourier coefficients.
%
%
%%%%%%%%%%%%%%%%%%%%%%%%%%%%%%%%%%%%%%%%%%%%%%%%%%%%%%%%%%%%%%
\subsection{Collapse of the two equilibrium equations}\label{sec:collapse}
%%%%%%%%%%%%%%%%%%%%%%%%%%%%%%%%%%%%%%%%%%%%%%%%%%%%%%%%%%%%%%
%
Consider the energy functional $E$ defined as the average of the elastic energies $E^{(k)}$ defined in Eq.~(\ref{eq:Elast_energy_k}):
\begin{equation}
E = \frac12( E^{(1)} + E^{(2)} ).
\label{eq:Elast_energy}
\end{equation}
As $E^{(k)}$ depends only on the displacement field $\dep^{(k)}(r)$, we have
\begin{equation}
\R \in k,  \quad \frac{\partial E}{\partial \dep^{(k)}(\R)} = \frac12 \, \frac{\partial E^{(k)}}{\partial \dep^{(k)}(\R)}.
\end{equation}
Therefore, the two variational problems stated in Eq.~(\ref{eq:Equilibrium_k_real}) are equivalent to a single one given by
\begin{equation}
\R \in (1)+(2), \quad \frac{\partial E}{\partial \dep(\R)} =  0,
\label{eq:Equilibrium_real}
\end{equation}
where the notation $(1)+(2)$ refers to the full SC displacement grids and where the field $\dep(\R)$ simply collects the two subgrid fields $\dep^{(k)}(\R)$, $k=1,2$, i.e.,
\begin{equation}
\label{eq:full_dep}
\text{if } \R \in k, \quad \dep(\R) = \dep^{(k)}(\R).
\end{equation}
 Whereas each variational problem expressed in Eq.~(\ref{eq:Equilibrium_k_real}) involves $\frac{N}{2}$ equations and $\frac{N}{2}$ unknowns,  Eq.~(\ref{eq:Equilibrium_real}) involves obviously $N$ equations and $N$ unknowns. 
  
The variational problem (\ref{eq:Equilibrium_real}) may be written in Fourier space. We first introduce the Fourier transform of the "full" displacement field $\dep(\R)$,
\begin{equation}
\dep(\Q)=\frac1N\sum_{\R \in (1)+(2) } \dep(\R) \exp(-i\Q.\R),
\label{eq:FT_u}
\end{equation} %
and the corresponding inverse Fourier transform
\begin{equation}
\dep(\R)=\sum_{\Q} \dep(\Q) \exp(i\Q.\R).
\label{eq:IFT_u}
\end{equation} 
The Fourier transform $\dep(\Q)$ defined in Eq.~(\ref{eq:FT_u}) involves the displacement field $\dep(\R)$ defined on all the sites of the SC displacement grid, as opposed to the restricted Fourier transforms introduced previously, that involve only the relevant FCC subgrids (see Eq.~(\ref{eq:FT_subgrid_dep})). Consequently, the $N$ Fourier components $\dep(\Q)$ are obviously independent, contrary to the situation encountered above for the coefficients $\dep^{(k)}(\Q)$, which were linked by the relations given in Eq.~(\ref{eq:Fourier_component_relations_1}). Therefore, equations (\ref{eq:Equilibrium_real}) are strictly equivalent to the following ones:
\begin{equation}
\forall \Q, \quad  \frac{\partial E}{\partial \dep\Q)} = 0.
\label{eq:Equilibrium_Q}
\end{equation}
Now, using the inverse Fourier transform (\ref{eq:IFT_u}) and the chain rule, the partial derivative of energy $E$ with respect to Fourier component $\dep(\Q)$ may be written as
\begin{equation}
\label{eq:der_E_dep_Q}
\frac{\partial E}{\partial \dep(\Q)} =  \sum_{\R \ \in (1)+(2)}  \frac{\partial E}{\partial \dep(\R)} \exp(i\Q.\R),
\end{equation}
which, using Eqs.~(\ref{eq:Elast_energy}) and (\ref{eq:full_dep}) becomes
\begin{equation}
\frac{\partial E}{\partial \dep(\Q)} = \frac{1}{2}\{ \sum_{\R \ \in (1)}  \frac{\partial E^{(1)}}{\partial \dep^{(1)}(\R)} \exp(i\Q.\R) + \sum_{\R \ \in (2)}  \frac{\partial E^{(2)}}{\partial \dep^{(2)}(\R)} \exp(i\Q.\R)\,\},
\end{equation}
where we have used the fact that energy $E^{(1)}$ (respectively, $E^{(2)}$) does not depend on displacement field $\dep^{(2)}(\R)$ (respectively, $\dep^{(1)}(\R)$). Next, using the inverse Fourier transforms of $\dep^{(1)}_i(\Q)$ and $\dep^{(2)}_i(\Q)$ given in Eq.~(\ref{eq:IFT_1_dep}), the previous equation becomes
\begin{equation}
\frac{\partial E}{\partial \dep(\Q)} = \frac{1}{2}\{ \frac{\partial E^{(1)} }{\partial \dep^{(1)}(\Q)} + \frac{\partial E^{(2)} }{\partial \dep^{(2)}(\Q)} \,\},
\end{equation}
which, using Eq.~(\ref{eq:delta_Energy_k_wt_u_Q}), leads to
\begin{equation}
\frac{\partial E}{\partial \dep(\Q)}  = \frac{Nd^3}{2}\,\sum_{k=1,2} \sum_{K=A,B}  \overline{\tend{\sigma}^{(Kk)}(\Q)}.\tens{D}^{(Kk)}(\Q).
\label{eq:der_E_dep_Q_bis}
\end{equation}
Then, using the expressions of the differential operators given in Eq.~(\ref{eq:DKkq}), we may recast this equation as
\begin{equation}
\frac{\partial E}{\partial \dep(\Q)}  = \frac{Nd^3}{2}\ \{\overline{\tend{\sigma}^{(T1)}(\Q)}.{\tens{D}(\Q)}  - \overline{\tend{\sigma}^{(T2)}(\Q)}.\overline{{\tens{D}(\Q)}}\},
\label{eq:der_E_dep_Q_ter}
\end{equation}
 where the quantities $\tend{\sigma}^{(T1)}(\Q)$ and $\tend{\sigma}^{(T2)}(\Q)$ are defined by
\begin{subequations}
\label{eq:sigma_T1_T2_Q}
\begin{align}
\stress^{(T1)}(\Q) &=  \stress^{(A1)}(\Q) + \stress^{(B2)}(\Q)            \label{eq:sigma_T1_Q}, \\
\stress^{(T2)}(\Q) &=  \stress^{(A2)}(\Q) + \stress^{(B1)}(\Q)            \label{eq:sigma_T2_Q}.
\end{align}
\end{subequations}
We will comment on the meaning of the superscripts $T1$ and $T2$ latter. Finally, the equilibrium in Fourier space given by Eqs.~(\ref{eq:Equilibrium_Q}) takes the form 
\begin{equation}
\forall \Q \ne (000)^*,   \quad  \tend{\sigma}^{(T1)}(\Q).\overline{\tens{D}(\Q)}  - \tend{\sigma}^{(T2)}(\Q).{\tens{D}(\Q)} = 0,
\label{eq:Equilibrium_Q_div_bis}
\end{equation}
where, as in the previous section, the point $\Q=(000)^*$ is excluded because of translational invariance in real space.

We now formulate the equilibrium equation (\ref{eq:Equilibrium_Q_div_bis}) into a Lippmann-Schwinger form. First, in a similar way to the introduction of the stress fields $\stress^{(T1)}(\Q)$ and $\stress^{(T2)}(\Q)$, we define the q-space fields 
\begin{subequations}
\label{eq:X_T1_T2_Q}
\begin{align}
\tend{X}^{(T1)}(\Q) &=  \tend{X}^{(A1)}(\Q) + \tend{X}^{(B2)}(\Q)            \label{eq:X_T1_Q} ,\\
\tend{X}^{(T2)}(\Q) &=  \tend{X}^{(A2)}(\Q) + \tend{X}^{(B1)}(\Q)            \label{eq:X_T2_Q} ,
\end{align}
\end{subequations}
where $\tend{X}$ stands for a strain $\eps$, a fluctuating strain $\deps$ or a polarization tensor $\pol$. Of course, these definitions are such that the polarization tensors $\pol^{(Ti)}(\Q)$ are linked to the corresponding stress and strain fields by the same form as in Eq.~(\ref{eq:pol}):
\begin{align}
\pol^{(Ti)}(\Q)  = \stress^{(Ti)}(\Q) - \tenq{\lambda}^0:\eps^{(Ti)}(\Q), \quad i=1,2.
\label{eq:pol_T1_T2_Q}
\end{align}
Also, using Eq.~(\ref{eq:deps_Q}) that relates the fluctuating subgrid strains to the fluctuating subgrid displacement fields, together with Eqs.~(\ref{eq:DKkq}) that gives the differential operators and the obvious relation $\dep(\Q)=\dep^{(1)}(\Q)+\dep(^{(2)}\Q)$, we have
\begin{subequations}
\label{eq:deps_dep_Q}
\begin{align}
\deps^{(T1)}(\Q) &=    \tens{D}(\Q) \otimes_s \dep(\Q),                     \label{eq:desp_T1_dep_Q} \\
\deps^{(T2)}(\Q) &=    -\overline{\tens{D}(\Q)} \otimes_s \dep(\Q).    \label{eq:desp_T2_dep_Q} 
\end{align}
\end{subequations}
Now, inserting the polarization tensors given in Eq.~(\ref{eq:pol_T1_T2_Q}) into Eq.~(\ref{eq:Equilibrium_Q_div_bis}) and using Eqs.~(\ref{eq:deps_dep_Q}) leads to the Fourier space version of a Lippmann-Schwinger equation for the displacement field $\dep(\Q)$, that is
\begin{equation}
\forall \Q \ne (000)^*, \quad \tend{\Omega}(\Q)^{-1}.\dep(\Q) = - \{\pol^{(T1)}(\Q).\overline{\tens{D}(\Q)} - \pol^{(T2)}(\Q).\tens{D}(\Q)\},
\label{eq:Lippmann_Schwinger_u_inverse}
\end{equation}
where, for clarity, the inverse of the Green function $\tend{\Omega}(\Q)^{-1}$, which has already been given above (see the earlier Eqs.~(\ref{eq:Green_tensor_1}) and (\ref{eq:Green_tensor_2})), is again reproduced here:
 \begin{equation}
\tend{\Omega}(\Q)^{-1} = \overline{\tens{D}(\Q)}.\tenq{\lambda}^0. \tens{D}(\Q) +  \tens{D}(\Q).\tenq{\lambda}^0. \overline{\tens{D}(\Q)}.
\label{eq:Green_tensor}
\end{equation}
For the same reason as above, we can safely set $\tend{\Omega}(\Q)=0$ when $\tend{\Omega}(\Q)^{-1}$ is equal to zero, i.e. for $\Q=(000)^*$ and $\Q=(\pi\pi\pi)^*$, and write the Lippmann-Schwinger equation under the form
\begin{equation}
\forall \Q, \quad   \dep(\Q) = \tend{\Omega}(\Q). \{- \pol^{(T1)}(\Q).\overline{\tens{D}(\Q)} + \pol^{(T2)}(\Q).\tens{D}(\Q)\},
\label{eq:Lippmann_Schwinger_u_direct}
\end{equation}
Note that the Fourier components $\dep(\Q)$ for the modes $\Q=(000)^*$ and $\Q=(\pi\pi\pi)^*$ are therefore chosen to be equal to zero, which is fully compatible with their purely translational character. 
The Lippmann-Schwinger equation (\ref{eq:Lippmann_Schwinger_u_direct}) can then be solved iteratively using, for exemple, the basic scheme, which reads
\begin{equation}
\forall \Q, \quad   \dep(\Q)_{(n+1)} = \tend{\Omega}(\Q). \{- \pol_{(n)}^{(T1)}(\Q).\overline{\tens{D}(\Q)} + \pol_{(n)}^{(T2)}(\Q).\tens{D}(\Q)\},
\label{eq:basic_scheme_dep_1}
\end{equation}
This algorithm can be further simplified by using the following property,
\begin{equation}
 \tend{\Omega}(\Q)\, . \,  \tenq{\lambda}^0 :\{ \tens{D}(\Q) \otimes_s \dep(\Q).\overline{\tens{D}(\Q)} + \overline{\tens{D}(\Q)} \otimes_s \dep(\Q).\tens{D}(\Q)\} = \dep(\Q),
\end{equation}
which, because we set $\dep(\Q)$ to zero when $\Omega(\Q)$ is equal to zero, is valid for any q-vector. The iterative scheme then becomes
\begin{equation}
\forall \Q, \quad   \dep(\Q)_{(n+1)} = \dep(\Q)_{(n)} - \tend{\Omega}(\Q). \{ \stress_{(n)}^{(T1)}(\Q).\overline{\tens{D}(\Q)} - \stress_{(n)}^{(T2)}(\Q).\tens{D}(\Q)\},
\label{eq:basic_scheme_dep_2}
\end{equation}
which, compared with the scheme (\ref{eq:basic_scheme_dep_1}), avoids the calculation and storage of the polarisation fields $\pol^{(T1)}(\Q)$ and  $\pol^{(T2)}(\Q)$. It also explicitly uses, in the last term on the right-hand side, a quantity which, measuring in q-space the deviation from mechanical equilibrium (see Eq.~(\ref{eq:Equilibrium_Q_div_bis})), is needed to test the convergence of the iterative scheme.

We now comment on the meaning of the superscripts $T1$ and $T2$. Let $X$ be any one of the 2nd-order tensors that appear in Eqs.~(\ref{eq:sigma_T1_T2_Q}) and (\ref{eq:X_T1_T2_Q}), i.e.  a stress $\stress$, a strain $\eps$, a fluctuating strain $\deps$ or a polarization tensor $\pol$. According to the definitions given in Eqs.~(\ref{eq:sigma_T1_T2_Q}) and (\ref{eq:X_T1_T2_Q}) and using Eqs.~(\ref{eq:IFT_1_A}) and (\ref{eq:IFT_2_A}), the inverse Fourier transforms of $\tend{X}^{(Ti)}(\Q)$, i.e.
\begin{equation}
\tend{X}^{(Ti)}(\R) =\sum_{\Q} \tend{X}^{(Ti)}(\Q) \exp(i\Q.\R), \quad i=1,2,    \\
\label{eq:InvFT_X}
\end{equation}
are given by:
\begin{subequations}
\label{eq:X_T1_T2_R}
\begin{align}
\tend{X}^{(T1)}(\R) &=  \tend{X}^{(A1)}(\R) \,\text{ if }\, \R \in A,  \quad  \tend{X}^{(T1)}(\R) = \tend{X}^{(B2)}(\R) \,\text{ if }\, \R \in B,           \\
\tend{X}^{(T2)}(\R) &=  \tend{X}^{(A2)}(\R) \,\text{ if }\, \R \in A,  \quad  \tend{X}^{(T2)}(\R) = \tend{X}^{(B1)}(\R) \,\text{ if }\, \R \in B.      
\end{align}
\end{subequations}
In other words, the superscripts $T1$ (respectively, $T2$) labels a quantity which is defined on each site of the SC strain grid and, concerning its link with the displacement field $\dep(\R)$, is associated to the sites of the SC displacement grid through only one type of tetrahedron, namely the blue tetrahedra (respectively, the red ones) that are displayed in Fig.~(\ref{fig:tetrahedral_stencil}). In other words, our tetrahedral stencil can be interpreted without recourse to the FCC displacements subgrids $k=1,2$ and FCC strain subgrids $K=A,B$ and without introducing the corresponding reduced Fourier transforms, even if its mathematical formulation must go through these lattice subdivisions, notably to show that, thanks to its tetrahedral character, the stencil is mathematically stable, i.e. automatically free from any oscillation or checkerboard effect that would be observed if the associated q-space operators were to vanish on any q-point not belonging to the relevant Reciprocal Lattice.

In summary, as a result of the present reformulation, our tetrahedral stencil leads to a single Lippmann-Schwinger equation that involves quantities defined on full SC grids. In addition, the Fourier transforms, which are now also defined on SC grids, are naturally suited to FFT algorithms. The underlying FCC character of the stencil is reflected in the appearance of two strain and stress fields, labelled above with the indices $T1$ and $T2$. Once the Lippmann-Schwinger equation has been solved, the final stress and strain fields are simply defined as the arithmetic averages of the corresponding fields:
\begin{subequations}
\label{eq:final_fields}
\begin{align}
\stress(\R) &= \frac12 \{ \stress^{(T1)}(\R) + \stress^{(T2)}(\R) \}  \label{eq:final_stress},\\
\eps(\R) &= \frac12 \{ \eps^{(T1)}(\R) + \eps^{(T2)}(\R) \} .      \label{eq:final_strain}
\end{align}
\end{subequations}
Of course, these two fields are linked by the relation
\begin{equation}
\label{eq:sigma_epsilon}
\stress(\R)=\stiff(\R):(\eps(\R)-\eps^0(\R)),
\end{equation}
which is inherited from Eqs.~(\ref{eq:linear_elasticity_Kk}) through Eq.~(\ref{eq:X_T1_T2_R}) and (\ref{eq:final_fields}). Note also that the same relation holds between $\stress^{(T1)}(\R)$ and $\eps^{(T1)}(\R)$ on the one hand, and between $\stress^{(T2)}(\R)$ and $\eps^{(T2)}(\R)$  on the other, and that $\eps^{(T1)}(\R)$ and $\eps^{(T2)}(\R)$ are related to their corresponding fluctuating strain fields through the same decomposition as in Eq.~(\ref{eq:eps}) :
\begin{subequations}
\label{eq:eps_T1_T2}
\begin{align}
\eps^{(T1)}(\R)= \bar \eps + \deps^{(T1)}(\R), \\
\eps^{(T2)}(\R)= \bar \eps + \deps^{(T2)}(\R) .
\end{align}
\end{subequations}
Finally, we note that, since the definitions given in Eqs.~(\ref{eq:final_fields}) only involve quantities defined on the same site $\R$, the averaged fields $\stress(\R)$ and $\eps(\R)$ are no less accurate than the initial ones, as might be observed if some spatial averaging process had been used (see for example \cite{Willot2008FastFT}).
%
%%%%%%%%%%%%%%%%%%%%%%%%%%%%%%%%%%%%%%%%%%%%%%%%%%%%%%%%%%%%%%%%%%%%%%%%%%%%%%%%%%%%%%%%%%%%%
\subsection{Iterative loop on the average strain when the system is subjected to an applied stress}
%%%%%%%%%%%%%%%%%%%%%%%%%%%%%%%%%%%%%%%%%%%%%%%%%%%%%%%%%%%%%%%%%%%%%%%%%%%%%%%%%%%%%%%%%%%%%
%
The boundary conditions considered above correspond to an applied strain, which is identified with the average strain $\bar \eps$. When the boundary conditions consists in applying an external stress, the elastic energy given in Eq.~(\ref{eq:Elast_energy}) must be supplemented by an external contribution,
\begin{equation}
E = \frac12( E^{(1)} + E^{(2)} )- Nd^3\stress^a:\bar \eps,
\label{eq:Grand_Elast_energy}
\end{equation}
where $\stress^a$ is the applied stress. Minimisation of this energy with respect to $\bar \eps$ leads to
\begin{equation}
\label{eq:equi_epsilon_bar}
<\stress(\R)> = \sigma^a ,
\end{equation}
where $\stress(\R)$ is the local stress defined in Eq.~(\ref{eq:final_stress}) and  $<X>$ represents the spatial average of $X$. This equation may be easily recast in a form that allows it to be solved using an iterative algorithm. This can be done as follows. Using successively Eqs.~(\ref{eq:sigma_epsilon}), (\ref{eq:final_strain}) and (\ref{eq:eps_T1_T2}), Eq.~(\ref{eq:equi_epsilon_bar}) leads to
\begin{equation}
<\stiff(\R)>:\bar \eps=\stress^a - <\stiff(\R):\deps(\R)> + <\stiff(\R):\eps^0(\R)>,
\end{equation}
where $\deps(\R)=\frac12(\deps^{(T1)}(\R)+\deps^{(T2)}(\R))$. Shifting the left member of this equation to the right and adding the product $\stiff^0 : \bar \eps$ to each member gives
\begin{equation}
 \stiff^0 : \bar \eps =  \stress^a + (\stiff^0 -  <\stiff(\R)>):\bar \eps - <\stiff(\R):\deps(\R)> + <\stiff(\R):\eps^0(\R)>,
\end{equation}
which suggests the following iterative scheme
\begin{equation}
\bar \eps_{(n+1)} =   (\stiff^0)^{-1} : \{\stress^a + (\stiff^0 -  <\stiff(\R)>):\bar \eps_{(n)} - <\stiff(\R):\deps(\R)_{(n)} > + <\stiff(\R):\eps^0(\R)> \},
\end{equation}
where $\bar \eps_{(n)} $ and $\deps(\R)_{(n)}$ are the average and fluctuating strain fields at iteration $n$.

%
%%%%%%%%%%%%%%%%%%%%%%%%%%%%%%%%%%%%%%%%%%%%%%%%%%%%%%%%%%%%%%%%%%%%%%%%%%%%%%%%%%%%%%%%%%%%%
\subsection{Summary}
%%%%%%%%%%%%%%%%%%%%%%%%%%%%%%%%%%%%%%%%%%%%%%%%%%%%%%%%%%%%%%%%%%%%%%%%%%%%%%%%%%%%%%%%%%%%%
%
To highlight the simplicity of the overall method, we show in Algorithm (\ref{alg:dep}) the algorithm which implements the displacement-based iterative procedure given in Eq.~(\ref{eq:basic_scheme_dep_2}).

%--------------------------------------------------------------------------------------------------------------------------------------------------------------------------------------------
\begin{algorithm}
\caption{: displacement-based iterative method for the tetrahedral stencil on a simple cubic discretization scheme}\label{alg:dep}
\begin{algorithmic}[1]

   \Statex initial guess for $\dep(\Q)$ \quad (e.g. $\dep(\Q)$  = 0)
   \Statex if applied stress: initial guess for $\bar \eps$  \quad (e.g. $\bar \eps =   (\stiff^0)^{-1} : \stress^a$ where $ \stress^a$ is the applied stress)

   \medskip
   \State   $\deps^{(T1)}(\Q) =    \tens{D}(\Q) \otimes_s \dep(\Q)$                    \label{al:deps}
   \Statex $\deps^{(T2)}(\Q) =   -\overline{\tens{D}(\Q)} \otimes_s \dep(\Q $  

   \medskip
   \State   $\deps^{(T1)}(\R) = \text{FT}^{-1}\{\deps^{(T1)}(\Q) \}$                
   \Statex $\deps^{(T2)}(\R) = \text{FT}^{-1}\{\deps^{(T2)}(\Q) \}$

 %  \State   $\eps^{(T1)}(\R) = \bar \eps + \deps^{(T1)}(\R)$            
%   \Statex $\eps^{(T2)}(\R) = \bar \eps + \deps^{(T2)}(\R)$

%   \State   $\stress^{(T1)}(\R) = \stiff(\R):( \eps^{(T1)}(\R) -\eps^0(\R)) $
%   \Statex $ \stress^{(T2)}(\R) = \stiff(\R):( \eps^{(T2)}(\R) -\eps^0(\R))$

   \medskip
   \State   $\stress^{(T1)}(\R) = \stiff(\R):( \bar \eps +  \deps^{(T1)}(\R) -\eps^0(\R))$
   \Statex $\stress^{(T2)}(\R) = \stiff(\R):(\bar \eps +  \deps^{(T2)}(\R) -\eps^0(\R))$
    
  % \State   $\pol^{(T1)}(\R) =  \stress^{(T1)}(\R) - \stiff^0:\eps^{(T1)}(\R) $ 
 %  \Statex $\pol^{(T2)}(\R) =  \stress^{(T2)}(\R) - \stiff^0:\eps^{(T2)}(\R) $  

   \medskip
   \State   $\stress^{(T1)}(\Q) = \text{FT}\{\stress^{(T1)}(\R)\}$
   \Statex $\stress^{(T2)}(\Q) = \text{FT}\{\stress^{(T2)}(\R)\}$

   \medskip
   \State  $\dep(\Q)  \gets \,\dep(\Q) - \tend{\Omega}(\Q) . \{ \stress^{(T1)}(\Q).\overline{\tens{D}(\Q)} - \stress^{(T2)}(\Q).\tens{D}(\Q)\}$.    \label{al:dep}

   \medskip
   \State  if applied stress : 
   \Statex $\deps(\R)=\frac12(\deps^{(T1)}(\R)+\deps^{(T2)}(\R))$
   \Statex $\bar \eps \gets   (\stiff^0)^{-1} : \{\stress^a + (\stiff^0 -  <\stiff(\R)>):\bar \eps - <\stiff(\R):\deps(\R) > + <\stiff(\R):\eps^0(\R)> \} $

   \medskip
   \State return to $1$ until Convergence

\end{algorithmic}
\end{algorithm}
\FloatBarrier 
%--------------------------------------------------------------------------------------------------------------------------------------------------------------------------------------------
%
\noindent The operators $\tens{D}(\Q)$ that appears in steps \ref{al:deps} and \ref{al:dep} are given in Eq.~(\ref{eq:Dq}) and the Green function $\tend{\Omega}(\Q)$ is given in Eq.~(\ref{eq:Green_tensor}). The symbols $\text{FT}$ and $\text{FT}^{-1}$ denote Fourier transform and Inverse Fourier transform, respectively. 

Convergence can be evaluated using the L2-norm of the quantity that measures deviation from equilibrium, i.e.
\begin{equation}
\label{eq:L2-norm}
L_2 = \sqrt{\frac1N \, \sum_{\R \in (1)+(2)} \, \left \vert   \, \frac{\partial E}{d^3 \, \partial \dep(\R)} \right\vert^2},
\end{equation}
where the quantity $\frac{1}{d^3}\frac{\partial E}{\partial \dep(\R)}$ represents, in unit of $1/d$, the local force density that acts on site $\R$. Using Parseval's theorem and Eq.~(\ref{eq:der_E_dep_Q}), we get
\begin{equation}
\label{eq:L2-norm}
L_2 = \sqrt{ \sum_{\Q} \, \left\vert   \, \frac{\partial E}{Nd^3 \, \partial \dep(\Q)} \right\vert^2}
\end{equation}
which, using Eq.~(\ref{eq:der_E_dep_Q_ter}), leads to 
\begin{equation}
\label{eq:L2-norm}
L_2 = \sqrt{ \sum_q \, \left\vert \, \frac{\tend{\sigma}^{(T1)}(\Q).\overline{\tens{D}(\Q)}  - \tend{\sigma}^{(T2)}(\Q).{\tens{D}(\Q)}}{2} \, \right\vert^2 }
\end{equation}
Finally, to make the convergence test insensitive to the amplitude of the average stress $<\stress>$, we define the convergence criterion by 
\begin{equation}
\label{eq:tolerance}
\epsilon =  \frac{L_2}{ \vert\vert <\stress>   \vert\vert }
\end{equation}
where $\vert\vert <\stress>   \vert\vert  = \sqrt{ \sum_{i,j} \vert <\sigma_{ij}>\vert^2}$ is the Frobenius norm of the average stress.

Finally, for the sake of completeness, we translate now the previous algorithm, that iterates on the displacement field, into an algorithm that iterates on strains. The resulting scheme will be of course strictly equivalent, at each iteration step, to the initial one. From the Lippmann-Schwinger equation (\ref{eq:Lippmann_Schwinger_u_direct}), which relates the displacement field to the polarisation fields $\pol^{T1}(\Q)$ and $\pol^{T2}(\Q)$, and Eqs.~(\ref{eq:desp_T1_dep_Q}) and (\ref{eq:desp_T2_dep_Q}), which relate the fluctuating deformation fields $\deps^{(T1)}(\Q)$ and $\deps^{(T2)}(\Q)$ to the displacement field $\dep(\Q)$, we obtain, after a straightforward algebra,
\begin{equation}
\label{eq:Lippmann_Schwinger_strain_1}
\begin{split}
\deps^{(T1)}_{(n+1)}(\Q) &=   -\tenq{G}^d(\Q) : \pol^{T1}_{(n)}(\Q)  - \tenq{G}^{nd}(\Q) : \pol^{T2}_{(n)}(\Q),                                                \\
\deps^{(T2)}_{(n+1)}(\Q) &=   -\overline{\tenq{G}^{nd}(\Q)} : \pol^{T1}_{(n)}(\Q)  - \overline{\tenq{G}^{d}(\Q)} : \pol^{T2}_{(n)}(\Q),         
\end{split}
\end{equation}
where the equations are presented in iterative form. The diagonal and non-diagonal Green operators $\tenq{G}^d(\Q)$ and $\tenq{G}^{nd}(\Q)$ are defined by:
\begin{equation}
\label{eq:strain_based_Green_functions_tetrahedral_stencil}
\begin{split}
G^d_{ijkl}(\Q) &= \frac14 \{ D_i(\Q) \Omega_{jk}(\Q) \overline{D_l(\Q)} + D_j(\Q) \Omega_{ik}(\Q)  \overline{D_l(\Q)} +   D_i(\Q) \Omega_{jl}(\Q)  \overline{D_k(\Q)} + D_j(\Q) \Omega_{il}(\Q)  \overline{D_k(\Q)} \},         \\
G^{nd}_{ijkl}(\Q) &=  -\frac{1}{4} \{ D_i(\Q) \Omega_{jk}(\Q) D_l(\Q) + D_j(\Q) \Omega_{ik}(\Q)  D_l(\Q) +   D_i(\Q) \Omega_{jl}(\Q)  D_k(\Q) + D_j(\Q) \Omega_{il}(\Q)  D_k(\Q) \}.             
\end{split}
\end{equation}
We note that, as a consequence of property (\ref{eq:displacement_Green_function_symmetry}), these Green operators display the following symmetry relations:
\begin{equation}
\begin{split}
G^d_{ijkl}(\Q) &=  G^d_{jikl}(\Q)  =  G^d_{ijlk}(\Q) =  \overline{G^d_{klij}(\Q)} ,                                     \\
G^{nd}_{ijkl}(\Q) &=  G^{nd}_{jikl}(\Q)  =  G^{nd}_{ijlk}(\Q) =  G^{nd}_{klij}(\Q).           
\end{split}
\end{equation}
Finally, using the following properties,
\begin{equation}
\nonumber
\begin{split}
\tenq{G}^d(\Q) : \tenq{\lambda}^0 : \deps^{(T1)}(\Q) + \tenq{G}^{nd}(\Q) : \tenq{\lambda}^0 : \deps^{(T2)}(\Q) &=   \deps^{(T1)}(\Q),                                         \\
\overline{\tenq{G}^{nd}(\Q)} : \tenq{\lambda}^0 : \deps^{(T1)}(\Q) + \overline{\tenq{G}^{d}(\Q)} : \tenq{\lambda}^0 : \deps^{(T2)}(\Q) &=   \deps^{(T2)}(\Q),         
\end{split}
\end{equation}
which, because  $\deps^{(T1)}(\Q)$ and $ \deps^{(T2)}(\Q)$ are equal to zero when $\Q=0$, are valid for any q-vector, the iterative scheme (\ref{eq:Lippmann_Schwinger_strain_1}) can be written as
\begin{equation}
\label{eq:Lippmann_Schwinger_strain_2}
\begin{split}
\deps^{(T1)}_{(n+1)}(\Q) &=  \deps^{(T1)}_{(n)}(\Q)  - \tenq{G}^d(\Q) : \stress^{T1}_{(n)}(\Q)  - \tenq{G}^{nd}(\Q) : \stress^{T2}_{(n)}(\Q)                                        \\
\deps^{(T2)}_{(n+1)}(\Q) &=  \deps^{(T1)}_{(n)}(\Q)  - \overline{\tenq{G}^{nd}(\Q)} : \stress^{T1}_{(n)}(\Q)  - \overline{\tenq{G}^{d}(\Q)} : \stress^{T2}_{(n)}(\Q)               
\end{split}
\end{equation}
which, as compared to (\ref{eq:Lippmann_Schwinger_strain_1}), avoids the computation and storage of the polarization fields. This strain-based algorithm may be further reformulated as
\begin{equation}
\label{eq:Lippmann_Schwinger_strain_3}
\begin{split}
\deps^{(T1)}_{(n+1)}(\Q) &=  \deps^{(T1)}_{(n)}(\Q)  - \tent{\Gamma}(\Q) . \{ \stress_{(n)}^{(T1)}(\Q).\overline{\tens{D}(\Q)} - \stress_{(n)}^{(T2)}(\Q).\tens{D}(\Q)\}                                      \\
\deps^{(T2)}_{(n+1)}(\Q) &=  \deps^{(T2)}_{(n)}(\Q)  + \overline{\tent{\Gamma}(\Q)} . \{ \stress_{(n)}^{(T1)}(\Q).\overline{\tens{D}(\Q)} - \stress_{(n)}^{(T2)}(\Q).\tens{D}(\Q)\}                                                   
\end{split}
\end{equation}
where the third-order tensor $\tent{\Gamma}(\Q)$ is defined by
\begin{equation}
\Gamma_{ijk}(\Q) = \frac12 \{ D_i(\Q) \Omega_{jk}(\Q) + D_j(\Q) \Omega_{ik}(\Q) \}.      
\end{equation}
Note that this tensor is symmetric with respect to the first two indices. The advantage of algorithm (\ref{eq:Lippmann_Schwinger_strain_3}) over (\ref{eq:Lippmann_Schwinger_strain_2}) is that only one Green operator is needed and that it explicitly uses the  quantity which, measuring in q-space the deviation from mechanical equilibrium (see Eq.~(\ref{eq:Equilibrium_Q_div_bis})), is needed to test the convergence of the iterative scheme. Finally, the strain-based algorithm that results from (\ref{eq:Lippmann_Schwinger_strain_3}) is summarized in Algorithm \ref{alg:strain}. If formulation (\ref{eq:Lippmann_Schwinger_strain_2}) is preferred, equations in step \ref{step:strain_step} are replaced by Eqs.~(\ref{eq:Lippmann_Schwinger_strain_2}))

%--------------------------------------------------------------------------------------------------------------------------------------------------------------------------------------------
\begin{algorithm}
\caption{: strain-based iterative method for the tetrahedral stencil on a simple cubic discretization scheme}\label{alg:strain}
\begin{algorithmic}[1]

   \Statex initial guess for $\deps^{(T1)}(\Q)$ and $\deps^{(T2)}(\Q)$  \quad (e.g. $\deps^{(T1)}(\Q)$  = 0, $\deps^{(T2)}(\Q)$  = 0)
   \Statex if applied stress: initial guess for $\bar \eps$  \quad (e.g. $\bar \eps =   (\stiff^0)^{-1} : \stress^a$ where $ \stress^a$ is the applied stress)

   \medskip 
   \State   $\deps^{(T1)}(\R) = \text{FT}^{-1}\{\deps^{(T1)}(\Q) \}$               
   \Statex $\deps^{(T2)}(\R) = \text{FT}^{-1}\{\deps^{(T2)}(\Q) \}$

   \medskip 
   \State   $\stress^{(T1)}(\R) = \stiff(\R):( \bar \eps +  \deps^{(T1)}(\R) -\eps^0(\R))$
   \Statex $\stress^{(T2)}(\R) = \stiff(\R):(\bar \eps +  \deps^{(T2)}(\R) -\eps^0(\R))$

   \medskip 
   \State   $\stress^{(T1)}(\Q) = \text{FT}\{\stress^{(T1)}(\R)\}$
   \Statex $\stress^{(T2)}(\Q) = \text{FT}\{\stress^{(T2)}(\R)\}$

   \medskip
   \State   $ \deps^{(T1)}(\Q) \gets \,  \deps^{(T1)}(\Q)   - \tent{\Gamma}(\Q) . \{ \stress^{(T1)}(\Q).\overline{\tens{D}(\Q)} - \stress^{(T2)}(\Q).\tens{D}(\Q)\}  $     \label{step:strain_step}
   \Statex $ \deps^{(T2)}(\Q) \gets \,  \deps^{(T2)}(\Q)  + \overline{\tent{\Gamma}(\Q)} . \{ \stress^{(T1)}(\Q).\overline{\tens{D}(\Q)} - \stress^{(T2)}(\Q).\tens{D}(\Q)\}  $  
   
   \medskip 
   \State  if applied stress : 
   \Statex $\deps(\R)=\frac12(\deps^{(T1)}(\R)+\deps^{(T2)}(\R))$
   \Statex $\bar \eps \gets   (\stiff^0)^{-1} : \{\stress^a + (\stiff^0 -  <\stiff(\R)>):\bar \eps - <\stiff(\R):\deps(\R) > + <\stiff(\R):\eps^0(\R)> \} $

   \medskip  
   \State return to $1$ until Convergence

\end{algorithmic}
\end{algorithm}
\FloatBarrier 
%--------------------------------------------------------------------------------------------------------------------------------------------------------------------------------------------

%
%%%%%%%%%%%%%%%%%%%%%%%%%%%%%%%%%%%%%%%%%%%%%%%%%%%%%%%%%%%%%%%%%%%%%%%%%%%%%%%%%%%%%%%%%%%%%
\section{Numerical experiments}\label{sec:numerical_exp}
%%%%%%%%%%%%%%%%%%%%%%%%%%%%%%%%%%%%%%%%%%%%%%%%%%%%%%%%%%%%%%%%%%%%%%%%%%%%%%%%%%%%%%%%%%%%%
%
In this section, we analyse the performance of the tetrahedral stencil approach for different applications, including comparisons with the original Moulinec-Suquet \cite{moulinec1994fast,moulinec1998numerical} method and the rotated finite difference scheme proposed by Willot \cite{willot2015fourier}. Our aim is to test the accuracy and efficiency of the method and also its ability to treat, in term of convergence rate,  extreme situations where the elastic contrast is large or infinite.

%All the simulations reported below have been performed with the iterative algorithm described above, which follows the original basic scheme, here adapted to our tetrahedral stencil. 

All the results presented below, with the exception of those in section (\ref{sec:hard_inclusion}) that concern accelerated schemes and that have been taken from the literature, were obtained with the basic scheme which, for our tetrahedral stencil, was implemented in the form described in Algorithm \ref{alg:dep} and the FFT implementation rely upon the  FFTW library \cite{FFTW05}.
%
%%%%%%%%%%%%%%%%%%%%%%%%%%%%%%%%%%%%%%%%%%%%%%%%%%%%%%%%%%%%%%%%%%%%%%%%%%%%%%%%%%%%%%%%%%%%%
\subsection{A cubic inclusion with eigenstrain in an elastically isotropic homogeneous material}\label{sec:inclusion_with_eigenstrain}
%%%%%%%%%%%%%%%%%%%%%%%%%%%%%%%%%%%%%%%%%%%%%%%%%%%%%%%%%%%%%%%%%%%%%%%%%%%%%%%%%%%%%%%%%%%%%
%
For our first application, following \cite{eloh2019development}, we consider a homogeneous, elastically isotropic medium containing a cubic inclusion characterised by a finite eigenstrain $\eps^0$. The interest of this situation is that mechanical equilibrium is known analytically (see for example \cite{li2001compact}), which permits to test the accuracy of our numerical method.

 In the present situation, we use the same material parameters as in \cite{eloh2019development}. The media is elastically homogeneous and its shear modulus and Poisson ratio are given by $\mu = 132.30$ GPa and $\nu = 0.26$. The eigenstrain in the inclusion is such that $\epsilon^0_{ij} =0$ if $(ij) \ne (33)$ whereas $\epsilon^0_{33}$ is finite. The periodic medium is discretized with $256^3$ cells and the cubic inclusion consists into $31^3$ cells.
\begin{figure}[h]
	\centering
	\includegraphics[width=0.49\textwidth]{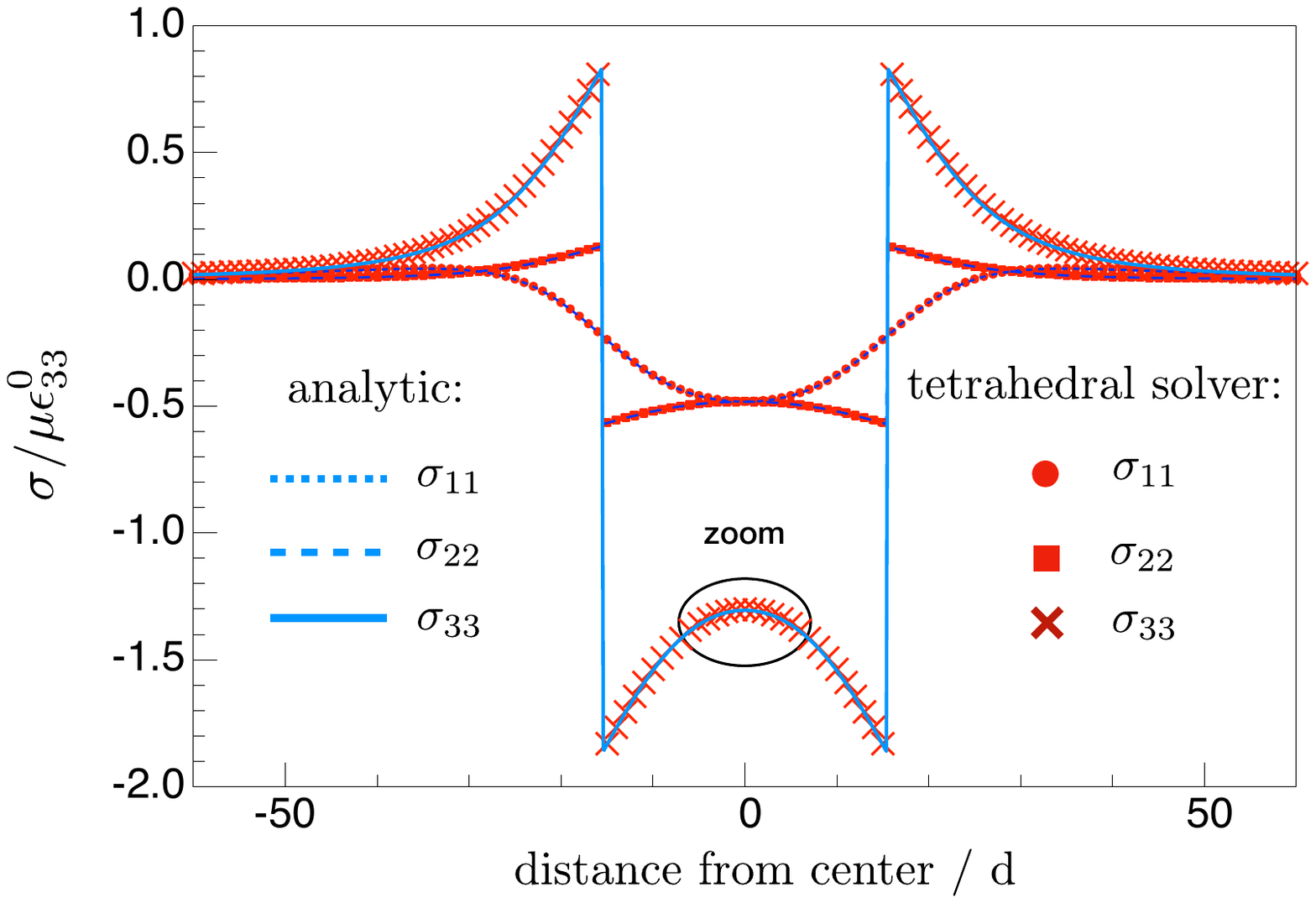}
	\includegraphics[width=0.49\textwidth]{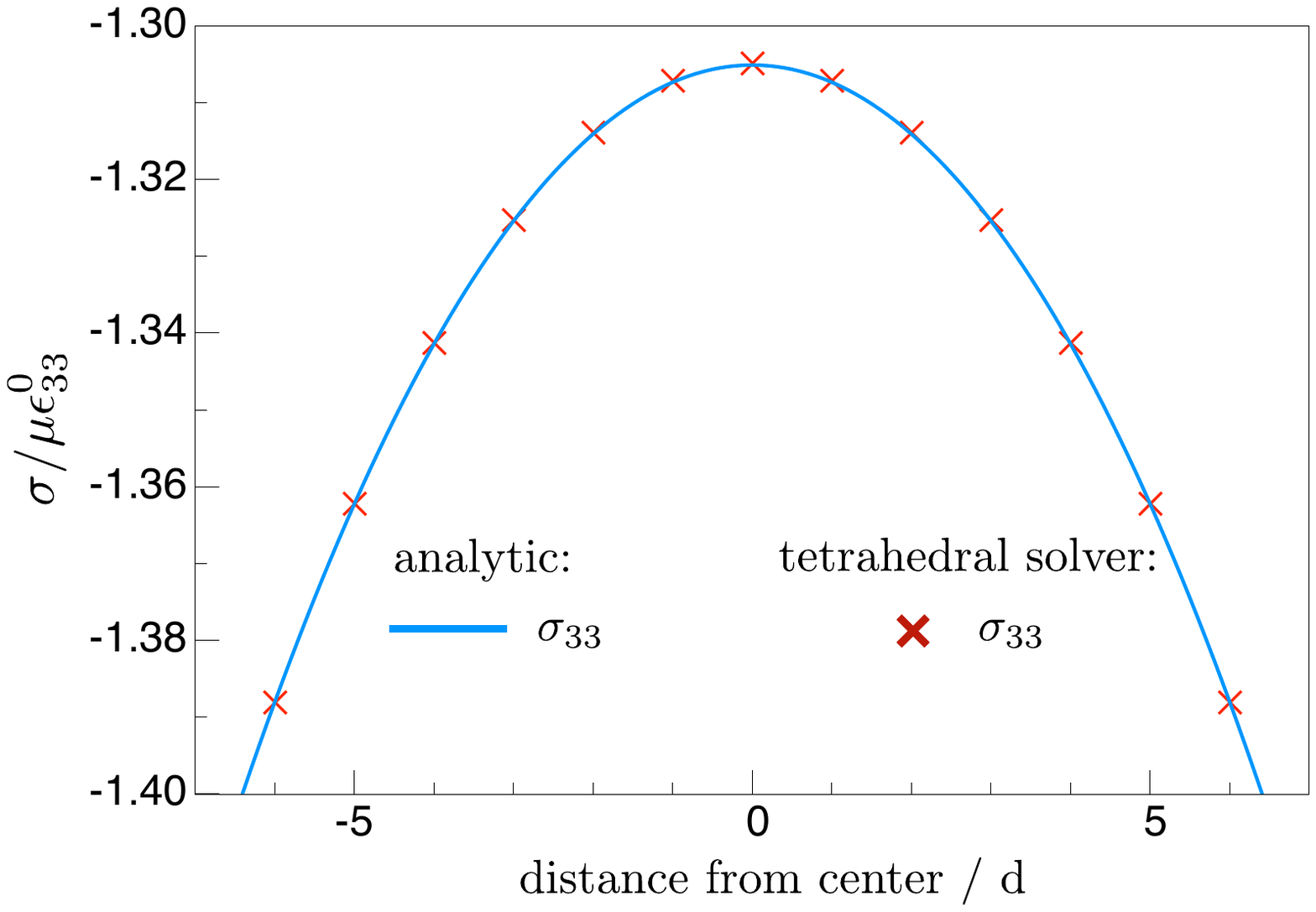}
	\caption{\textcolor{blue}{A cubic inclusion with eigenstrain $\epsilon^0_{33}$ in an isotropic medium with a zero applied stress: internal stresses $\sigma_{11}$,  $\sigma_{22}$ and  $\sigma_{33}$, in units of $\mu\epsilon^0_{33}$, along the line parallel to direction $(100)$ and that goes through the center of the cell. The numerical results obtained with the tetrahedral stencil (red symbols) are compared to the analytic solution (blue lines); the right panel shows a zoom on the small zone marked in the left panel.}
}
\label{fig:cubic_inclusion}
\end{figure}
The boundary condition used here is a zero external stress. The results obtained with the tetrahedral stencil are compared with the exact solution (obtained from \cite{li2001compact}) in Fig.~(\ref{fig:cubic_inclusion}), where we report the stress components $\sigma_{11}, \sigma_{22}, \sigma_{33}$  along a line parallel to the cubic direction $(100)$ and that goes through the center of the cubic inclusion. The comparison with the analytical solution shows that the tetrahedral stencil leads to very accurate numerical results, as observed in particular in the right panel of Fig.~(\ref{fig:cubic_inclusion}), where we display a zoom on a small region inside the inclusion.
%
%%%%%%%%%%%%%%%%%%%%%%%%%%%%%%%%%%%%%%%%%%%%%%%%%%%%%%%%%%%%%%%%%%%%%%%%%%%%%%%%%%%%%%%%%%%%%
\subsection{A cubic void within an elastically isotropic homogeneous matrix}\label{sec:cubic_void}
%%%%%%%%%%%%%%%%%%%%%%%%%%%%%%%%%%%%%%%%%%%%%%%%%%%%%%%%%%%%%%%%%%%%%%%%%%%%%%%%%%%%%%%%%%%%%
%
Next, we consider a cubic void embedded into an elastically isotropic matrix. Our aim is to test the behaviour of the tetrahedral stencil in the presence of infinite elastic contrast, a situation which is known to be challenging, both in term of convergence of the iterative algorithms and, if convergence is achieved, in term of accuracy. Indeed, it has long been observed that, in the original Moulinec-Suquet method, the high frequency cut-off generated by the discretization of the Lippmann-Schwinger equation generates oscillations whose amplitude increases with the elastic contrast. It has also been observed that the iterative basic scheme converges slowly for large contrast and does not converge if the contrast is infinite. 

We consider a material with the same elastic parameters as in the previous case. The unit cell is discretized into $64^3$ voxels and we insert a cubic cavity which consists in $31^3$ cells and whose faces are perpendicular to the cartesian axis (this relatively small size has been selected since it may typically correspond to linear dimensions of heterogeneities in a complex microstructure). The overall system is subjected to an hydrostatic pressure $P=300$ MPa.

In the following, we will compare the results obtained with the tetrahedral stencil to those obtained with the Moulinec-Suquet method and the rotated scheme used by Willot. We mention that, for these two methods and when needed, we use the regularisation recipes proposed by the corresponding authors to fix problems encountered with the behaviour of the respective Green operators. Indeed, within the Moulinec-Suquet discretisation, if the number of voxels along cartesian axis $i$ is even, conjugaison relations on the Fourier representation of the Green function are lost when the component $q_i$ of a q-vector reaches the highest frequency, i.e. when $q_i=\pi$, which makes the imaginary part of backward Fourier transforms of the strain field non zero. Following the recipe proposed by Moulinec-Suquet  \cite{moulinec1998numerical}, we then fix the Green function equal to the inverse of the elastic tensor $\tenq{\lambda}^0$. As for the rotated scheme, if the number of voxels along any two of the cartesian axis, say $i$ and $j$, are even, all the components of the q-space differential operator $\tens{D}(\Q)$ are equal to zero when components $q_i$ and $q_j$ of the q-vector reach the highest frequency, i.e. when $q_i=q_j=\pi$, which makes the Green function undefined at the corresponding q-points. Following the recipe proposed by Willot, we then fix the Green function to zero. We recall that this problem does not exist with our discretization scheme, because the tetrahedral stencil has the property that the q-space differential operator $\tens{D}(\Q)$ never vanish for $\Q\ne\tens{0}$, except at a q-point equivalent to the origin of Fourier space (see Eq. (\ref{eq:Lippmann_Schwinger_2}) and corresponding discussion). 
\begin{figure}[h]
	\centering
	\includegraphics[width=0.65\textwidth]{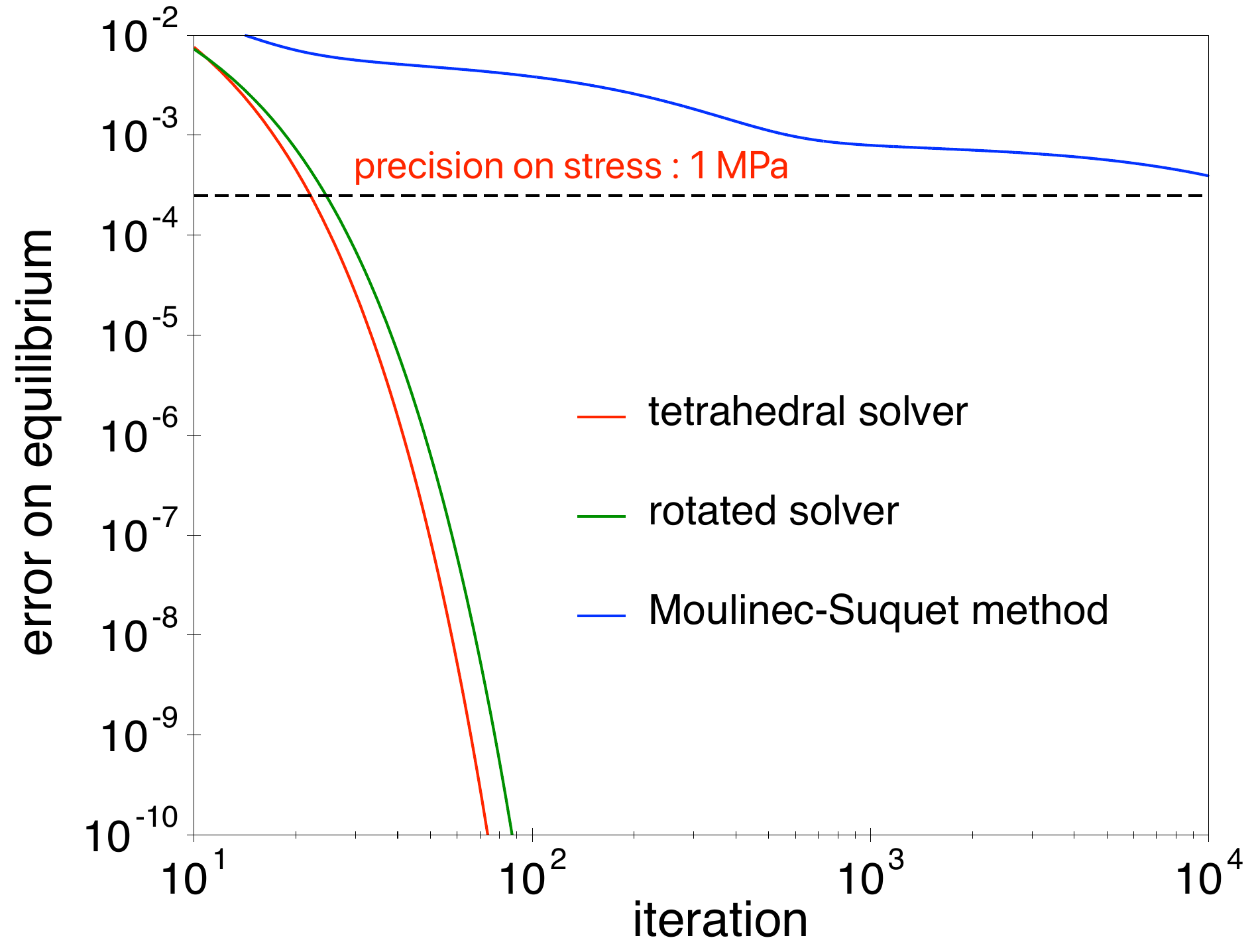}
	\caption{\textcolor{blue}{A cubic void within an isotropic medium subjected to an hydrostatic pressure $P=300$ MPa: error to equilibrium, defined in (\ref{eq:tolerance}),  as a function of iterates; comparison between the tetrahedral stencil (red), the rotated stencil (green) and Moulinec-Suquet method (blue). The dotted line indicates the precision required to reach, with the tetrahedral stencil, an accuracy of 1 MPa on the stress field.}}
\label{fig:cubic_void_error_to_equilibrium}
\end{figure}

First, we report in Fig.~(\ref{fig:cubic_void_error_to_equilibrium}) the evolution of the error on equilibrium, defined in Eq. (\ref{eq:tolerance}), as a function of the iterate index for the tetrahedral stencil and, for comparaison, for the rotated and Moulinec-Suquet methods. The stiffness tensor $\tenq{\lambda}^0$ of the reference medium has not been optimised and, for the three methods, has been simply set to $\tenq{\lambda}^0=0.8\,\tenq{\lambda}^m$ where $\tenq{\lambda}^m$ is the stiffness of the matrix. We observe that the tetrahedral and rotated schemes converge very rapidly, with the highest precision investigated here, $\epsilon = 10^{-10}$, being reached with less than $100$ iterates, with a slight advantage for the tetrahedral stencil, whereas the error with the Moulinec-Suquet method decreases very slowly and, in fact, does not appear to converge to zero, a fact that has already been observed when the elastic contrast is infinite, as in the present case. The highest precision investigated here is much more than necessary in common applications. In Fig.~(\ref{fig:cubic_void_error_to_equilibrium}), the dotted line indicates the precision required to reach, with the tetrahedral stencil, an accuracy of 1 MPa on the stress field, accuracy been defined as the maximum difference observed between the calculated stress field and that obtained in a high-precision reference calculation  ($\epsilon =10^{-10}$). We observe that only approximately 20 iterations are needed to achieve the required accuracy. 
\begin{figure}[h]
	\centering
	\includegraphics[width=0.68\textwidth]{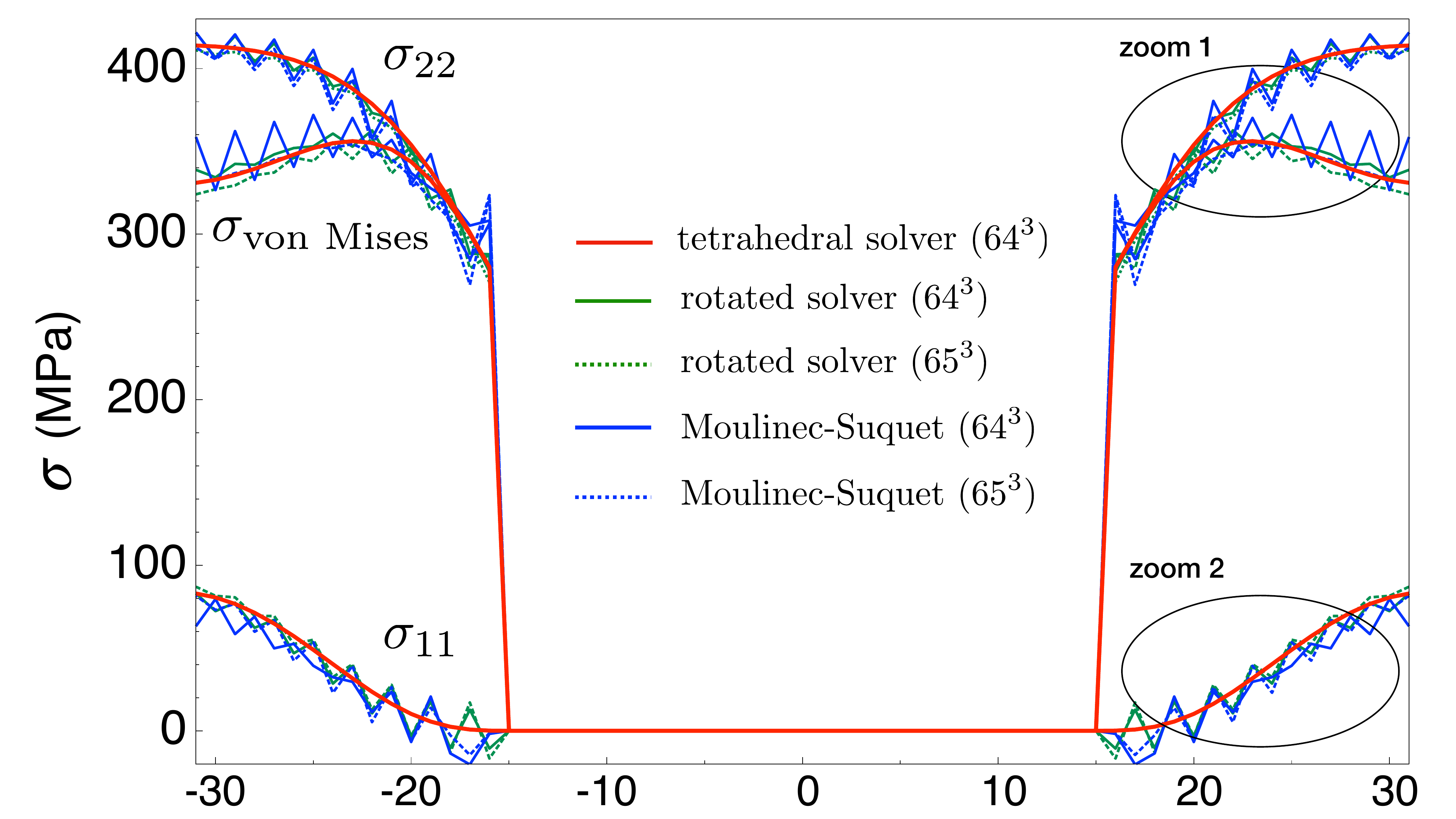}
	\includegraphics[width=0.68\textwidth]{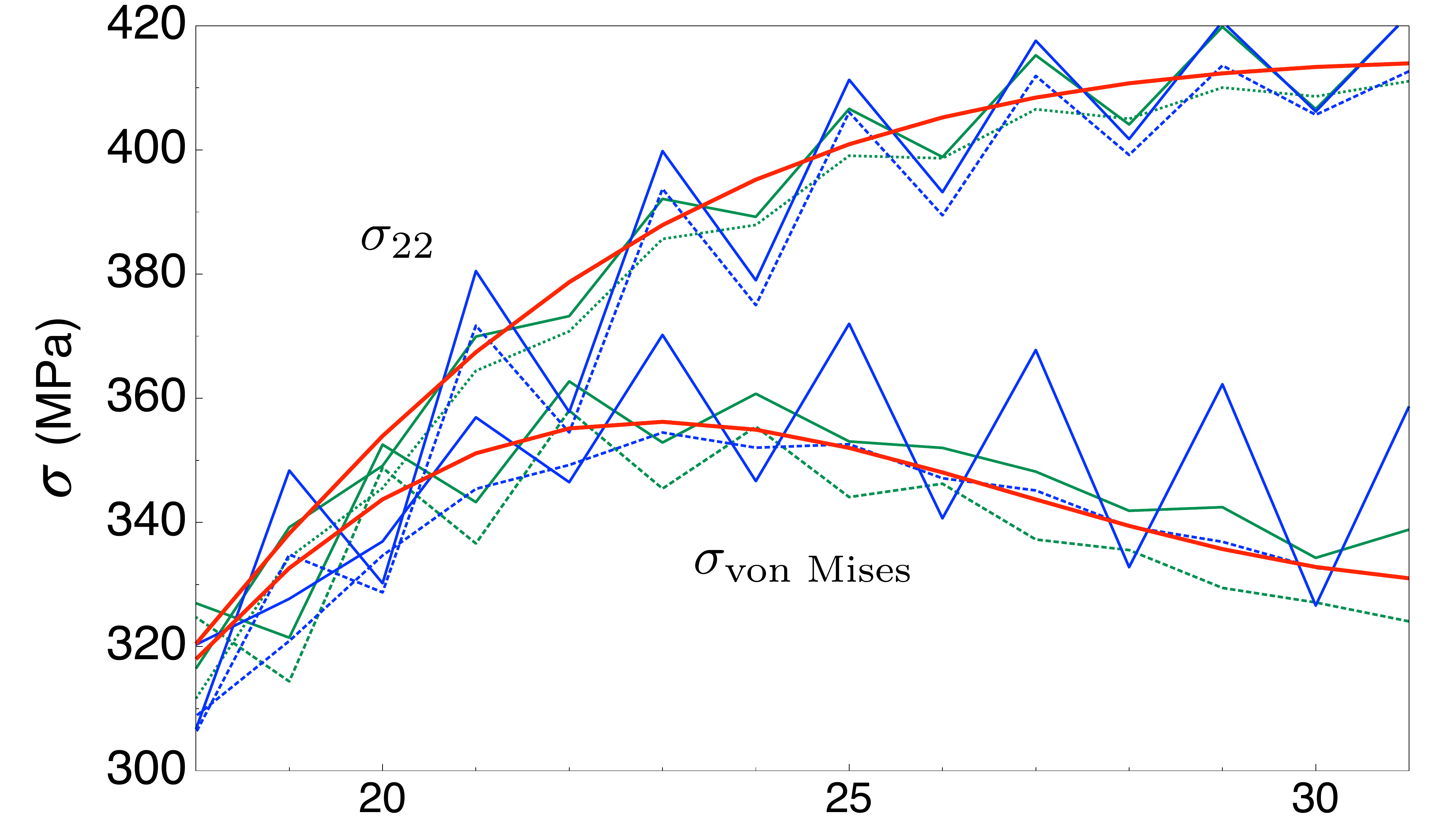}
	\includegraphics[width=0.73\textwidth]{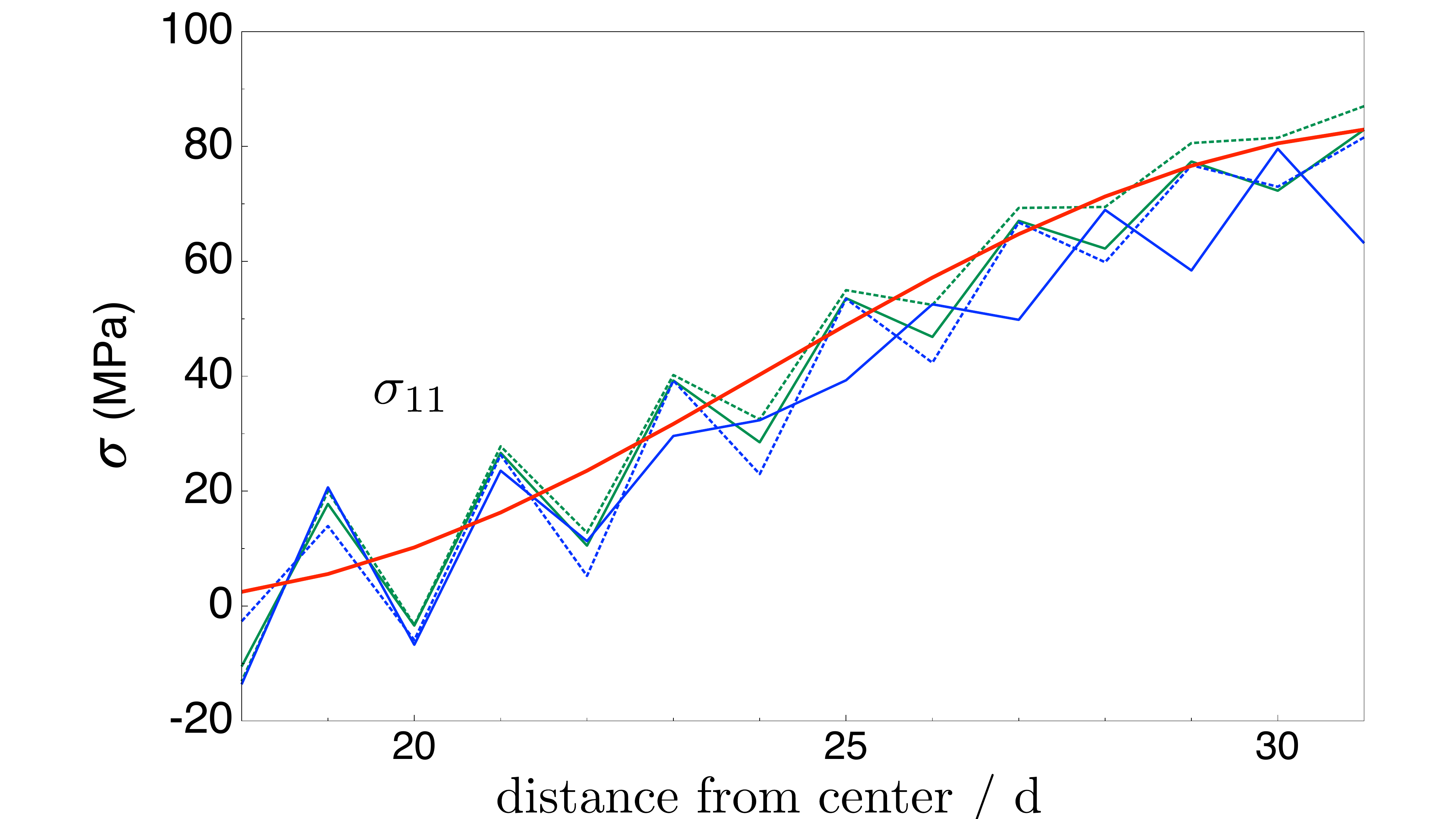}
	\caption{\textcolor{blue}{A cubic void within an isotropic medium subjected to an hydrostatic pressure $P=300$ MPa : internal stresses $\sigma_{11}$,  $\sigma_{22}$ and von Mises stress (MPa) along the line parallel to direction $(100)$ and passing through the center of the void. The numerical results obtained with the tetrahedral stencil (red lines) are compared to the results obtained with the rotated scheme (green lines) and with the Moulinec-Suquet scheme (blue lines); the cubic cavity consists in $31^3$ voxels; full lines represent results obtained with a unit cell containing $64^3$ voxels and doted lines those obtained with a unit cell containing $65^3$ voxels; the middle and bottom panels show zooms on the zones marked 1 and 2, respectively, in the top panel.}
}
\label{fig:cubic_void_1_middle_line}
\end{figure}
\FloatBarrier

Next, we display in Fig.~(\ref{fig:cubic_void_1_middle_line}) the stress components $\sigma_{11}$, $\sigma_{22}$ and the von Mises stress along a line parallel to direction $(100)$ and passing through the center of the cubic void. Results obtained with the tetrahedral stencil are compared to those obtained with the Moulinec-Suquet and rotated schemes.  For the last two methods, computations have also been performed with a computational box whose linear dimensions along the cartesian axis are odd, namely a box that consists into $65^3$ voxels, for which no regularisation recipes are needed. For all the computations, the stiffness tensor $\tenq{\lambda}^0$ of the reference homogeneous medium has been chosen equal $80\%$ of that of the matrix. The results reported for the tetrahedral and rotated stencils are those obtained when the precision $\epsilon$ defined in (\ref{eq:tolerance}) has reached $10^{-10}$, whereas, because of the slow convergence, those obtained with the Moulinec-Suquet method corresponds to a precision equal to $10^{-4}$. As seen in Fig.~(\ref{fig:cubic_void_1_middle_line}), the results obtained with the tetrahedral stencil are exempt of any artefact, whereas significant oscillations are observed with the Moulinec-Suquet method and also with the rotated scheme, whether the linear dimensions of the computational box are even or odd.

Clearly, the recipe proposed by Willot to correct for the lack of definition of the Green function when the box dimensions are even is unsatisfactory, as it leaves the Green function singular at q-points where it should be analytic, which in turn generates oscillations in real space. As seen in Fig.~(\ref{fig:cubic_void_1_middle_line}), these oscillations are hardly suppressed when the box dimensions are odd, even though the Green function in that case is defined for any finite q-vector. The reason why these oscillations are always present, whether the linear dimensions are even or odd, is explained in \ref{sec:appendixB}.

\begin{figure}[h]
	\centering
	\includegraphics[width=0.49\textwidth]{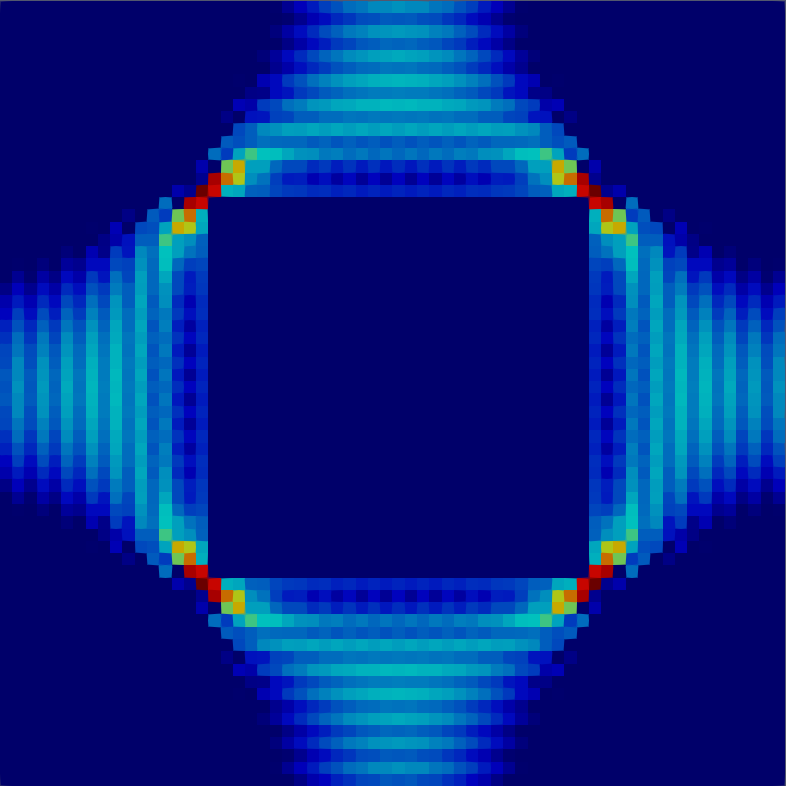} 
	\hspace*{1mm}
	\includegraphics[width=0.49\textwidth]{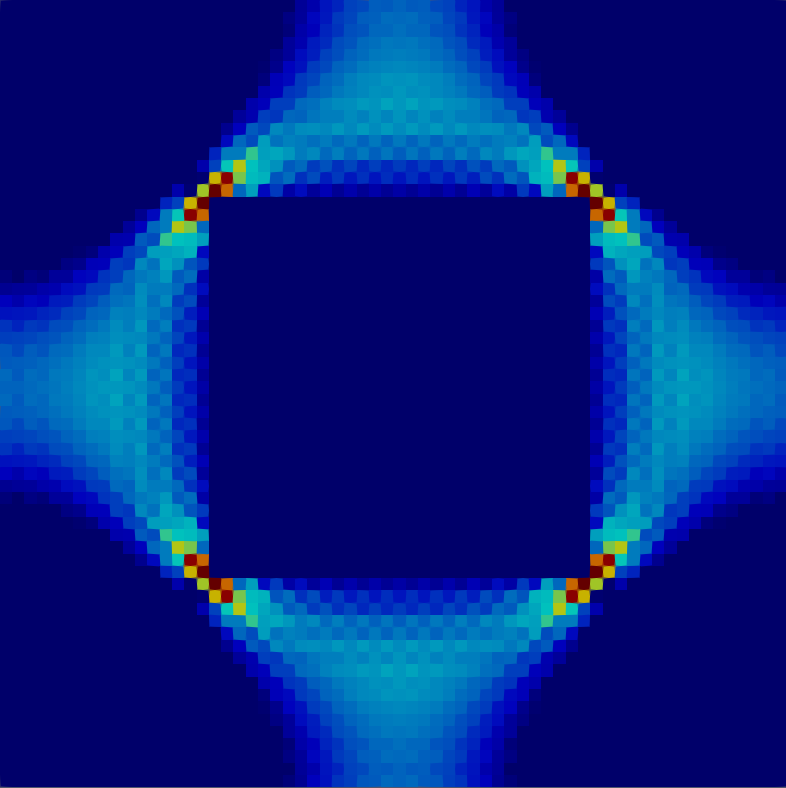}\\
	\vspace*{3mm}
	\includegraphics[width=0.49\textwidth]{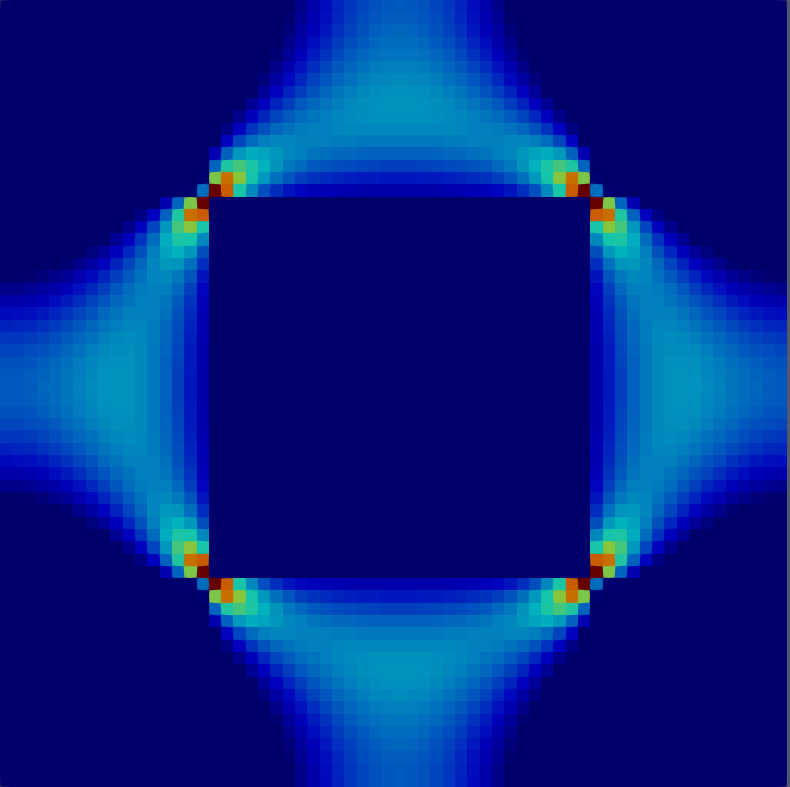}\\
        \vspace*{3mm}
        \includegraphics[width=0.35\textwidth]{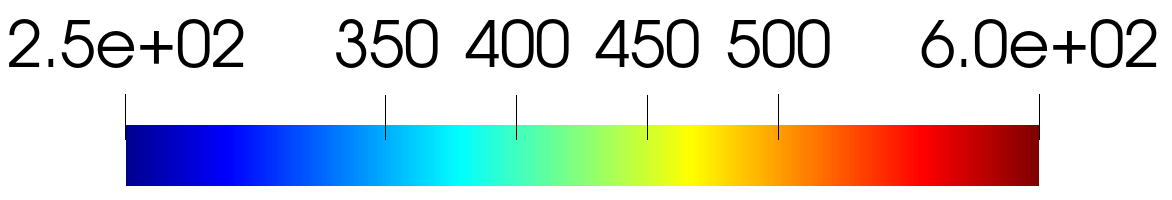}
	\caption{\textcolor{blue}{A cubic void within an isotropic medium subjected to an hydrostatic pressure $P=300$ MPa : von Mises stress maps along the (001) plane passing through the center of computational box; comparaison between the results obtained with the Moulinec-Suquet method (top left panel), the rotated scheme (top right panel) and the tetrahedral solver (bottom panel). The computational box is discretized into $64^3$ cells and the cubic void consists in $31^3$ voxels .}}
\label{fig:cubic_void_von_Mises_maps}
\end{figure}

Finally, we show in Fig.~(\ref{fig:cubic_void_von_Mises_maps}) the map of the von Mises stress obtained with the tetrahedral stencil and, for comparaison, those obtained with Moulinec-Suquet and Willot methods. We clearly observed oscillations with the Moulinec-Suquet method, as well as checkerboard effects with the rotated stencil, whereas the stress map obtained with the tetrahedral stencil is again exempt of any artefact. 
\FloatBarrier 
%
%%%%%%%%%%%%%%%%%%%%%%%%%%%%%%%%%%%%%%%%%%%%%%%%%%%%%%%%%%%%%%%%%%%%%%%%%%%%%%%%%%%%%%%%%%%%%
\subsection{A hard spherical inclusion within a soft matrix}\label{sec:hard_inclusion}
%%%%%%%%%%%%%%%%%%%%%%%%%%%%%%%%%%%%%%%%%%%%%%%%%%%%%%%%%%%%%%%%%%%%%%%%%%%%%%%%%%%%%%%%%%%%%

We consider here a situation opposite to the previous one, i.e. a soft matrix reinforced by hard precipitates. The microstructure consists in a single glass inclusion, with Young modulus and Poisson coefficient given by $E=72$ GPa and $\nu=0.22$, embedded into a polyamide matrix characterised by $E=2.1$ GPa and $\nu=0.3$. The periodic medium consists in $128^3$ cells and the radius of the spherical inclusion is, in cell units, R=40. The overall system is subjected to the uniaxial external strain $\bar \epsilon_{11}=0.05$.

\begin{figure}[ht]
	\centering
	\includegraphics[width=0.6\textwidth]{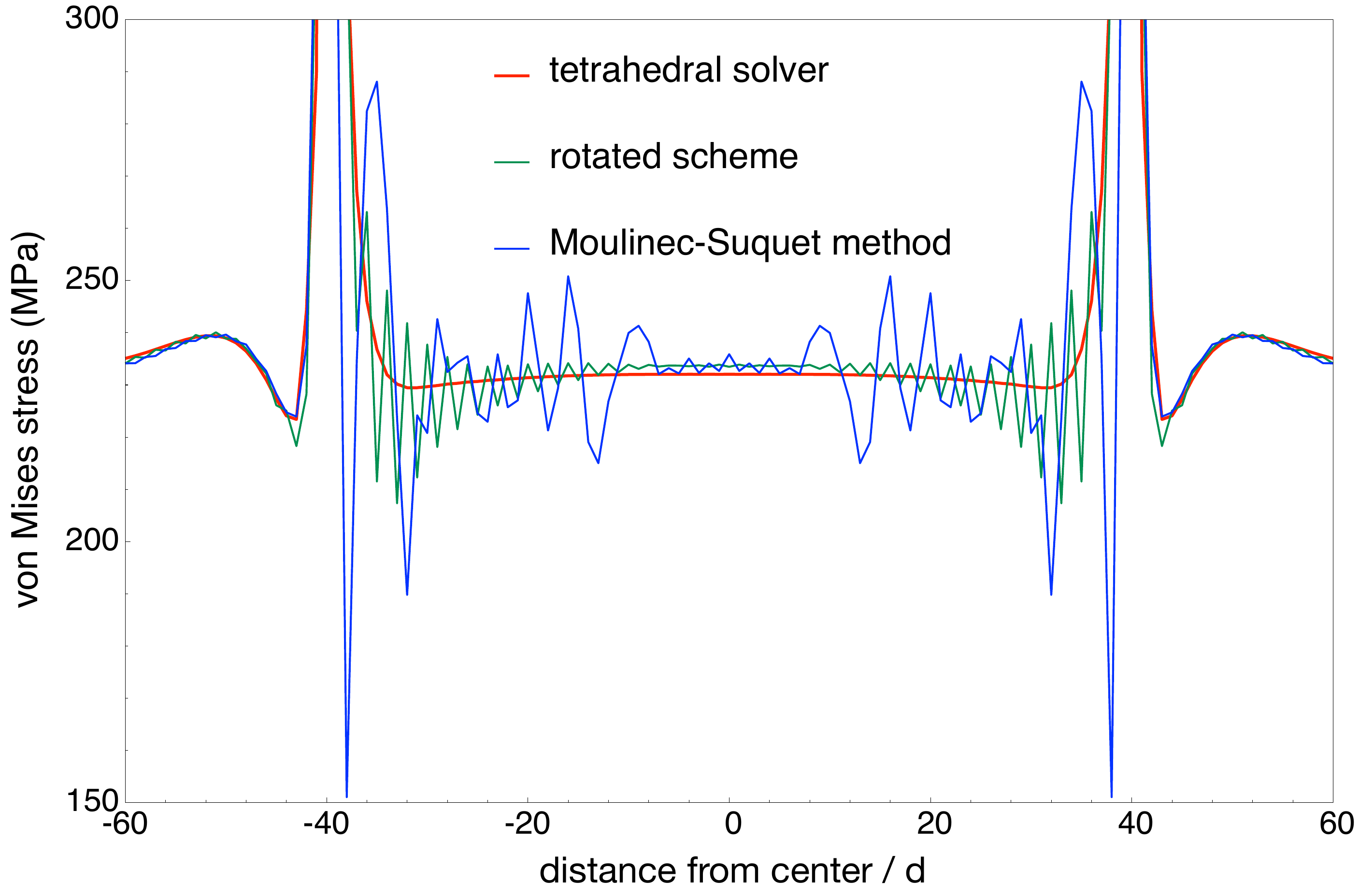}
	\caption{\textcolor{blue}{A hard spherical inclusion within a soft isotropic matrix subjected to the uniaxial external strain $\bar \epsilon_{11}=0.05$ : internal von Mises stress (MPa) along the line parallel to direction $(100)$ and passing through the center of the inclusion. Numerical results obtained with the tetrahedral stencil (red lines) are compared to the results obtained with the rotated scheme (green lines) and with the Moulinec-Suquet scheme (blue lines)}}
\label{fig:spherical_inclusion_von_Mises_along_100}
\end{figure}

\begin{figure}[h]
	\centering
	\includegraphics[width=0.49\textwidth]{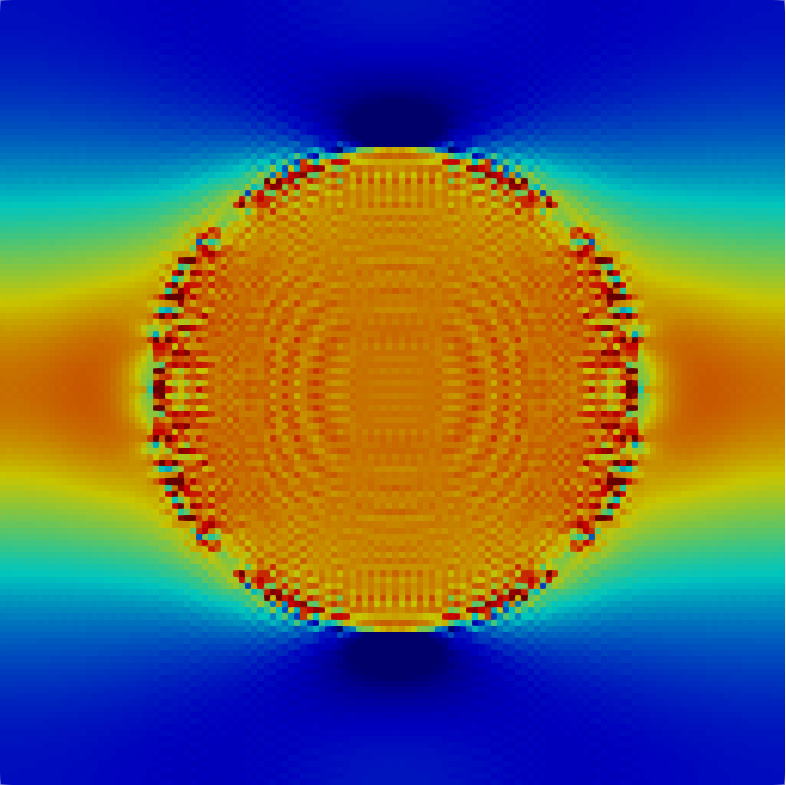} 
	\hspace*{1mm}
	\includegraphics[width=0.49\textwidth]{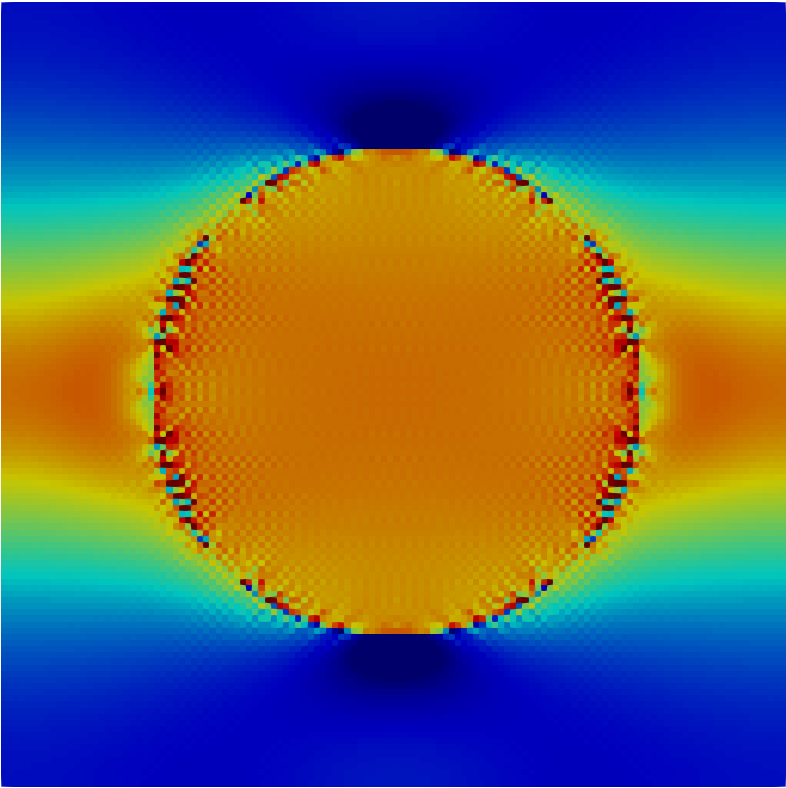}\\
	\vspace*{3mm}
	\includegraphics[width=0.49\textwidth]{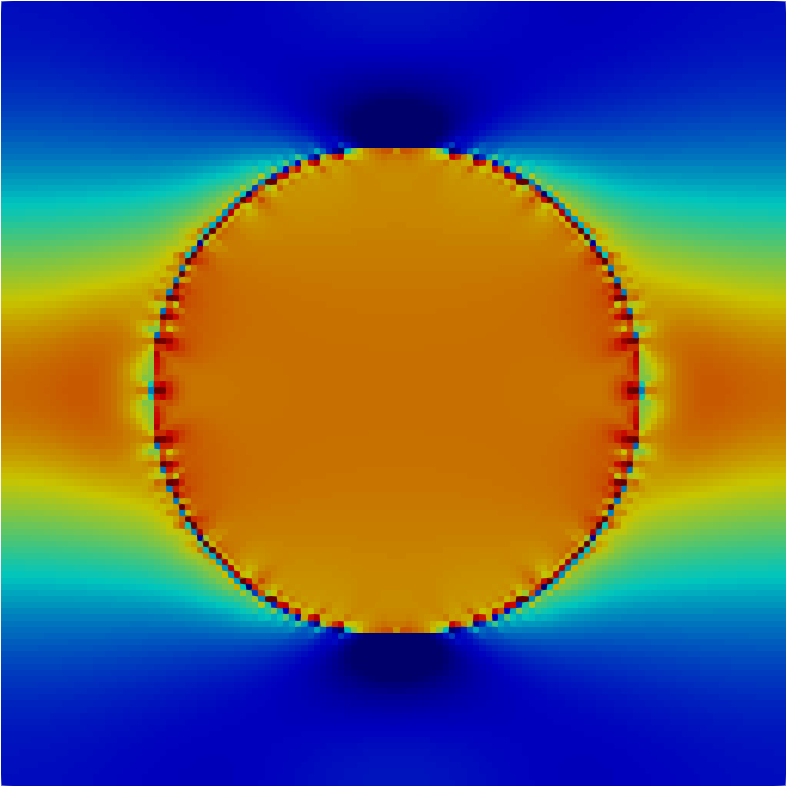}\\
        \vspace*{3mm}
        \includegraphics[width=0.35\textwidth]{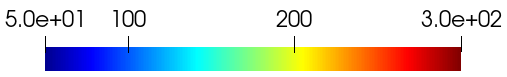}
	\caption{\textcolor{blue}{A cubic void within an isotropic medium subjected to an hydrostatic pressure $P=300$ MPa : von Mises stress maps along the (001) plane that goes through the center of computational box; comparaison between the results obtained with the Moulinec-Suquet solver (top left panel), the rotated scheme (top, right panel) and the tetrahedral solver (bottom panel). The computational box consists into $64^3$ voxels.}}
\label{fig:spherical_inclusion_von_Mises_maps}
\end{figure}

We show in Fig.~(\ref{fig:spherical_inclusion_von_Mises_along_100}) results concerning the von Mises stress obtained with the tetrahedral stencil and, for comparison, the stresses obtained with Moulinec-Suquet and rotated schemes. The stiffness tensor $\tenq{\lambda}^0$ of the reference medium has been  set to $\tenq{\lambda}^0=0.5\,(\tenq{\lambda}^m+\tenq{\lambda}^p)$ where $\tenq{\lambda}^m$ and  $\tenq{\lambda}^p$ are the stiffness matrices of matrix and precipitate, respectively. The results reported for the tetrahedral and rotated stencils are those obtained when the precision $\epsilon$ defined in (\ref{eq:tolerance}) has reached $10^{-10}$, whereas, because of the slow convergence, those obtained with the Moulinec-Suquet method corresponds to a precision equal to $10^{-8}$. Strong oscillations within the inclusion are observed with the Moulinec-Suquet method, which are hardly reduced with the rotated scheme, whereas the results obtained with the tetrahedral stencil are again exempt of any artefacts. The same conclusion are obtained when we compare 2D maps, as shown in Fig.~(\ref{fig:spherical_inclusion_von_Mises_maps}), where we observe ringing artefacts with the Moulinec-Suquet method and checkerboarding with the rotated scheme, whereas the tetrahedral stencil map is free of these artefacts.
\FloatBarrier 
%
%%%%%%%%%%%%%%%%%%%%%%%%%%%%%%%%%%%%%%%%%%%%%%%%%%%%%%%%%%%%%%%%%%%%%%%%%%%%%%%%%%%%%%%%%%%%%
\subsection{Rate of convergence vs contrast}\label{sec:rate_of_convergence}
%%%%%%%%%%%%%%%%%%%%%%%%%%%%%%%%%%%%%%%%%%%%%%%%%%%%%%%%%%%%%%%%%%%%%%%%%%%%%%%%%%%%%%%%%%%%%

Finally, we present a brief analysis of the convergence rate of the tetrahedral solver as a function of the elastic contrast. Our aim is to test the performance of the tetrahedral stencil, associated to the basic iterative scheme, in presence of large contrats. 

The microstructure consists in a single spherical inclusion whose radius, in cell units, is $R=17.6$, embedded into a periodic medium which consists in $44^3$ cells. Matrix and inclusion are elastically isotropic with the same Poisson coefficient $\nu=0.125$ and a large range of elastic contrasts $C=k_{inclusion} / k_{matrix}$ between the bulk moduli of inclusion and matrix is considered. 

\begin{figure}[ht]
	\centering
	\includegraphics[width=0.9\textwidth]{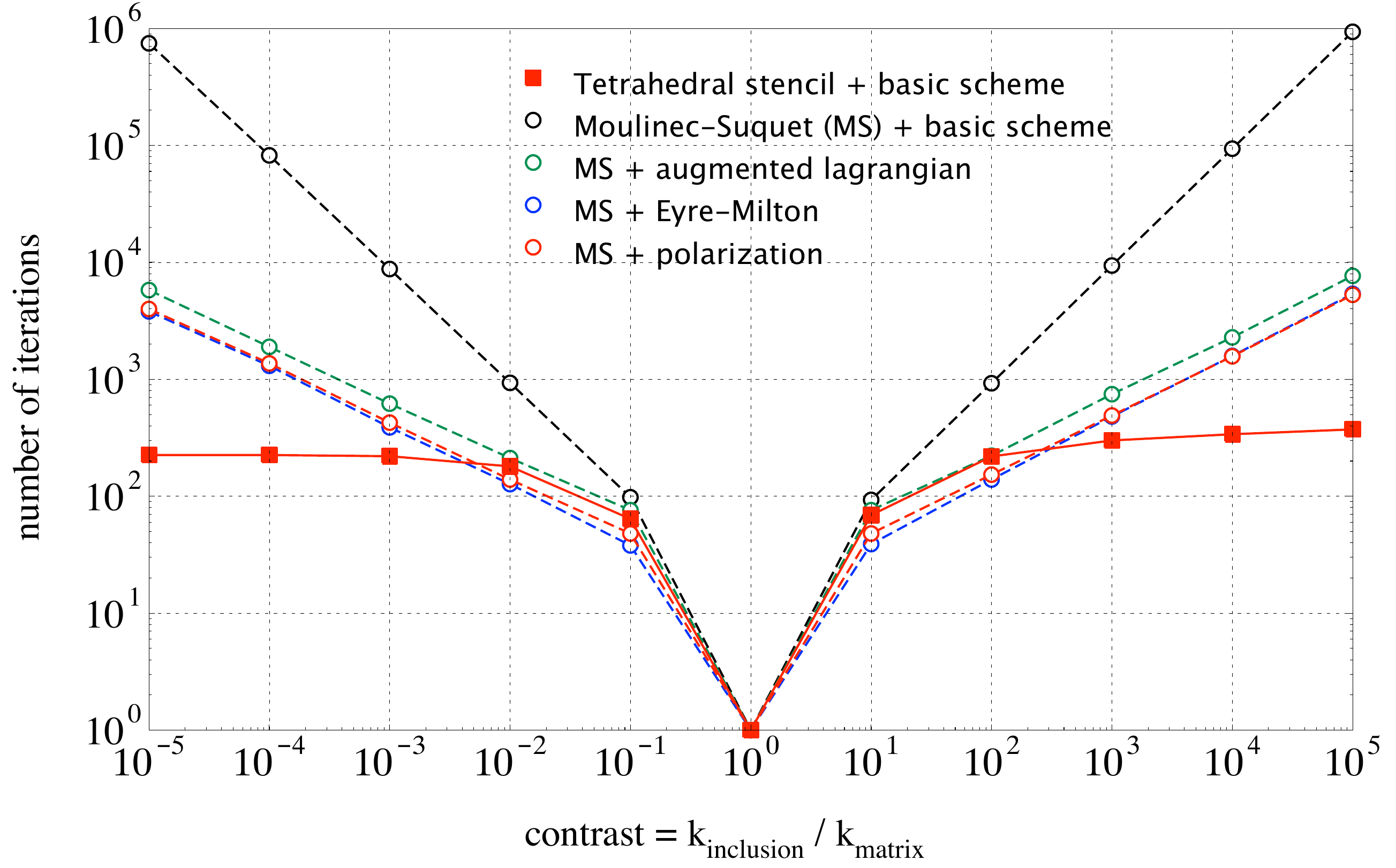}
	\caption{\textcolor{blue}{Number of iterations at convergence (tolerance = $10^{-10}$) for different contrasts between inclusion and matrix bulk moduli. Results from the tetrahedral stencil and the Moulinec-Suquet (MS) basic scheme have been obtained with a reference medium that has the same Poisson coefficient as the inclusion and matrix and a bulk modulus given by $k_0=\alpha\,k_{matrix}+(1-\alpha)\,k_{inclusion}$ with $\alpha=0.55$ for $C<1$ and $\alpha=0.45$ for $C>1$. Results for the MS method accelerated by the augmented lagrangian scheme, Eyre-Milton scheme and the Monchiet-Bonnet polarization method have been taken from Ref. \cite{moulinec2014comparison}.}}
\label{fig:comparaison_convergence_rates}
\end{figure}

We report in Fig~(\ref{fig:comparaison_convergence_rates}) the number of iterates needed to reach equilibrium as a function of the contrast $C$. Precision on equilibrium has been fixed to $10^{-10}$. Results obtained with the tetrahedral stencil are compared to those obtained with the original Moulinec-Suquet method and three accelarated schemes that have been proposed to overcome the slow convergence observed with the Moulinec-Suquet basic scheme, namely the augmented lagrangian method, the Eyre-Milton approach and the polarization-based method proposed by Monchiet and Bonnet. The results concerning the last three methods have been taken from Ref. \cite{moulinec2014comparison}. For the tetrahedral stencil and Moulinec-Suquet method, we used an isotropic reference medium that has the same Poisson coefficient as the inclusion and matrix and a bulk modulus given by $k_0=\alpha\,k_{matrix}+(1-\alpha)\,k_{inclusion}$ with $\alpha=0.55$ for $C<1$ and $\alpha=0.45$ for $C> 1$. 
\FloatBarrier 

The main conclusion is that, whereas the number of iterates increases linearly with $C$ (or with $1/C$ if $C < 1$) for the Moulinec-Suquet method and displays square-root dependancies for the three accelerated iterative schemes, the tetrahedral stencil converges rapidly and, most importantly, requires a number of iterates that essentially does not depend on the contrast. Therefore, microstructures with voids or infinitely rigid inclusions can be handled as easily as situations where the contrast is finite, even though the iterative scheme used here follows the simple basic scheme. 

%%%%%%%%%%%%%%%%%%%%%%%%%%%%%%%%%%%%%%%%%%%%%%%%%%%%%%%%%%%%%%%%%%%%%%%%%%%%%%%%%%%%%%%%%%%%%
\section{Conclusions}\label{sec:conclusion}
%%%%%%%%%%%%%%%%%%%%%%%%%%%%%%%%%%%%%%%%%%%%%%%%%%%%%%%%%%%%%%%%%%%%%%%%%%%%%%%%%%%%%%%%%%%%%
%
We have introduced a new FFT-based method for the computation of the macroscopic mechanical properties of heterogeneous microstructures. The main characters and advantages of the method are the following:
\begin{itemize} 
\item the discretization scheme, which is based on a tetrahedral stencil, has been chosen in such a way that mechanical equilibrium is mathematically stable, i.e. automatically free from any unphysical oscillations or checkerboard effects, such as those observed with the original Moulinec-Suquet method or the discrete rotated scheme;

\item all the components (diagonal and non-diagonal) of a tensorial quantity live on the same location, which avoids the difficulty that arises when the discretization of a constitutive laws links together tensor components that do not sit on the same location;

\item  convergence to equilibrium is fast and requires a number of iterates that, for large contrasts, is essentially independent on the contrast, so that infinite contrast is obtained  without additional computational cost;

\item the method uses the same simple fixed point algorithm as the original basis scheme and, therefore, is easy to implement.

Finally, we mention that the tetrahedral stencil proposed here may be transposed to other fields, such as electrokinetics or multiferroics.

\end{itemize}
%

%%%%%%%%%%%%%%%%%%%%%%%%%%%%%%%%%%%%%%%%%%%%%%%%%%%%%%%%%%%%%%%%%%%%%%%%%%%%%%%%%%%%%%%%%%%%%
%      ANNEXES
%%%%%%%%%%%%%%%%%%%%%%%%%%%%%%%%%%%%%%%%%%%%%%%%%%%%%%%%%%%%%%%%%%%%%%%%%%%%%%%%%%%%%%%%%%%%%

\appendix
%
%%%%%%%%%%%%%%%%%%%%%%%%%%%%%%%%%%%%%%%%%%%%%%%%%%%%%%%%%%%%%%%%%%%%%%%%%%%%%%%%%%%%%%%%%%%%%
\section{Derivative of discrete elastic energy $E^{(k)}$ with respect to Fourier component  $\underline{u}^{(k)}(\underline{q})$ of the displacement field. } \label{sec:appendixA}
%%%%%%%%%%%%%%%%%%%%%%%%%%%%%%%%%%%%%%%%%%%%%%%%%%%%%%%%%%%%%%%%%%%%%%%%%%%%%%%%%%%%%%%%%%%%%
%
In this section, we demonstrate the relationship given in Eq. (\ref{eq:delta_Energy_k_wt_u_Q}). This is done simply as follows. Using the links between the Fourier components of the fluctuating strains $\deps^{(Kk)}(\Q)$ and the displacement fields $\dep^{(k)}(\Q)$ given in Eq.~(\ref{eq:deps_Q}), we have
\begin{equation}
\label{eq:A1}
\frac{\partial E^{(k)} }{\partial u^{(k)}_i(\Q)} = \sum_{K=A,B }\, \frac{\partial E^{(k)} }{\partial \epsilon^{(Kk)}_{ij}(\Q)} \,D^{(Kk)}_j(\Q).
\end{equation}
Using the inverse Fourier transform (\ref{eq:IFT_1_A}) applied to $\eps^{(Kk)}(\Q)$ and the chain rule, the first term in the right-hand side may be written as
\begin{equation}
\label{eq:A2}
\frac{\partial E^{(k)} }{\partial  \epsilon^{(Kk)}_{ij}(\Q)} = \sum_{\R \in K }\, \frac{\partial E^{(k)} }{\partial \epsilon^{(Kk)}_{ij}(\R)} \, \exp(i\Q.\R).
\end{equation}
On the other hand, the stress $\sigma_{ij}^{(Kk)}(\R)$ is, by definition, related to the total energy $E^{(k)}$ given in Eq.~(\ref{eq:Elast_energy_k}) by 
\begin{equation}
\label{eq:A3}
\sigma_{ij}^{(Kk)}(\R) =  \frac{1}{d^3}\, \frac{\partial E^{(k)} }{\partial \epsilon^{(Kk)}_{ij}(\R)}.
\end{equation}
Altogether, Eqs.~(\ref{eq:A1}-\ref{eq:A3}) lead to
\begin{equation}
 \frac{\partial E^{(k)} }{\partial u^{(k)}_i(\Q)} = N d^3 \sum_{K=A,B}   \overline{\sigma_{ij}^{(Kk)}(\Q)}\,{D_{j}^{(Kk)}(\Q)},
\end{equation}
where $\stress^{(Kk)}(\Q)$ is the Fourier transform of the stress $\stress^{(Kk)}(\R)$.
%
%%%%%%%%%%%%%%%%%%%%%%%%%%%%%%%%%%%%%%%%%%%%%%%%%%%%%%%%%%%%%%%%%%%%%%%%%%%%%%%%%%%%%%%%%%%%%
\section{About the non-analyticity of the strain Green function within the rotated scheme} \label{sec:appendixB}
%%%%%%%%%%%%%%%%%%%%%%%%%%%%%%%%%%%%%%%%%%%%%%%%%%%%%%%%%%%%%%%%%%%%%%%%%%%%%%%%%%%%%%%%%%%%%
%

Within the rotated scheme, the displacement and, consequently, strain Green functions are undefined at some finite q-points in Fourier space. These singularities generate mathematical artefacts in the real space strain and stress fields that cannot be suppressed by any recipe, as we explain here.

The rotated scheme may be summarised as follows. The discretization scheme involves a SC grid for the strain field and a staggered SC grid for the displacement field. The Fourier transforms of the finite difference operators, that link the fluctuating strain and displacement fields through $\deps(\Q)= \tens{D}(\Q) \otimes_s \dep(\Q))$, are given by \footnote{If the displacement and strain fields are referred to the same origin in real space, as in Fig.~(\ref{fig:tetrahedral_stencil}), each component $D_i(\Q)$ should be multiplied by $\exp(-i(q_1+q_2+q_3)/2)$. However, this prefactor has no impact on the computed strain field, since, in the inverse of the displacement Green function $\tend{\Omega}(\Q)$ given in (\ref{eq:inverse_displacement_Green_function_rotated}) and in the strain Green function $\tenq{G}(\Q)$ introduced below in (\ref{eq:stain_Lippmann_Schwinger_rotated}) and (\ref{eq:strain_Green_function_rotated}), the components $D_i(\Q)$ appear always by product of the form $D_i(\Q)\overline{D_j(\Q)}$. }
\begin{equation}
\label{eq:Dq_rotated}
\begin{split}
D_1(\Q) &= \frac14 \{ \exp(iq_1) - 1 \} \{ \exp(iq_2) + 1 \} \{ \exp(iq_3) + 1 \},  \\ 
D_2(\Q) &= \frac14 \{ \exp(iq_1) + 1 \} \{ \exp(iq_2) - 1 \} \{ \exp(iq_3) + 1 \} , \\ 
D_3(\Q) &= \frac14 \{ \exp(iq_1) + 1 \} \{ \exp(iq_2) + 1 \} \{ \exp(iq_3) - 1 \} . \\ 
\end{split}
\end{equation}
The displacement-based Lippmann-Schwinger equation reads
\begin{equation}
\forall \Q \ne (000)^*, \quad \tend{\Omega}(\Q)^{-1}.\dep(\Q) = - \pol(\Q).\overline{\tens{D}(\Q)},
\label{eq:displacement_Lippmann_Schwinger_rotated}
\end{equation}
where $\pol(\Q) = \stress(\Q) - \tenq{\lambda}^0:\eps(\Q)$ is the polarization tensor and $\stress(\Q)$ the Fourier transform of the stress field which, in a small deformation setting, is linked to the strain  $\eps(\R)=\bar \eps + \deps(\R)$ by $\stress(\R) = \stiff(\R):(\eps(\R) -\eps^0(\R))$. $\bar \eps$ is the spatial average of the strain and the field $\eps^0(\R)$ is a local eigenstrain. The inverse $\tend{\Omega}(\Q)^{-1}$ of the displacement Green function is given by
\begin{equation}
 \tend{\Omega}(\Q)^{-1} = \overline{\tens{D}(\Q)}.\tenq{\lambda}^0. \tens{D}(\Q).
\label{eq:inverse_displacement_Green_function_rotated}
\end{equation}
As noted in \cite{willot2015fourier}, if the numbers of voxels along any two cartesian axis, say $i$ and $j$, are even, all the components of the q-space differential operator $\tens{D}(\Q)$ are equal to zero when components $q_i$ and $q_j$ of the q-vector reach the highest frequency, i.e. when $q_i=q_j=\pi$, which makes the Green function, and therefore the Fourier transforms of displacement, strain and stress fields, undefined at the corresponding q-points. In fact, this non-definitiveness is one of the manifestation of the singular behaviour displayed by the Green function in the neighbourhoods of the three lines that emerge from the point $(\pi\pi\pi)^*$ and that are aligned along the directions $(100)^*, (010)^*$ and $(001)^*$, to which we refer below as the lines $L_{100}$,  $L_{010}$ and  $L_{001}$, respectively.

To illustrate this behaviour, we consider the neighbourhood of the line $L_{100}$. We note $\Q=(q_1,\pi+\epsilon\sin\theta,\pi+\epsilon\cos\theta)^*$ a point within this neighbourhood, with $\epsilon \ge 0$ and $\theta \in [0,2\pi[ $. Note that we introduce continuous notations for $\Q$, as required to analyze analytical properties in q-space. For any fixed $q_1 \ne \pi$ and in the limit $\epsilon \rightarrow 0$, the differential operators $D_i(\Q)$ behave, to the lowest order in $\epsilon$, as
\begin{equation}
\label{eq:Dq_limit}
\begin{split}
D_1(q_1,\pi+\epsilon\sin\theta,\pi+\epsilon\cos\theta) &\sim \frac{\epsilon^2}{2}  \,  \cos(\frac{q_1-\pi}{2}) \, \exp(i\frac{q_1-\pi}{2}) \, \cos \theta \,  \sin \theta, \\ 
D_2(q_1,\pi+\epsilon\sin\theta,\pi+\epsilon\cos\theta) &\sim  \epsilon  \,  \sin(\frac{q_1-\pi}{2}) \, \exp(i\frac{q_1-\pi}{2}) \,  \cos \theta,  \\ 
D_3(q_1,\pi+\epsilon\sin\theta,\pi+\epsilon\cos\theta) &\sim \epsilon  \,  \sin(\frac{q_1-\pi}{2}) \, \exp(i\frac{q_1-\pi}{2}) \,  \sin \theta. \\ 
\end{split}
\end{equation}
This limiting behaviour is reflected in the inverse Green operator $\tend{\Omega}(\Q)^{-1}$ given in Eq.~(\ref{eq:inverse_displacement_Green_function_rotated}) which, in the limit $0 \le  \epsilon \ll \vert  \tan\frac{q_1-\pi}{2} \vert$ and  for $q_1 \ne \pi$, behaves as
\begin{gather*}
 \tend{\Omega}(q_1,\pi+\epsilon\sin\theta,\pi+\epsilon\cos\theta)^{-1}   \sim  \epsilon^2 	 \\
\times
\begin{pmatrix} 
C_{44} \sin^2\frac{q_1-\pi}{2}                                                                             &\frac{\epsilon}{4}  (C_{12}\!+\!C_{44})  \sin(q_1\!-\!\pi) \cos^2\theta \sin\theta  &\frac{\epsilon}{4}  (C_{12}\!+\!C_{44})  \sin(q_1\!-\!\pi) \cos\theta \sin^2\theta    \\
\frac{\epsilon}{4}  (C_{12}\!+\!C_{44})  \sin(q_1 \!-\! \pi) \cos^2\theta \sin\theta  &(C_{11}\cos^2\theta \!+\! C_{44}\sin^2\theta) \sin^2\frac{q_1-\pi}{2}                &(C_{12}\!+\!C_{44}) \sin\theta\cos\theta \sin^2\frac{q_1-\pi}{2}                           \\
\frac{\epsilon}{4}  (C_{12}\!+\!C_{44})  \sin(q_1 \!-\! \pi) \cos\theta \sin^2\theta  &(C_{12}\!+\!C_{44}) \sin\theta\cos\theta \sin^2\frac{q_1-\pi}{2}                        &(C_{11}\sin^2\theta \!+ \!C_{44}\cos^2\theta) \sin^2\frac{q_1-\pi}{2})
\end{pmatrix},
\end{gather*}
where we have considered a cubic elastic anisotropy. As a result, in the limit $0 \le \epsilon \ll \vert  \tan\frac{q_1-\pi}{2} \vert$ and for $q_1 \ne \pi$, the Green function $\tend{\Omega}(\Q)$ behaves as
\begin{gather}
\nonumber
\tend{\Omega}(q_1,\pi+\epsilon\sin\theta,\pi+\epsilon\cos\theta) \sim  \epsilon^{-2} \, \sin^{-2}(\frac{q_1-\pi}{2}) \\
\times
\begin{pmatrix} 
        \omega_{11}(\theta)  &   \epsilon \sin^{-2}\frac{q_1-\pi}{2} \sin(q_1 \!-\! \pi) \omega_{12}(\theta)                         &    \epsilon \sin^{-2}\frac{q_1-\pi}{2} \sin(q_1 \!-\! \pi) \omega_{13}(\theta)                                 \\
           \epsilon \sin^{-2}\frac{q_1-\pi}{2} \sin(q_1 \!-\! \pi) \omega_{12}(\theta)                  &  \omega_{22}(\theta)      &   \omega_{23}(\theta)   \\
         \epsilon \sin^{-2}\frac{q_1-\pi}{2} \sin(q_1 \!-\! \pi) \omega_{13}(\theta)          &  \omega_{32}(\theta)       &   \omega_{33}(\theta) 
\end{pmatrix}.
\label{eq:Green_function_rotated_limit}
\end{gather}
The explicit  $\theta$-dependance of the functions $\omega_{ij}(\theta)$ are irrelevant for the rest of the discussion.

Now, the impact of the limiting behaviour of $\tend{\Omega}(\Q)$ on the strain and stress fields is readily seen through the examination of the strain Green function $\tenq{G}(\Q)$ which, when $\tend{\Omega}(\Q)$ is defined, controls the strain-based Lippmann-Schwinger equation equivalent to (\ref{eq:displacement_Lippmann_Schwinger_rotated}), i.e.
\begin{equation}
\forall \Q \ne (000)^*, \quad  \deps(\Q) =   -\tenq{G}(\Q) : \pol(\Q),       
\label{eq:stain_Lippmann_Schwinger_rotated}
\end{equation}
with
\begin{equation}
G_{ijkl}(\Q) = \frac14 \{ D_i(\Q) \Omega_{jk}(\Q) \overline{D_l(\Q)} + D_j(\Q) \Omega_{ik}(\Q)  \overline{D_l(\Q)} +   D_i(\Q) \Omega_{jl}(\Q)  \overline{D_k(\Q)} + D_j(\Q) \Omega_{il}(\Q)  \overline{D_k(\Q)} \}.    
\label{eq:strain_Green_function_rotated}  
\end{equation}
From Esq.~(\ref{eq:Dq_limit}) and (\ref{eq:Green_function_rotated_limit}), and using Voigt notations, it is easy to see that, for any angle $\theta$ and for $q_1 \ne \pi$, the Green function component $G_{IJ}(q_1,\pi+\epsilon\sin\theta,\pi+\epsilon\cos\theta) \rightarrow 0$ when $\epsilon \rightarrow 0$ except when the pair of indices $IJ$ (or the pair obtained by permutation of $I$ and $J$) belongs to the set $S_{100}$ defined by
\begin{equation}
\label{eq:De}
S_{100} = \{ 22,23,24,33,34,44,55,56,66 \}, 
\end{equation}
where the notation $S_{100}$ is a reminder that we are considering the behaviour of the Green function in the neighbourhood of the q-line $L_{100}$. However, when $IJ$ belongs to $S_{100}$ and for $q_1 \ne \pi$, component $G_{IJ}(q_1,\pi+\epsilon\sin\theta,\pi+\epsilon\cos\theta)$  stays finite when we approach the line. More precisely, within the finite neighbourhood of $L_{100}$ defined by $0 \le \epsilon \ll \vert  \tan\frac{q_1-\pi}{2} \vert$ and $q_1 \ne \pi$, $G_{IJ}(q_1,\pi+\epsilon\sin\theta,\pi+\epsilon\cos\theta)$ becomes independent of $q_1$ and of $\epsilon$, but varies continuously with the angle $\theta$ around the $L_{100}$ line, i.e.
\begin{equation}
\text{if } IJ \in S_{100}, \,\, 0  \le \epsilon \ll \vert  \tan\frac{q_1-\pi}{2} \vert \text{ and } q_1 \ne \pi,  \quad G_{IJ}(q_1,\pi+\epsilon\sin\theta,\pi+\epsilon\cos\theta) \simeq g_{IJ}(\theta),
\label{eq:G_ij_theta}
\end{equation}
where the explicit $\theta$-dependance of the functions $g_{IJ}(\theta)$ are irrelevant. The important consequence of this behaviour is that, when $IJ$ belongs to $S_{100}$, there is no analytic continuation for $G_{IJ}(\Q)$ along the line $L_{100}$, which explained why the non-definitiveness of $\tenq{G}(\Q)$ along $L_{100}$, and the ensuing mathematical artefacts in real space, cannot be solved by any recipe. Conversely, when $IJ$ does not belong to $S_{100}$,  and  except on the point $(\pi\pi\pi)^*$, the Green function component $G_{IJ}(\Q)$ may be analytically continued on the line $L_{100}$ by setting $G_{IJ}(\Q)=0$.

A similar analysis along lines $L_{010}$ and $L_{001}$ leads to the same results as above, replacing $S_{100}$ by $S_{010}$ or $S_{001}$ as appropriate, with $S_{010}$ and $S_{001}$  defined by
\begin{equation}
\nonumber
\begin{split}
S_{010} = \{ 11,13,15,33,35,44,46,55,66 \}, \\
S_{001} = \{ 11,12,16,22,26,44,45,55,66 \}.
\end{split}
\end{equation}
When the couple of indices $IJ$ does not belong to any one of the sets $S_{100}$, $S_{010}$ and $S_{001}$, the Green function component $G_{IJ}(\Q)$ may be analytically continued along the lines $L_{100}$, $L_{010}$ and $L_{001}$ using the recipe $G_{IJ}(\Q)=0$ (except at the point $q=(\pi\pi\pi)^*$, which is discussed in the next paragraph). This is the case for $IJ$ (or $JI$) $=14, 25$ and $36$. However, for any other $IJ$, the corresponding  $G_{IJ}(\Q)$ is non-analytic along at least one of the lines $L_{100}$, $L_{010}$ and $L_{001}$.

The point $q=(\pi\pi\pi)^*$ was excluded from the previous analysis because the mathematical behaviour of the strain Green function in the neighbourhood of this point differs from the one observed above. Indeed, a simple analysis shows that, when we approach $q=(\pi\pi\pi)^*$ along any ray and for any couple $IJ$,  $G_{IJ}(\Q)$ stays finite and equal to a value that depends only on the direction of the ray, i.e.
\begin{equation}
\forall  \,\, IJ,  \quad 0  \le \epsilon \ll 1,\quad G_{IJ}(\pi + \epsilon\sin\theta\cos\phi,\pi + \epsilon\sin\theta\sin\phi,\pi + \epsilon\cos\theta) \simeq h_{IJ}(\theta,\phi)
\label{eq:G_ij_theta_phi}
\end{equation}
where the angles $\theta$ and $\phi$ define the direction of the ray in a spherical coordinate system around $\Q=(\pi\pi\pi)^*$. Therefore, for any $IJ$, $G_{IJ}(\Q)$  cannot be analytically continued to $q=(\pi\pi\pi)^*$.

These non-analytical behaviours have important consequences in real space which are best highlighted by analysing the inverse Fourier transform of the Green function $\tenq{G}(\Q)$, defined by
\begin{equation}
\tenq{G}(\R) = \frac{1}{N} \,\sum_{\Q} \tenq{G}(\Q) \exp(i\Q.\R).
\label{eq:G_ij_r}
\end{equation}
Consider for instance the Green function component $G_{24}(\Q)$ which, according to the previous analysis, is only singular along the line $L_{100}$. When we approach $L_{100}$ along a ray perpendicular to the line, $G_{24}(\Q)$ becomes constant and equal to a value that depends only on the direction of the ray defined by the angle $\theta$ (see Eq.~(\ref{eq:G_ij_theta})). Now, because of this "one-dimensional" non-analytical  behaviour, we make the hypothesis that, in (\ref{eq:G_ij_r}), the contribution of  the neighbourhood of the line $L_{100}$ to the real space Green function $G_{24}(\R)$ will dominate and lead to a very slow decrease of $G_{24}(\R)$ with $\vert\vert \R \vert\vert$. This hypothesis will be supported later, when we will observe that the numerical computation of $G_{24}(\R)$ through the exact inverse Fourier transform (\ref{eq:G_ij_r}) will be in agreement with the result of the ongoing approximate analysis. For now,  we proceed by extending (\ref{eq:G_ij_theta}) to the entire q-space, which leads to the approximation 
\begin{equation}
G_{24}(l_1,l_2,l_3) \approx  \frac{1}{N} \,\exp\{i\pi(l_2+l_3)\} \, \sum_{q_1,Q_2,Q_3} g_{24}(\frac{Q_2}{\vert\vert \tens{Q} \vert\vert},\frac{Q_3}{\vert\vert \tens{Q} \vert\vert})  \exp\{i(q_1l_1+Q_2l_2+Q_3l_3)\},
\label{eq:G24_approxi}
\end{equation}
where the integers $l_1,l_2$ and $l_3$ are the coordinates of the vector $\R$ and where the vector $\Q=(q_1,q_2,q_3)^*$ has been decomposed into $\Q=(q_1,\pi,\pi)^* + \tens{Q}$ with $\tens{Q}=(0,q_2-\pi,q_3-\pi)^*$. Since the vector $\tens{Q}$ is perpendicular to the line $L_{100}$ and univocally defines the angle $\theta$, the $\theta$-dependent function $g_{24}(\theta)$ defined in Eq.~(\ref{eq:G_ij_theta}) has been equivalently replaced by the function $g_{24}(Q_2/\vert\vert \tens{Q} \vert\vert,Q_3/\vert\vert \tens{Q} \vert\vert)$, where $\vert\vert \tens{Q} \vert\vert$ is the norm of $\tens{Q}$. Now, we proceed by considering $\R$ vectors that lie within the plane $(100)$ that contains the origin of real space, i.e. vectors of type $\R=(0,l_2,l_3)$. Since the q-space function $g_{24}(Q_2/\vert\vert \tens{Q} \vert\vert,Q_3/\vert\vert \tens{Q} \vert\vert)$ does not depend on $q_1$, we get
\begin{equation}
G_{24}(0,l_2,l_3) \approx  \frac{1}{N_2N_3} \,\exp\{i\pi(l_2+l_3)\} \, \sum_{Q_2,Q_3} g_{24}(\frac{Q_2}{\vert\vert \tens{Q} \vert\vert},\frac{Q_3}{\vert\vert \tens{Q} \vert\vert})  \exp\{i(Q_2l_2+Q_3l_3)\}.
\end{equation}
Now, we consider the continuous limit of the remaining discrete sums and also extend to infinity the domain of integration, which leads to
\begin{equation}
G_{24}(0,l_2,l_3) \approx  \frac{1}{(2\pi)^2} \,\exp\{i\pi(l_2+l_3)\} \, \int_{Q_2,Q_3} \, dQ_2dQ_3 \,\,  g_{24}(\frac{Q_2}{\vert\vert \tens{Q} \vert\vert},\frac{Q_3}{\vert\vert \tens{Q} \vert\vert})  \exp\{i(Q_2l_2+Q_3l_3)\}.
\label{eq:G_ij_r_continuous}
\end{equation}
We expect these approximations to be valid if $\vert\vert \R \vert\vert$ is sufficiently smaller than the system size $L$, in order to minimize the impact of the periodic boundary conditions which limit the continuous setting of the q-space, and sufficiently larger than $d$, in order to minimize the effects of the finite grid spacing which limits the extension of the q-space domain. We immediately see from Eq.~(\ref{eq:G_ij_r_continuous}) that we have the following property:
\begin{equation}
G_{24}(0,\lambda l_2,\lambda l_3) \approx \frac{1}{\lambda^2} \, \exp\{ i\pi\lambda(l_2+l_3)\} \, G_{24}(0, l_2, l_3)
\end{equation}
where, $l_2$ and $l_3$ being integers, the exponential term is equal to $\pm 1$. We conclude that, in the plane (100) which contains the origin of real space, the Green function $G_{24}(\R)$ oscillates and decreases as $1/\vert\vert \R \vert\vert^2$. 

This slow, oscillating decay will be observed for any $G_{IJ}(\R)$ component if $IJ$ belongs to one (and only one) of the sets $S_{100}$, $S_{010}$ or $S_{001}$ (for $\R$ vectors in the $(100)$, $(010)$ or $(001)$ plane which contains the origin of real space, as appropriate).

If $IJ$ does not belong to any one of the sets $S_{100}$, $S_{010}$ and $S_{001}$, we have seen above that $G_{IJ}(\R)$ may be analytically continued along the lines $L_{100}$, $L_{010}$ and $L_{001}$, except at the point $\Q=(\pi\pi\pi)^*$. A simple analysis in the same spirit as above, with Eq.~(\ref{eq:G_ij_theta}) replaced by Eq.~(\ref{eq:G_ij_theta_phi}), shows that this remaining singularity generates in real space an oscillatory decrease in $1/\vert\vert \R \vert\vert^3$. 

These oscillatory decays are very different from the behaviour we would observe if the q-space Green function $\tenq{G}(\Q)$ was regular for any finite $\Q$. In fact, the only singularity that should be observed is that associated with translational invariance, which therefore lies at the point $\Q=(000)^*$ and which, when it is alone, generates a monotonous decrease in $1 / \vert\vert \R \vert\vert^3$ for all the Green function component $G_{IJ}(\R)$.

We now present numerical results that will confirm the approximate mathematical analysis that has been outlined above. Green functions obtained with the rotated scheme will be compared to those obtained with the tetrahedral stencil. For the latter, we will consider the real part of the sum of the Green functions $\tenq{G}^d(\Q)$ and $\tenq{G}^{nd}(\Q)$ defined in Eq.~(\ref{eq:strain_based_Green_functions_tetrahedral_stencil}). Referring to the mathematical nature of strain Green functions, which is to generate a strain as a response to a given polarization field, and taking that the final strain for the tetrahedral stencil is the arithmetic average defined in Eq.~(\ref{eq:final_strain}), the collapse of the two Lippmann-Schwinger equations given in Eq.~(\ref{eq:Lippmann_Schwinger_strain_1}), in which we identify the two polarization fields with a single given field, indeed leads to the real part of this sum. The elastic tensor $\tenq{\lambda}^0$ is isotropic, with shear modulus and Poisson coefficient given by $\mu=132.30$ GPa and $\nu=0.26$, and the periodic medium is discretized with $L=128$ cells along each cartesian coordinates.

\begin{figure}[h]
	\centering
	\includegraphics[width=1\textwidth]{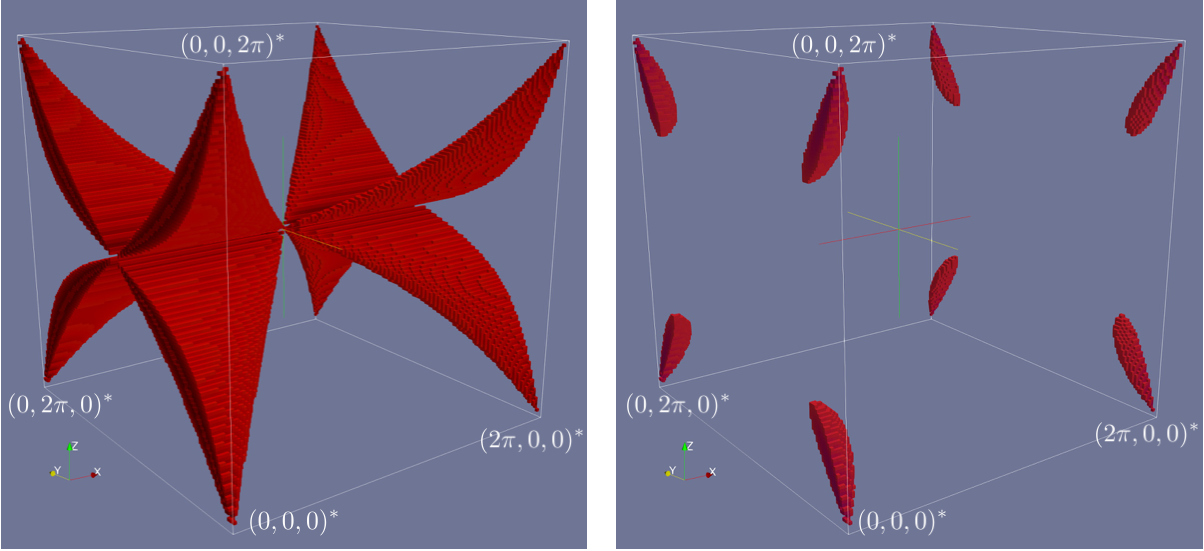} 
	\caption{\textcolor{blue}{Isosurfaces of the modulus of the q-space Green functions $G_{24}(\Q)$ obtained with the rotated scheme (left panel) and with the tetrahedral stencil (right panel) (the isosurfaces have been  defined by $G_{24}(\Q)=0.27 \pm 0.2$, in units of $\mu^{-1}$).}}
\label{fig:G_24_q_rotated_scheme_vs_tetrahedral_stencil}
\end{figure}

\begin{figure}[h]
	\centering
	\hspace*{5mm} $G_{24}(\R)  \times (\vert\vert \R \vert\vert/L)^2$   \hspace*{48mm}  $G_{24}(\R)  \times (\vert\vert \R \vert\vert/L)^3$ \\
	\vspace*{3mm}
	\begin{turn}{90} 
	\hspace*{35mm} Rotated scheme
	\end{turn}
	\hspace*{0.5mm}
	\includegraphics[width=0.48\textwidth]{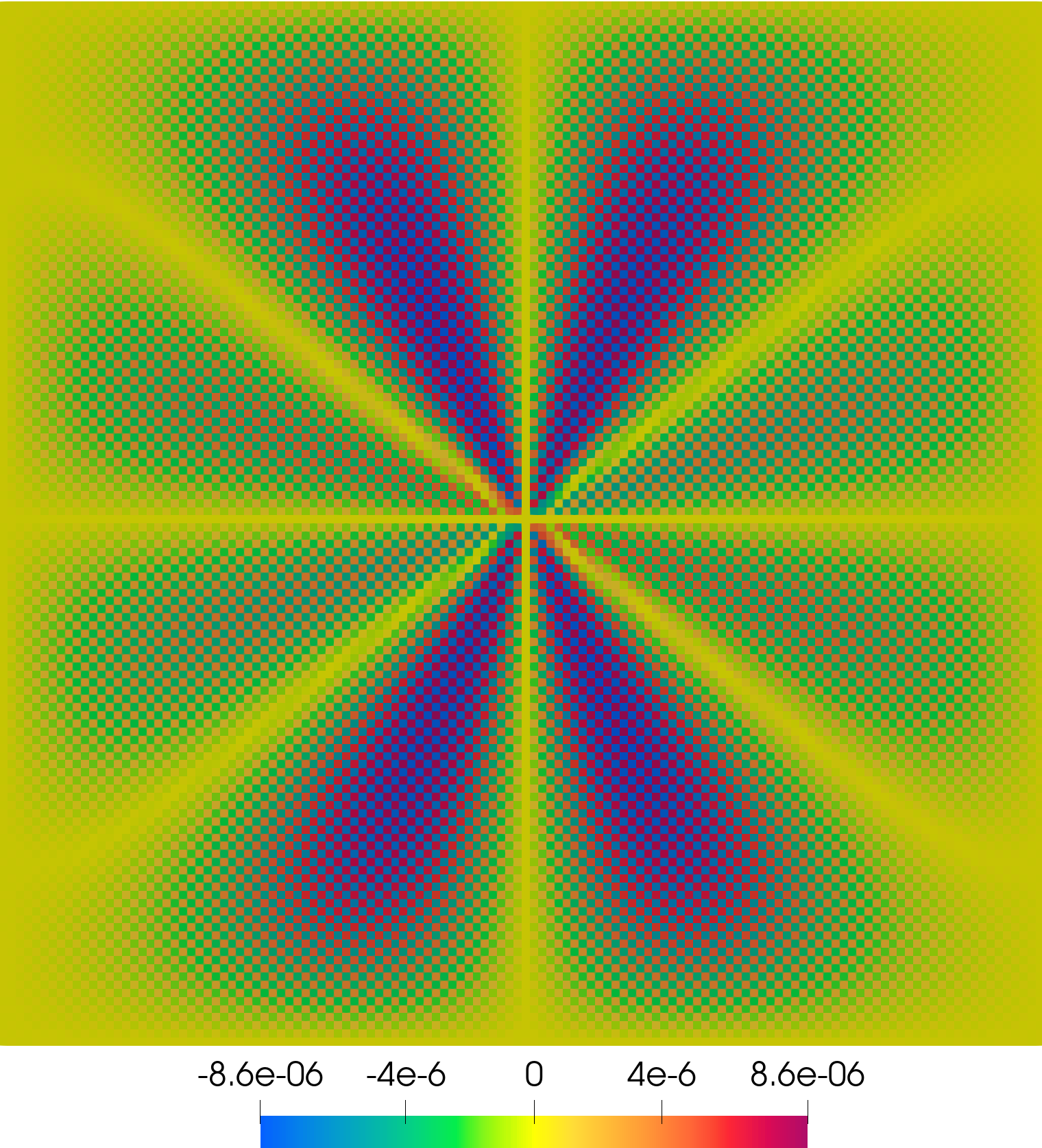} 
	\includegraphics[width=0.48\textwidth]{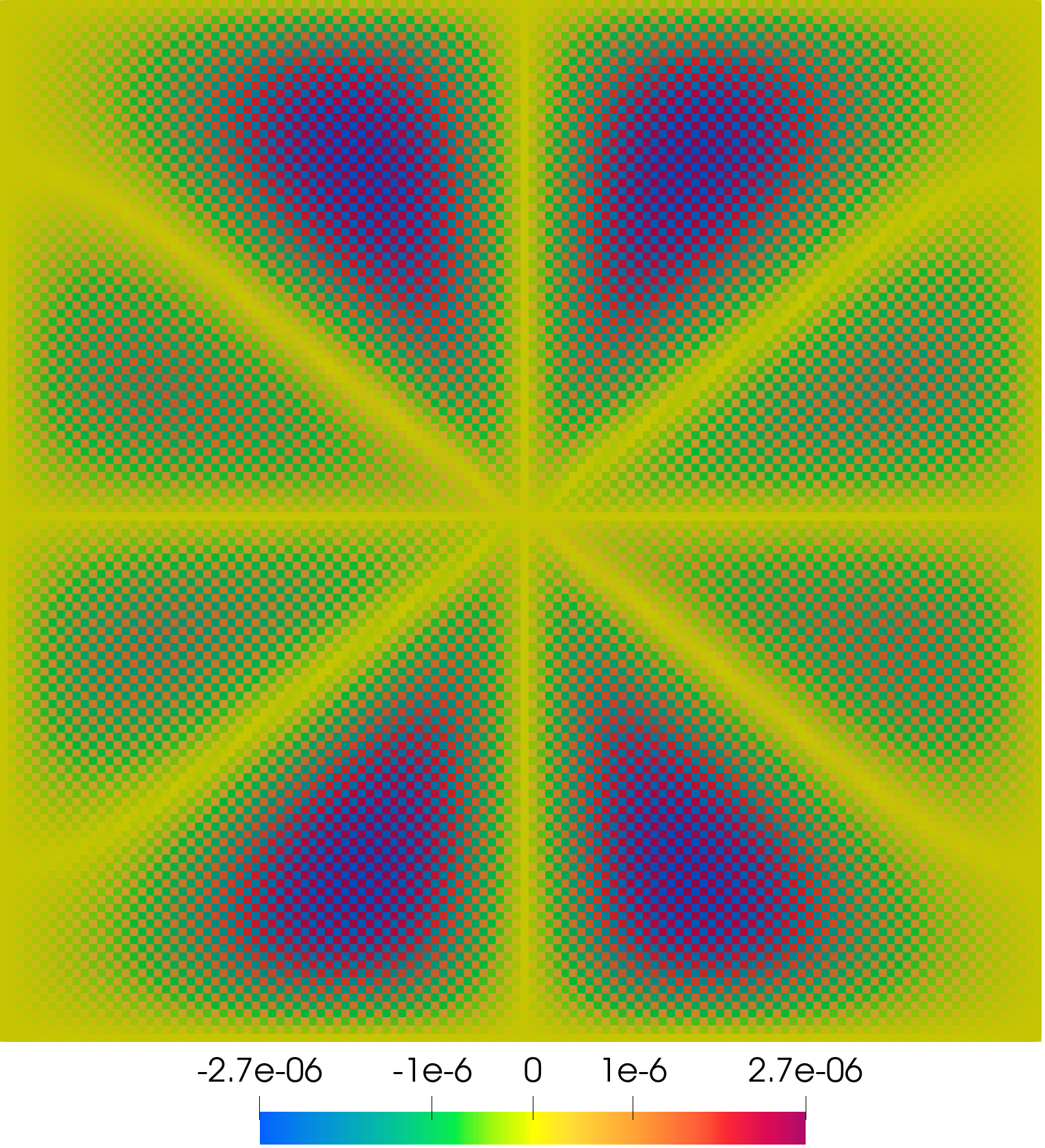}\\
	\vspace*{10mm}
	\begin{turn}{90} 
	\hspace*{35mm} Tetrahedral stencil
	\end{turn}
	\hspace*{1mm}
	\includegraphics[width=0.48\textwidth]{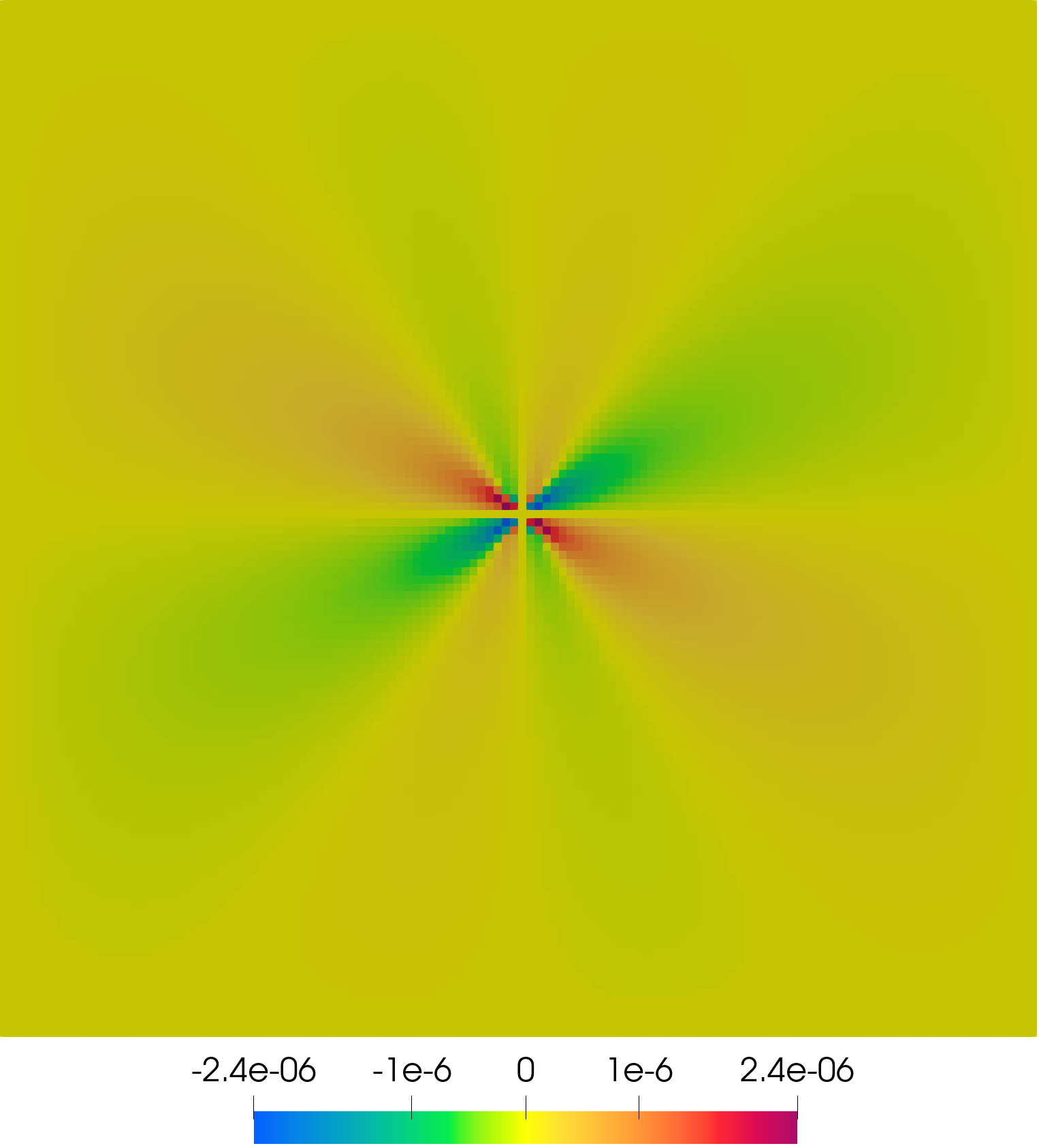} 
	\includegraphics[width=0.48\textwidth]{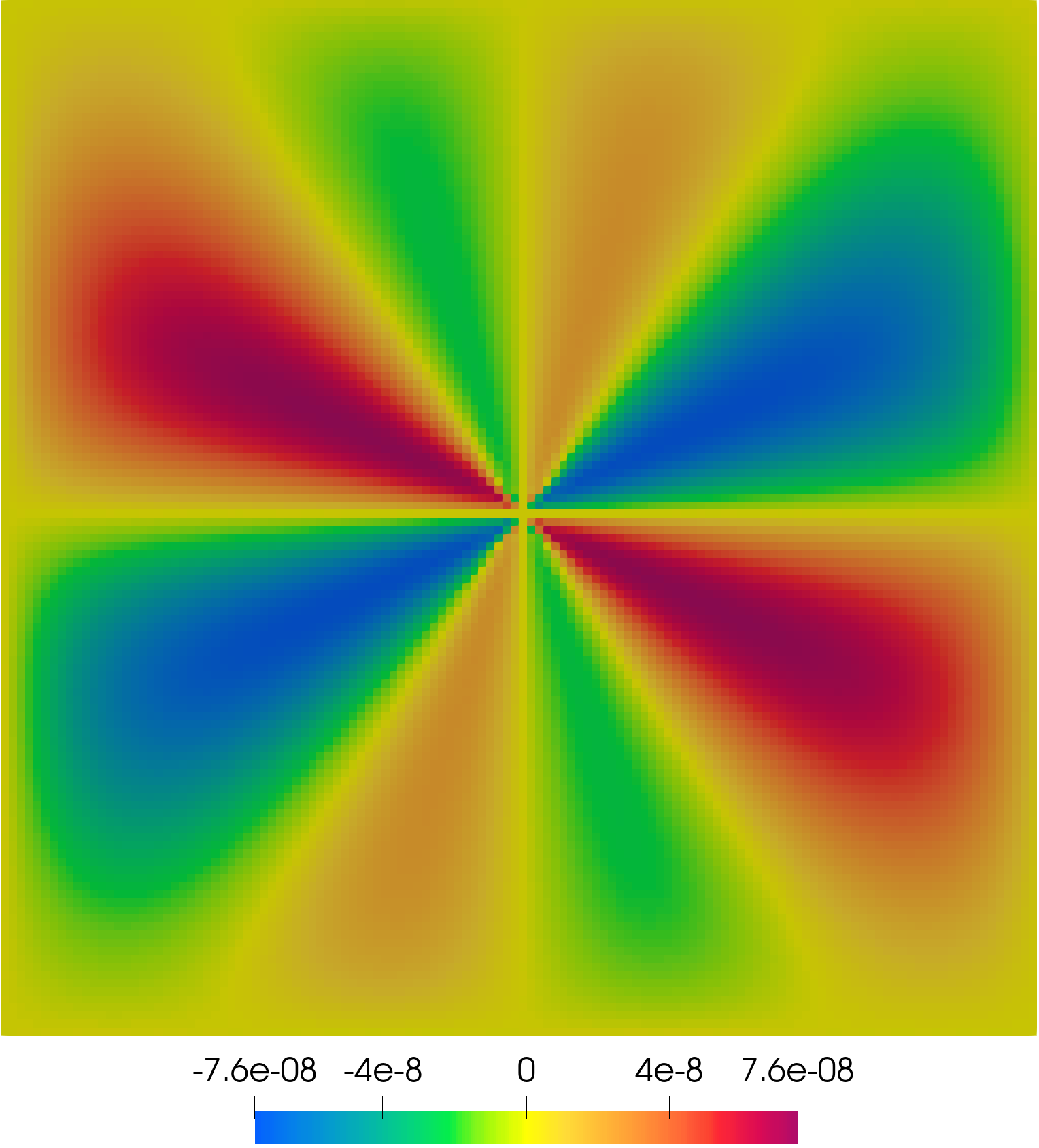}\\
	\caption{\textcolor{blue}{Green functions $G_{24}(\R)$ for $\R=(0,l_2,l_3)$ within the rotated scheme (top row) and the tetrahedral stencil (bottom row) :  scaling by $(\vert\vert \R \vert\vert/L)^2$ (left panels) and by $(\vert\vert \R \vert\vert/L)^3$ (right panels) (units : $\mu^{-1}$); integers $l_2$ and $l_3$ run horizontally and vertically, respectively, from $-L/2$ to $L/2-1$, with $L=128$.}}
\label{fig:G_24_r_rotated_scheme_vs_tetrahedral_stencil}
\end{figure}

First, we show in Fig.~(\ref{fig:G_24_q_rotated_scheme_vs_tetrahedral_stencil}) isosurfaces of the modulus of q-space Green functions $G_{24}(\Q)$ obtained with the rotated and tetrahedral stencils (left and right panels, respectively). For the rotated scheme, we clearly observe the singular behaviour along the line $L_{100}$ discussed above: sufficiently close to the line, $G_{24}(q_1,q_2,q_3)$ does not depend on $q_1$ or on the distance perpendicular to the line. For the tetrahedral stencil, we only observe, as expected, the isolated singularity associated with the point $\Q=(000)^*$, reproduced by periodicity at the height corners of the q-space domain defined in Eq.~(\ref{eq:q-vectors}).
\FloatBarrier
Next, we discuss the real space Green functions $G_{24}(\R)$, obtained numerically by inverse Fourier transform, as expressed in Eq.~(\ref{eq:G_ij_r}). In Fig.~(\ref{fig:G_24_r_rotated_scheme_vs_tetrahedral_stencil}), the Green functions $G_{24}(\R)$, for $\R=(0,l_2,l_3)$, have been scaled by $\vert\vert \R \vert\vert^2$ in the left panels and by $\vert\vert \R \vert\vert^3$ in the right ones. The results obtained with the rotated scheme and those obtained with the tetrahedral stencil are shown in the top and bottom rows, respectively. In line with the analysis presented above, we observe that, with the rotated scheme and if $\vert\vert \R \vert\vert$ is sufficiently smaller than the system size $L$, $G_{24}(\R)$ exhibits scaling by $\vert\vert \R \vert\vert^2$ and not by $\vert\vert \R \vert\vert^3$. Conversely, with the tetrahedral stencil and again if $\vert\vert \R \vert\vert$ is sufficiently smaller than the system size, we clearly observe scaling by $\vert\vert \R \vert\vert^3$ and no scaling by $\vert\vert \R \vert\vert^2$.

Note that the $\vert\vert \R \vert\vert^2$-scaling within the rotated scheme observed here numerically justifies \emph{a posteriori} the hypothesis we made above to write Eq.~(\ref{eq:G24_approxi}). Indeed, it confirms the dominant role of the $L_{100}$ singularity in the behaviour of $G_{24}(\R)$, since in the absence of this one-dimensional singularity, the slow $1/\vert\vert \R \vert\vert^2$  decay would be replaced by a the faster $1\vert\vert \R \vert\vert^3$ oscillatory decay associated to the singularity at $\Q=(\pi\pi\pi)^*$ supplemented by the $1/\vert\vert \R \vert\vert^3$ monotonous decay associated to the $\Q=(000)^*$ singularity.

In brief, within the rotated scheme, a Green function component $G_{IJ}(\R)$ that is singular along one of the lines $L_{100}$, $L_{010}$ or $L_{001}$ generates, in real space, oscillations that decrease as $1/\vert\vert \R \vert\vert^2$, which, in a 3D space, represents a very slow decrease. This is certainly at the source of the long-ranged stress oscillations observed with the rotated scheme.  As shown here, the mathematical origin of these artefacts is the behaviour of the Green function in the neighbourhood of the lines $L_{100}$, $L_{010}$ and $L_{001}$, which excludes any analytic continuation of the Green function along those lines. In numerical applications, finite boxes are used, which discretizes the Fourier space on a set of discrete q-points, as in (\ref{eq:q-vectors}). If the box dimensions are odd, none of these q-points sit on the lines $L_{100}$, $L_{010}$ and $L_{001}$. However, the underlying singular behaviour in the continuous q-space is still there and manifests itself through the discrete q-points that are close enough to the lines. Therefore, there is no recipe to correct for the real space artefacts, whether the linear dimensions of the simulation box are odd or even.

In summary, within the rotated scheme, the strain Green function components $G_{IJ}(\R)$ exhibit unphysical oscillations. For the few that do not belong to any of the sets $S_{100}$, $S_{010}$ and $S_{001}$, these oscillations decrease as $1\vert\vert \R \vert\vert^3$. However, for the many that belong to one and only one of the sets $S_{100}$, $S_{010}$ or $S_{001}$, these oscillatory artefacts decay as $1\vert\vert \R \vert\vert^2$, which is a very slow decay.

\FloatBarrier 
%
%%%%%%%%%%%%%%%%%%%%%%%%%%%%%%%%%%%%%%%%%%%%%%%%%%%%%%%%%%%%%%%%%%%%%%%%%%%%%%%%%%%%%%%%%%%%%
%%%%%%%%%%%%%%%%%%%%%%%%%%%%%%%%%%%%%%%%%%%%%%%%%%%%%%%%%%%%%%%%%%%%%%%%%%%%%%%%%%%%%%%%%%%%%
%
\bibliography{mybibfile}

\end{document}